\documentclass[bj]{imsart}

\RequirePackage{amsthm,amsmath,amsfonts,amssymb,subcaption,dsfont,color,mathtools}
\RequirePackage[numbers]{natbib}
\RequirePackage[colorlinks,citecolor=blue,urlcolor=blue]{hyperref}
\RequirePackage{graphicx}

\startlocaldefs
\theoremstyle{plain}

\newtheorem{theorem}{Theorem}[section]
\newtheorem{lemma}[theorem]{Lemma}
\newtheorem{definition}[theorem]{Definition}
\newtheorem{corollary}[theorem]{Corollary}
\newtheorem{proposition}[theorem]{Proposition}
\newtheorem{remark}[theorem]{Remark}
\newtheorem{assumption}[theorem]{Assumption}

\newcommand{\mP}{\mathbb{P}}
\newcommand{\mE}{\mathbb{E}}

\newcommand{\mV}{\mathrm{Var}}

\newcommand{\tr}{\text{tr}}

\newcommand{\convP}{\overset{p}{\longrightarrow}}
\newcommand{\convD}{\overset{d}{\longrightarrow}}

\DeclarePairedDelimiter\floor{\lfloor}{\rfloor}

\endlocaldefs

\begin{document}

\begin{frontmatter}
\title{Dimension-agnostic inference\\ using cross U-statistics}
%\title{A sample article title with some additional note\thanksref{t1}}
\runtitle{Dimension-agnostic inference using cross U-statistics}
%\thankstext{T1}{A sample additional note to the title.}

\begin{aug}
\author[A]{\inits{F.}\fnms{Ilmun} \snm{Kim}\ead[label=e1]{ilmun@yonsei.ac.kr}} \and
\author[B]{\inits{S.}\fnms{Aaditya} \snm{Ramdas}\ead[label=e2]{aramdas@cmu.edu}}
%\and
%\author[B]{\inits{T.}\fnms{Third} \snm{Author}\ead[label=e3,mark]{third@somewhere.com}\ead[label=u1,url]{www.foo.com}}
%%%%%%%%%%%%%%%%%%%%%%%%%%%%%%%%%%%%%%%%%%%%%%
%% Addresses                                %%
%%%%%%%%%%%%%%%%%%%%%%%%%%%%%%%%%%%%%%%%%%%%%%
\address[A]{Departments of Statistics \& Data Science, Applied Statistics, 
Yonsei University, Seoul, South Korea.
\printead{e1}}

\address[B]{Departments of Statistics \& Data Science, Machine Learning, Carnegie Mellon University, Pittsburgh, USA.
\printead{e2}}
\end{aug}

\begin{abstract}
Classical asymptotic theory for statistical inference usually involves calibrating a statistic by fixing the dimension $d$ while letting the sample size $n$ increase to infinity. Recently, much effort has been dedicated towards understanding how these methods behave in high-dimensional settings, where $d$ and $n$ both increase to infinity together. This often leads to different inference procedures, depending on the assumptions about the dimensionality, leaving the practitioner in a bind: given a dataset with 100 samples in 20 dimensions, should they calibrate by assuming $n \gg d$, or $d/n \approx 0.2$? This paper considers the goal of \emph{dimension-agnostic inference}---developing methods whose validity does not depend on any assumption on $d$ versus $n$. We introduce an approach that uses variational representations of existing test statistics along with sample splitting and self-normalization to produce a refined test statistic with a Gaussian limiting distribution, regardless of how $d$ scales with $n$. The resulting statistic can be viewed as a careful modification of degenerate U-statistics, dropping diagonal blocks and retaining off-diagonal blocks. We exemplify our technique for some classical problems including one-sample mean and covariance testing, and show that our tests have minimax rate-optimal power against appropriate local alternatives. In most settings, our cross U-statistic matches the high-dimensional power of the corresponding (degenerate) U-statistic up to a $\sqrt{2}$ factor.

\end{abstract}

\begin{keyword}
\kwd{sample splitting}
\kwd{studentization}
\kwd{degenerate U-statistics}
\kwd{high-dimensional limits}
\kwd{minimax optimality}
\end{keyword}

\end{frontmatter}

\section{Introduction}

This paper deals with a widespread and rarely challenged class of assumptions in statistical inference. To clarify, by ``statistical inference'', we mean the basic and core problems of 
\begin{itemize}
\item[(A)] testing composite hypotheses while controlling their type I error uniformly over the composite null at a given level $\alpha \in (0,1)$, at least asymptotically, and 
\item[(B)] forming $1-\alpha$ confidence intervals (CIs) that are asymptotically valid. 
\end{itemize}
One typically studies the asymptotic behavior of methods in a sequence\footnote{High-dimensional analyses proceed by considering a triangular array of problems. In the $n$-th row, we have $n$ observations from some distribution that can change with $n$. Thus all distributional quantities are indexed by $n$, in particular the dimension $d_n$ (see~\cite{chernozhukov2017central} and~\cite{donoho2022optimal}, for example). 
We initially make this dependence explicit for rigor, and later drop it for conciseness.} 
of problems indexed by $n$, where we have $n$ samples from some distribution in $d_n$ dimensions, and we want to understand the coverage of CIs or level of tests asymptotically as $n$ increases to infinity.

The class of assumptions that we will challenge --- we call it a class because it is not a single specific assumption --- concerns the term ``asymptotically'' innocently used above. One must typically pre-decide some regime of asymptotics, and specify how the dimensionality $d_n$ scales relative to $n$ in the aforementioned sequence of problems: does $d_n$ stay fixed as $n$ increases (or grow much slower than $n$), or does it grow proportionally or polynomially with $n$, or does $d_n$ grow exponentially with $n$? 

Classical inference methods assumed $d$ to be fixed as $n$ grows (and thus $d/n$ or $d^2/n$ goes to zero).  Portnoy~\cite{portnoy1984asymptotic,portnoy1985asymptotic} laid the seeds of what is today called ``high-dimensional inference'', by studying the regime where $d_n^2/n$ is large, meaning $d_n$ increases to infinity with $n$ at a quadratic rate. 

The agenda of high-dimensional inference begun by Portnoy has been hugely successful, despite it still being an active area of research and thus not close to achieving its varied goals. As only a handful of influential examples, Bean et al.~\cite{bean2013optimal}, El Karoui and Purdom~\cite{el2018can} and Sur and Cand\`es~\cite{sur2019modern} bring to light issues with popular methods like the bootstrap, regression using M-estimation, and classification using logistic regression in high-dimensional settings, often studying regimes where $d_n/n$ approaches a constant. Some call the regime where $d_n$ scales exponentially with $n$ as ``ultra-high-dimensional'' inference~\citep{hao2014interaction}. In sum, (ultra-)high-dimensional inference is here to stay, and our paper does not question its importance, but instead propose a new line of inquiry. 
We dare to ask an arguably bolder question: 

\begin{quote}
\emph{Is it possible to design inference procedures that are completely dimension-agnostic, making no assumptions about the relative scaling of $d_n$ and $n$? Can a single test or CI to be asymptotically level~$\alpha$ simultaneously in all the above regimes of fixed- (or low-), high- and ultra-high dimensionality? If yes, does the power of these methods automatically adapt to the regime, being nearly as powerful as methods tuned to specific regimes?}
\end{quote}
This paper provides a positive answer to all these questions for some classical, yet nontrivial, problems. In doing so, it identifies some general strategies (like exploiting variational characterizations of test statistics) for enabling dimension-agnostic inference that could apply to a wider class of problems. 

We see this as the first work in a longer line of investigation for what appears to be both a practically important and theoretically intriguing problem. In our limited search, it seems we are the first to ask the above question explicitly and provide some new mathematical techniques (cross U-statistics) to provide a positive answer for some problems. Let us elaborate more below.

\subsection{What does it mean for inference to be dimension-agnostic?}

One of the cornerstones of parametric statistics is the likelihood ratio test. For a composite null, when certain regularity conditions hold for Wilks' theorem~\citep{wilks1938large} to apply, it states that the log-likelihood ratio has an asymptotic chi-squared distribution (under fixed-dimensional asymptotics). 

However, test statistics like the log-likelihood ratio have a different limiting distribution in fixed- or low-dimensional regimes (sample size $n\to\infty$, dimension $d$ fixed or growing very slowly) and in high-dimensional regimes ($n,d_n \to \infty$ together at some relative rate). Depending on the regime, one must calculate different calibration thresholds for a level $\alpha$ test, and even if the test statistic is itself unaltered, the overall rejection rule is different in the two regimes \citep[cf.][]{jiang2013central,sur2019modern}. The main advantage of these sophisticated analyses is the calculation of asymptotically sharp thresholds, depending on the regime of asymptotics assumed. 
However, this leads to a difficult dilemma in practice: 

\begin{quote}
\emph{Given a fixed dataset with (say) 100 samples in 20 dimensions, should the practitioner calibrate their test/CI by assuming an asymptotic regime where $n\to\infty$ with $n \gg d$, or with $d_n/n \to 0.2$? }
\end{quote}
For a particular statistic, whether the test rejects or not may critically depend on which regime is assumed. Such tests are \emph{not} dimension-agnostic since their type I error control statements depend on specific (and fine-grained) assumptions about how $d_n$ scales to infinity with $n$. Note that such assumptions are entirely unverifiable since we typically have only one dataset at hand. But the practice of calibrating tests or building CIs, by implicitly assuming that we have a sequence of datasets of particular sizes, is extremely common in statistics. Such purely theoretical assumptions, with direct practical implications, may sometimes be questionable.

In this paper, we challenge the very necessity for \emph{any} such assumption about fixed- or high- or ultra-high-dimensional asymptotics.
We revisit a handful of well-studied testing problems and develop a test statistic that has a single limiting distribution as $n \to \infty$ regardless of how $d_n$ scales (whether it stays fixed or increases to infinity, and if it does increase to infinity then we are agnostic to the rate at which it does so). We label such a test as being ``dimension-agnostic'', and highlight that the dimension-agnostic property is about validity of a test (i.e.~type I error control). Like any other valid procedure, the power of a dimension-agnostic test will depend on the (ambient or intrinsic) dimension, see Remark~\ref{rem:dim}.

Formally, our interest is in achieving the following goal in composite null hypothesis testing: given a sequence of null hypotheses indexed by the sample size $H_{0,n}: P \in \mathcal{P}_{0,n,d_n}$ for some sequence $(d_n)_{n=1}^\infty$, we would like to construct a $p$-value $p_n$ such that 
\begin{equation}\label{eq:formal-goal}
	\limsup_{n\rightarrow\infty} \sup_{P \in \mathcal{P}_{0,n,d_n}} P (p_n \leq \alpha) \leq \alpha \text{ for any } \alpha \in [0,1], \text{ regardless of  } (d_n)_{n=1}^\infty.
\end{equation}
In words, given $n$ samples from some distribution $P$ in $d_n$ dimensions, we would like to test if $P$ lies in some composite class $\mathcal{P}_{0,n,d_n}$ or not, while \emph{uniformly} controlling the type I error.

If our goal is estimation, then we desire a confidence set for some functional  $\theta_n \equiv \theta_n(P)$ of the underlying distribution $P$, whose validity does not depend on how $d_n$ scales. In particular, given a class of distributions $\mathcal{P}_{\theta_n,n,d_n}$ associated with $\theta_n$ (meaning that there could be many distributions $P$ with the same functional $\theta_n$), we call $\mathcal{C}_n$ a dimension-agnostic confidence set for $\theta_n$ if
\begin{align} \label{eq:formal-goal2}
	\liminf_{n\rightarrow \infty} \inf_{P \in \mathcal{P}_{\theta_n,n,d_n}} P(\theta_n \in \mathcal{C}_n) \geq 1 - \alpha, \ \text{regardless of } (d_n)_{n=1}^\infty.
\end{align}
For conciseness of notation, we often drop the dependence on $n$ and refer to $\mathcal{P}_{\theta_n,n,d_n}$ as $\mathcal{P}_{\theta,d}$.

It is our experience that the statistical inference literature is rife with examples of tests or confidence sets that are \emph{not} dimension-agnostic, though we mention some exceptions later in the introduction.

\begin{remark}[Relationship to distribution-free inference] \normalfont We note briefly that ``distribution-free'' inference is not the same as dimension-agnostic inference, though the terms are related. Sometimes, the (asymptotically) distribution-free property is only established under specific dimensionality regimes. Indeed, several \emph{univariate} rank tests are famously distribution-free~\citep{lehmann1975nonparametrics} so the former term is already used to describe methods that applies when $d=1$ and $n \to \infty$. 
\end{remark}

\begin{remark}[Ambient vs.\ intrinsic dimension] \label{rem:dim}
\normalfont
	Strictly speaking, the asymptotic behavior of  test statistics depends on the intrinsic dimension, henceforth denoted by $D$, rather than the ambient dimension $d$ of the data. In some settings, these two dimensions coincide. In order to maintain simplicity (and also as is common in the literature), we make no distinction between these two notions when handling finite-dimensional parameters (e.g.~Section~\ref{Section: Mean Testing} and Section~\ref{Section: One-sample covariance testing}). However, when dealing with objects in an infinite-dimensional Hilbert space (e.g.~Appendix~\ref{Section: One-sample MMD and minimax optimality}), this distinction is crucial. We make it clear what we mean by the intrinsic dimension in our context and discuss its role for classical U-statistics in Remark~\ref{Remark: intrinsic dimension}. Importantly, the validity of our proposed procedures is not affected by either the ambient dimension or the intrinsic dimension. That is, we design an inference procedure whose type I error guarantees are independent of both  $D$ and $d$. As mentioned earlier, their power will depend on $D$.
\end{remark}

This paper chooses to highlight the dimension-agnostic property, and cast it as a worthy central goal, as opposed to a side benefit. We describe a new technique to derive dimension-agnostic tests that applies to a variety of problems that involve test statistics that are degenerate U-statistics under the null, or employ empirical estimates of integral probability metrics. The main idea is simple, yet powerful, and hence somewhat broadly (but not universally) applicable. To exemplify the power of this technique, we study three classical parametric and nonparametric testing problems in some detail. 

Our methods have several desirable properties worth highlighting in advance of presenting them:

\begin{itemize}
	\item[(A)] \textbf{Weaker assumptions:} Our dimension-agnostic asymptotics hold under weaker moment assumptions than past works that study the same problems as we do. This is due to a combination of techniques from sample splitting and studentization that we elaborate on later. In addition to the obvious practical benefit (lighter assumptions are more likely to hold in practice), there is also a mathematical benefit of lesser assumptions: cleaner, streamlined proofs. This latter benefit holds by design (our refined constructions of the test statistic), and not by coincidence, and these shorter (though nontrivial) proofs are to be seen as a benefit and not a drawback of this work. \vskip .5em
	
	\item[(B)] \textbf{Explicit, simple null distribution:} 
	For many classical problems, standard test statistics have chi-squared (or mixtures thereof) null distributions in low dimensions, but Gaussian null distributions in high dimensions~\citep[cf.][]{jiang2013central}. For other problems, the test statistic is  a degenerate U-statistic whose null distribution (fixed $d$) is very complicated: it is an infinite weighted sum of chi-squared distributions, where the weights depend on the unknown underlying distribution~\citep{serfling2009approximation}. However, our new sample split variant for every problem in the paper, including the one sample kernel maximum mean discrepancy (MMD), has a standard Gaussian asymptotic null distribution as $n\to\infty$ regardless of the assumption on the (intrinsic) dimension. 
	This has a computational benefit, because it avoids the need for permutation/bootstrap based inference. \vskip .5em
	
	\item[(C)] \textbf{Fast and uniform convergence to the null:} High dimensional convergence results in the existing nonparametric testing literature often have a slow rate of convergence to the Gaussian null, with the rate depending on $d_n$ (if such a rate has been identified in the first place). In contrast, we establish Berry--Esseen bounds which demonstrate that our sample-split variants have a fast $1/\sqrt{n}$ convergence rate to the Gaussian null, independent of $d_n$, which may come as a surprise. Importantly, our convergence holds \emph{uniformly} over the composite null set, unlike several existing results in the high-dimensional literature which only prove pointwise convergence \citep[cf.][]{srivastava2008test,chen2010two}. Together with the weaker assumptions, our tests instill even greater confidence in type I error control at practical sample sizes and dimensions.  \vskip .5em	
	
	\item[(D)] \textbf{Minimax rate-optimal power:} 
	Last, our new sample-split tests have good power in both theory and practice. We prove that they are minimax rate-optimal in all the problems we study. In addition, we show that in some high-dimensional settings, their power is identical in all problem dependent parameters ($n,d,\alpha$ and appropriate signal-to-noise ratios)  compared to  dedicated high-dimensional methods that do not have the dimension-agnostic property. Further, our asymptotic relative efficiency against appropriate local alternatives is only worse by a $\sqrt{2}$ factor (arising from sample splitting) compared to those methods. We think that this is a reasonably small price to pay (indeed, it could hardly be smaller) for our stronger type I error guarantee under our weaker assumptions. Nevertheless, we do not know if this $\sqrt 2$ factor is a necessary price for dimension-agnostic inference to be possible --- we tried, and failed, to eliminate it --- and we leave this as a question for others to pursue as followup work.
\end{itemize}

%To summarize, we present a technique to construct and analyze new test statistics with many advantages for a handful of classical and modern testing problems. Our statistics are all new, sample splitting variants of well-known, often standard, test statistics

To summarize, we present a novel approach for constructing and analyzing test statistics based on sample splitting with many advantages for a handful of classical and modern testing problems. With one exception of the test statistic for mean testing that bears similarities to those in \cite{huang2015projection} and \cite{liu2022projection}, our proposed statistics are all new and regarded as sample splitting variants of well-known, often standard, test statistics. %Although our test statistic for mean testing bears similarities to those in \cite{huang2015projection} and \cite{liu2022projection}, it is important to note that our analysis is original and provides rigorous theoretical guarantees. 
In each of these examples, we provide several new and strong technical results, such as dimension-indepedent Berry--Esseen bounds for studentized statistics. We only delve into a few examples in depth, but by the end of this paper, we expect it to be clear to the reader that these advantages extend well beyond the settings studied in depth here. But before jumping directly into the results, we provide some intuition that guided us from the start, that may aid the development of other dimension-agnostic methods.

\subsection{Some guiding intuition for developing dimension-agnostic tests} \label{Section: Some guiding intuition for developing dimension-agnostic tests}
We first describe the general idea used in all the examples in this paper: \emph{sample splitting as a tool for dimensionality reduction}. We develop the idea based on integral probability metrics but it can be further generalized to any quantity that has a variational representation. Suppose that $\mathcal{F}$ is a class of real-valued bounded measurable functions on a measurable space $\mathbb{S}$. Then the integral probability metric between two probability measures $P$ and $Q$ on $\mathbb{S}$ is defined as
\begin{align*}
	\gamma_{\mathcal{F}}(P,Q):= \sup_{g \in \mathcal{F}} \big| \mE_P[g(X)] - \mE_Q[g(Y)] \big|,
\end{align*}
where $\mE_{P}[\cdot]$ (and $\mE_{Q}[\cdot]$) denotes the expectation with respect to $P$ (and $Q$). Given samples from $P$ (and also from $Q$ if it is unknown), suppose we want to make an inference based on an empirical estimate of $\gamma_{\mathcal{F}}(P,Q)$. Our strategy can be summarized in two steps. 
\begin{enumerate}
	\item In the first step, we calculate an estimated optimal function $\widehat{g}^\ast$ using the first half of the data and return $\widehat{g}^\ast$ to the second step.
	\item Then, we estimate $\mE_P[\widehat{g}^\ast(X)] - \mE_Q[\widehat{g}^\ast(Y)]$ using the second half of the data, for example by taking the sample mean, using it as a final estimator of $\gamma_{\mathcal{F}}(P,Q)$.
\end{enumerate}
This second step can be viewed as dimensionality reduction: we  project the second half of the data onto the direction $\widehat{g}^\ast$, to get a set of one-dimensional points. If $\widehat{g}^\ast$ was a ``good'' direction to project onto, then the sample mean of these univariate points will be large. But under the null where $\gamma_{\mathcal{F}}(P,Q)=0$, we effectively picked a random direction to project onto, and the univariate sample mean will simply have a (conditional) Gaussian null distribution regardless of the dimension of the original space. 

\vskip .8em

\noindent \textbf{Computationally tractable classes.} The above strategy is general but may be difficult to implement depending on the class of functions $\mathcal{F}$. To demonstrate our idea further, let us 
restrict our attention to specific $\mathcal{F}$ for which explicit estimation of $\gamma_{\mathcal{F}}(P,Q)$ is possible. To this end, let $\mathcal{H}$ be a reproducing kernel Hilbert space with its reproducing kernel $h(x,y)$. When $\mathcal{F} = \{g : \|g\|_{\mathcal{H}} \leq 1 \}$, $\gamma_{\mathcal{F}}(P,Q)$ becomes the MMD between $P$ and $Q$ associated with $h(x,y)$. Within this function class, the optimal $g^\ast$ is called the ``witness function''~\citep{gretton2012kernel} and is proportional to 
\begin{align*}
	g^\ast(x) \propto \mE_P[h(X,x)] - \mE_Q[h(Y,x)].
\end{align*}
Due to this simple representation, $g^\ast$ can be easily estimated using the sample mean. The normalizing constant is not of interest since our final test statistic will be scale-invariant. %Throughout this paper, we focus on the one-sample setting where we only have a sample from $P$ and $\mE_Q[h(Y,x)]$ is known, and discuss an extension to two-sample cases in Section~\ref{Section: Discussion}. 
Specific examples of $h(x,y)$ that we focus on to study are as follows. 
\begin{itemize}
	\item \textbf{Bi-linear kernel (Section~\ref{Section: Mean Testing}).} In the first example we choose $h(x,y) = x^\top y$ and let $P$ and $Q$ be distributions in $\mathbb{R}^d$ with $\mE_P[X] =\mu$ and $\mE_Q[Y]=0$. Then the population MMD exactly equals the $L_2$ norm of $\mu$, denoted by $\|\mu\|_2$. In this case, $g^\ast(x)$ is proportional to $\mE_P[h(X,x)] = \mu^\top x$. Following our two-step strategy described above, our statistic for this problem becomes $\widehat{\mu}_1^\top \widehat{\mu}_2$ where $\widehat{\mu}_1$ and $\widehat{\mu}_2$ are the sample means based on the first half and second half of the data, respectively. From a different perspective, $\|\mu\|_2$ or more generally  for  $p,q\in[1,\infty]$, the $L_p$ norm of $\mu$, denoted by $\|\mu\|_p$, can be obtained by solving the following optimization problem: 
	\begin{equation}\label{eq:lp-duality}
		\|\mu\|_p = \sup_{v \in \mathbb{R}^d: \|v\|_q \leq 1} v^\top \mu \text{ where } 1/p+1/q=1 .
	\end{equation}
	In this case, the estimator of $\|\mu\|_p$ based on the two-step procedure is $\widehat{v}^\top \widehat{\mu}_2$ where $\widehat{v}$ is the solution of $\sup_{v \in \mathbb{R}^d: \|v\|_q \leq 1} v^\top \widehat{\mu}_1$. When $p=q=2$, we again have the estimator $\widehat{\mu}_1^\top \widehat{\mu}_2$ as before.  \vskip .5em
	\item \textbf{Covariance kernel (Section~\ref{Section: One-sample covariance testing}).} Next we choose $h(x,y) = \tr\{ (xx^\top - I) (yy^\top - I) \}$ where $I$ is the identity matrix. Further let $P$ and $Q$ be distributions in $\mathbb{R}^d$ with $\mE_P[XX^\top] = \Sigma$, $\mE_Q[YY^\top] = I$ and $\mE_P[X] = \mE_Q[Y] = 0$. In this case the population MMD becomes the Frobenius norm of $\Sigma - I$, denoted by $\|\Sigma - I\|_F$, and $g^\ast(x)$ is proportional to $\mE_P[h(X,x)] = \tr\{ (\Sigma - I)(xx^\top - I) \}$. Following our two-step strategy, the statistic for this example becomes $\tr\{ (\widehat{\Sigma}_1 - I)(\widehat{\Sigma}_2 - I) \}$ where $\widehat{\Sigma}_1$ and $\widehat{\Sigma}_2$ are the sample covariances based on the first half and second half of the data, respectively.  \vskip .5em
	\item \textbf{Gaussian kernel (Appendix~\ref{Section: One-sample MMD and minimax optimality}).} In the third example, $h(x,y)$ is chosen as a Gaussian kernel defined later in (\ref{Eq: Gaussian kernel0}) and the resulting $\gamma_{\mathcal{F}}(P,Q)$ is known as Gaussian MMD. Unlike the previous examples, a closed-form expression of $g^\ast$ may not be available in this case but it can be easily estimated by taking the sample mean based on the first half of the data. Then the final estimator of $\gamma_{\mathcal{F}}(P,Q)$ is estimated using the second half of the data similarly as in the previous examples. 
\end{itemize}

\begin{figure}[t]
	\begin{center}
		\includegraphics[width=0.8\textwidth]{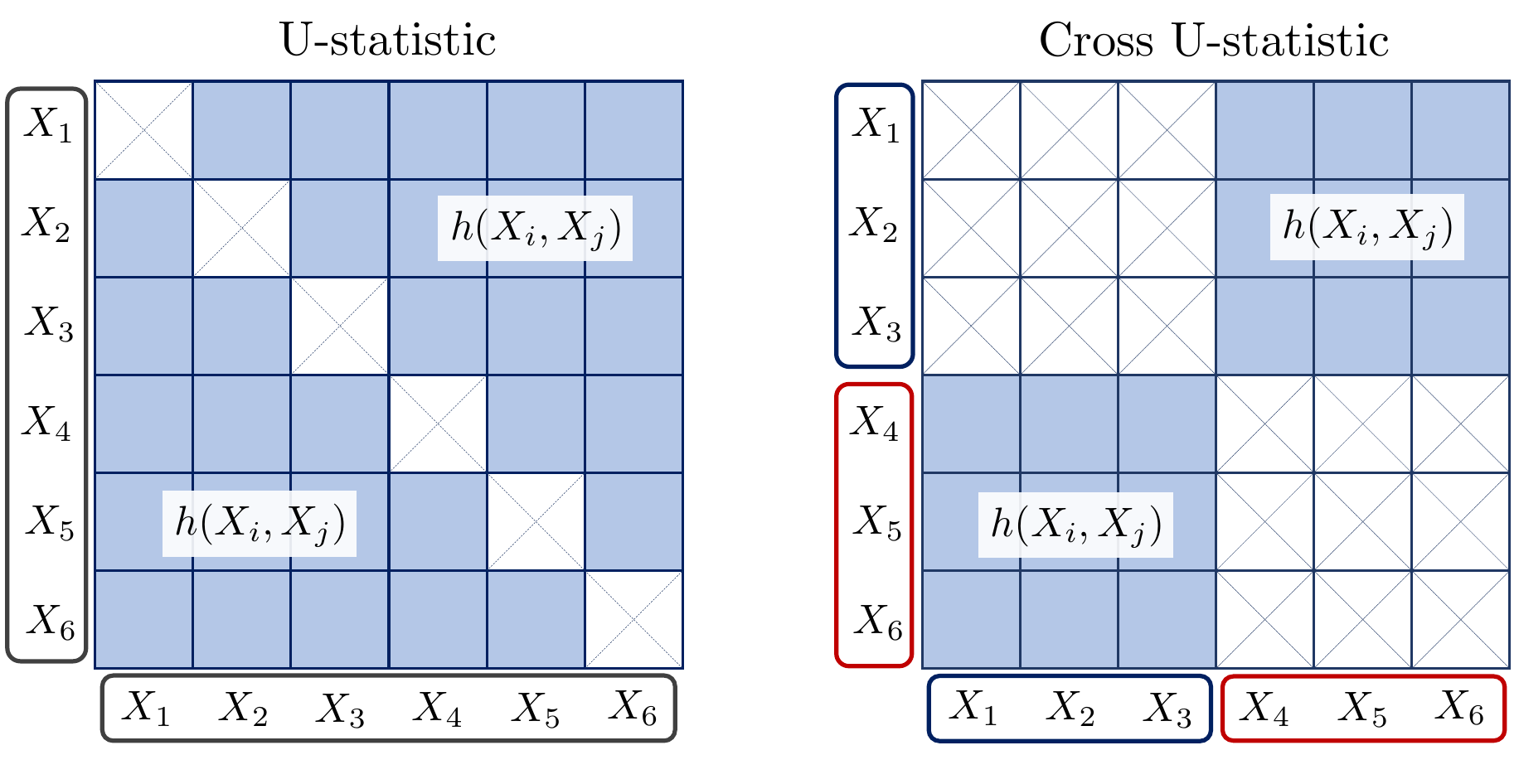}
	\end{center}
	\caption{Pictorial illustration of the difference between the U-statistic~(\ref{Eq: U-statistic}) and the proposed sample-split counterpart~(\ref{Eq: proposed statistic}) based on the same kernel $h(x,y)$. The U-statistic is defined as the average of all pairwise distances among observations, corresponding to all elements in the $6 \times 6$ kernel matrix except the diagonal components. On the other hand, the proposed statistic via sample-splitting is defined as the average of pairwise distances between observations from two disjoint subsets, corresponding to the upper-right and lower-left components in the $6 \times 6$ kernel matrix.  
 % The above visualization is useful in all examples in the paper. 
 We prove that this change yields a dimension-agnostic null, while retaining minimax rate optimal power.} 
	\label{Figure: visualization} 
\end{figure}

In all of the above examples, the final statistic can be viewed as a simple linear statistic conditional on the first half of the data. Therefore one can expect that a properly studentized statistic has a Gaussian limiting distribution. The goal of this paper is to formalize our intuition with these concrete examples but our strategy is clearly more general as described before: it would apply to any problem where the underlying quantity of interest, usually the specific notion of `signal' in the problem, has a variational representation. We would perform the optimization on one half of the data to learn an `estimator', and evaluate the quality of the estimator on the other half. In fact, this switch from in-sample to out-of-sample estimates of quality is also a central piece of technical intuition that guided our search for dimension-agnostic techniques. Informally, in the first example, the estimator $\|\widehat \mu\|^2$, where $\widehat \mu$ is the sample mean from both halves of the data, can be viewed as estimating optimal $g \in \mathcal{F}$ from the data and then evaluating its quality on the same data. The problem that arises is that the amount of `overfitting' (especially its tail behavior) is dimension-dependent, causing the null distribution to change depending on the setting. However, combined with the aforementioned intuition about dimensionality reduction, the independence of the two halves causes the quality estimates to be well-behaved out-of-sample, regardless of the dimension. 

It is worth pointing out that the proposed statistic is closely related to a U-statistic, and in fact it is an example of incomplete U-statistics~\citep{blom1976some,lee1990u}. In contrast to the (complete) U-statistic that takes into account all non-diagonal elements in a kernel matrix, our approach carefully selects \emph{half} the entries in the same kernel matrix and achieves the dimension-agnostic property, while preserving a good power property. For reasons that are apparent from the visual illustration in Figure~\ref{Figure: visualization}, we call this a \emph{cross U-statistic}, and consider it to be an object that is worthy of independent study due to its rather different limiting behavior from U-statistics, especially in the degenerate case. We lay the foundations of such a study in this work; see Section~\ref{Section: General results for degenerate U-statistics} for a formal explanation.

\subsection{Showing asymptotic normality under sample splitting} \label{Section: Showing asymptotic normality under sample splitting}

\begin{figure}[h]
	\begin{center}
		\includegraphics[width=0.85\textwidth]{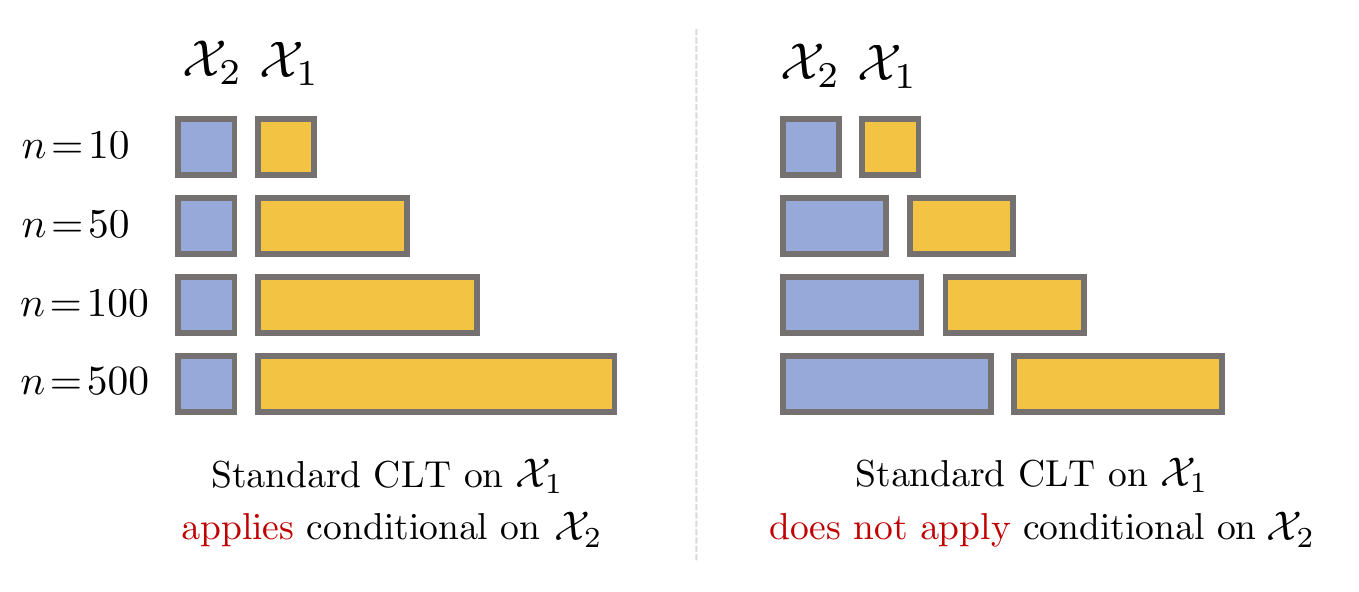}
	\end{center}
	\caption{Illustration of the data settings with a fixed conditioning set $\mathcal{X}_2$ (left) and an increasing conditioning set $\mathcal{X}_2$ (right). When $\mathcal{X}_2$ is fixed, the standard CLT on $\mathcal{X}_1$ applies conditional on $\mathcal{X}_2$. However, when $\mathcal{X}_2$ increases with $n$, proving the asymptotic normality becomes nontrivial and requires  additional technical effort as explained in Section~\ref{Section: Showing asymptotic normality under sample splitting}.} 
	\label{Figure: CLT under sample splitting}
\end{figure}

The principle outlined in Section~\ref{Section: Some guiding intuition for developing dimension-agnostic tests} elucidates the role of sample splitting in obtaining a Gaussian limiting distribution. It turns out that establishing the asymptotic normality under sample splitting is subtle and requires an extra layer of sophisticated analysis. This section briefly describes this subtlety. For concreteness, suppose that we observe i.i.d.~random objects $X_1,\ldots,X_n$ of size $n$, and set $\mathcal{X}_1 := \{X_i\}_{i=1}^{m_1}$ and $\mathcal{X}_2 := \{X_i\}_{i=m_1+1}^{n}$. The underlying principle of our proposal is to learn an optimal function $g^\ast$ based on $\mathcal{X}_2$, and project the data in $\mathcal{X}_1$ onto the estimated direction $\widehat{g}^\ast$. As emphasized earlier, the projected univariate variables $\widehat{g}^\ast(X_1),\ldots,\widehat{g}^\ast(X_{m_1})$ are i.i.d.~conditional on $\mathcal{X}_2$, and hence the corresponding studentized statistic may approximate a Gaussian distribution. This is indeed ensured by the standard central limit theorem (CLT) when $\mathcal{X}_2$ is fixed (as in the left plot of Figure~\ref{Figure: CLT under sample splitting}). However, once we allow the size of $\mathcal{X}_2$ to increase with $n$ (as in the right plot of Figure~\ref{Figure: CLT under sample splitting}), $\widehat{g}^\ast(X_1),\ldots,\widehat{g}^\ast(X_{m_1})$ should be treated as a triangular array conditional on $\mathcal{X}_2$, and thus the standard CLT --- assuming a fixed sequence of i.i.d.~random variables --- may fail (see Appendix~\ref{Section: Asymptotic analysis with an increasing conditioning set} for a concrete example).

This subtle difference is often overlooked in the literature~\cite[e.g.,][]{huang2015projection,kubler2022witness,liu2022projection}, and one should rely on a more advanced technical argument such as the Lyapunov CLT to bypass this issue. It is also worth pointing out that the ultimate goal we would like to achieve is the asymptotic normality unconditional on $\mathcal{X}_2$. A direct use of the CLT under sample splitting only gives us a conditional statement, and their conditions are probabilistic depending on $\mathcal{X}_2$. Another technical challenge is to convert these probabilistic conditions to deterministic (and more interpretable) moment conditions under which the unconditional CLT can be verified. This highlights the nontriviality of our technical contributions and brings out the underrated challenge regarding establishing CLTs under sample splitting.

\subsection{Some related work on dimension-agnostic inference}
Our work is related to but conceptually more general than Kim et al.~\cite{kim2016classification} which investigates two-sample tests based on classification accuracy. Similar to our approach, their test statistic is constructed based on sample splitting where the first half of the data is used to train a classifier---this step can be viewed as dimensionality reduction---and the second half of the data is used to estimate the accuracy of the trained classifier. They show that the resulting test statistic has a Gaussian limiting null distribution, which holds independent of the dimension. In addition, focusing on Fisher's LDA classifier, they prove that the accuracy-based test achieves an asymptotic relative efficiency (ARE) of $1/\sqrt{\pi}$ compared to Hotelling's $T^2$ test in some high-dimensional setting. Our work establishes similar results for different problems, showing that the proposed tests have the ARE of $1/\sqrt{2}$ relative to other existing tests without sample splitting which are not dimension-agnostic (Theorem~\ref{Theorem: asymptotic power} and Theorem~\ref{Theorem: asymptotic power of covariance test}). 

Another related work is that of Paindaveine and Verdebout~\cite{paindaveine2016high} which studies the dimension-agnostic property of multivariate sign tests. The authors prove that the fixed $d$-multivariate sign tests are asymptotically valid in any high-dimensional regime under some symmetry condition. They illustrate this property based on several applications including testing for uniformity and independence. Nevertheless their focus is on type I error control, and the techniques considered in their paper may not lead to the dimension-agnostic property for different types of statistics. On the other hand our focus is on both type I and type II error control (indeed, a large part of our analysis is focused on the latter) and our underlying idea is more broadly applicable than \cite{paindaveine2016high} as described in Section~\ref{Section: Some guiding intuition for developing dimension-agnostic tests}.

Tests based on center-outward ranks, recently proposed by \cite{hallin2017distribution,chernozhukov2017monge,hallin2021distribution}, have the dimension-agnostic property under a mild continuity assumption. In particular, the null distribution of the center-outward ranks is shown to be invariant to the underlying data generating distribution, which has motivated the development of distribution-free independence and two-sample tests~\citep{shi2020distribution,deb2021multivariate,deb2021efficiency}. Nevertheless, the computation of the center-outward ranks is costly in large sample sizes and their theoretical properties are largely unknown especially in high-dimensional scenarios. Indeed, as far as we are aware, the power of center-outward ranks-based tests has been studied only in a fixed-dimensional setting.

Rinaldo et al.~\cite{rinaldo2019bootstrapping} present an assumption-lean inference method for linear models based on sample splitting and bootstrapping. Again, in contrast to our work, the focus of \cite{rinaldo2019bootstrapping} is mainly on the asymptotic validity (type I error) rather than the computational and statistical efficiency (type II error).

The recent universal inference technique based on the split likelihood ratio \citep{wasserman2019universal} is non-asymptotically valid, and hence dimension-agnostic. However, these authors had a different motivation: the central goal of their work was to remove all \emph{regularity conditions} required to employ the usual asymptotics of the likelihood ratio test, thus massively expanding the scope of likelihood based methods to many irregular problems for which no known inferential technique was previously known. Unlike the majority of our current paper, their method applies to situations where a likelihood can be calculated and maximized, which were mostly parametric settings (like mixtures) and some nonparametric settings where densities still exist (like testing for monotonicity or log-concavity). We also mention that there is currently limited (but promising) understanding of its power~\citep{dunn2021gaussian}. In contrast, we handle nonparametric settings and prove that our methods are minimax rate optimal in power.

More broadly, any work that applies permutation/randomization tests~\citep[e.g.~Chapter 15.2 of][]{lehmann2006testing} is also immediately non-asymptotically valid, and thus dimension-agnostic. However, the class of problems for which non-asymptotic inference is possible is large but limited. For example, permutation tests require a group structure on the null, while universal inference requires well-defined likelihood ratios, both of which we do not need in our work. Further, our understanding the power of permutation methods is still rather incomplete~\citep{kim2020minimax}. The current paper derives dimension-agnostic tests for problems where ensuring finite-sample validity is extremely hard without involving stringent assumptions, and establishes minimax rate optimality even in high-dimensional settings.

We would not be surprised if there are also other works with dimension-agnostic tests. Nevertheless, it appears that our work is conceptually new in defining this as an explicit goal, and proposing a new, general, powerful and somewhat broadly applicable technique.

\begin{remark}[Dimension-agnostic inference is easy at the expense of power] \normalfont
In some applications, there is a simple way of achieving the goal of dimension-agnostic inference, but at the expensive of power. In the context of the MMD mentioned earlier, for example, one can consider the sample mean of the MMD statistics (or kernels) computed over disjoint blocks of data. This type of statistics is often called block-wise or linear-time statistics \citep[Lemma 14]{gretton2012kernel} --- referring back to Figure~\ref{Figure: visualization}, the block-wise method corresponds to using the elements only within small blocks along the diagonal, while ignoring off-diagonal elements outside a small band. 

This is in sharp contrast to our method. Since the constructed statistic is a sum of i.i.d.~random variables, it clearly approximates a Gaussian distribution irrespective of the dimension. However, as pointed out by \cite{ramdas2015adaptivity} and \cite{kim2018robust}, the power of the block-wise test is typically worse by more than a constant factor compared to the corresponding V- or U-statistic; it has a worse rate and is not minimax optimal. 

Another way of obtaining a Gaussian limit for the MMD statistic has been suggested by \cite{makigusa2020asymptotic}. Building on the ideas of Ahmad~\cite{ahmad1993modification}, they propose a modified V-statistic that is asymptotically Gaussian under the null. Nevertheless, we expect that the power of their approach is not rate optimal, given that the convergence rate of the proposed statistic is $1/\sqrt{n}$, which is much slower than the $1/n$ rate of the corresponding V-statistic under the null. Our approach, on the other hand, yields a rate optimal test, while maintaining the dimension-agnostic property. 
\end{remark}

\subsection{Paper outline and technical notation} \label{Section: Notation}

The rest of this paper is organized as follows. Section~\ref{Section: Mean Testing} focuses on one-sample mean testing and formally develops our intuition. In particular we show that the proposed method for one-sample mean testing has the dimension-agnostic property, which is in contrast to the existing method using a U-statistic. We also illustrate that the proposed method possesses good power properties against dense or sparse alternatives. We end Section~\ref{Section: Mean Testing} by describing dimension-agnostic confidence sets for a mean vector. Section~\ref{Section: One-sample covariance testing} provides similar results in the context of one-sample covariance testing. In Section~\ref{Section: General results for degenerate U-statistics}, we identify general conditions under which a sample-splitting analogue of a degenerate U-statistic has a Gaussian limiting distribution. We end with a discussion and directions for future work in Section~\ref{Section: Discussion}. Additional results are in the appendices. Appendix~\ref{Section: Multiple sample-splitting} discusses the approach based on multiple sample-splitting, while Appendix~\ref{Section: Further discussions on power} presents a general strategy for studying the asymptotic power of the proposed method. Appendix~\ref{Section: Asymptotic analysis with an increasing conditioning set} provides an example that demonstrates the non-triviality of the asymptotic normality under sample splitting. In Appendix~\ref{Section: Dimension-agnostic confidence set}, we construct dimension-agnostic confidence sets for mean vectors by inverting dimension-agnostic $p$-values. Appendix~\ref{Section: One-sample MMD and minimax optimality} illustrates our main results using Gaussian MMD and studies minimax power against nonparametric alternatives. We support our theoretical findings with simulations in Appendix~\ref{Section: Simulations}, and all proofs are provided in Appendix~\ref{Section: Proofs}.

Throughout this paper, we adopt the following notation. Let $\Phi(\cdot)$ be the cumulative distribution function of the standard normal random variable $N(0,1)$ and $z_\alpha$ be the $\alpha$ quantile of $N(0,1)$. The symbol $\mathds{1}[\cdot]$ is used for indicator functions and we write $P[\cdot]=\mE_{P} \mathds{1}[\cdot]$. For any real sequences $\{a_n\}$ and $\{b_n\}$, we write $a_n = O(b_n)$ if there exists a constant $C$ such that $|a_n| \leq C|b_n|$ for sufficiently large $n$ and $a_n = o(b_n)$ if for any $c>0$, $|a_n| \leq c|b_n|$ holds for sufficiently large $n$. The symbol $a_n \asymp b_n$ means that $a_n = O(b_n)$ and $b_n = O(a_n)$. We say $\{X_1,\ldots,X_n\} \overset{i.i.d.}{\sim} P$ if random variables $\{X_1,\ldots,X_n\}$ are independent and identically distributed with the common distribution $P$. For a random vector $X \in \mathbb{R}^d$, we denote the $i$th component of $X$ by $X^{(i)}$ for $i=1,\ldots,d$.

\section{One-sample mean testing} \label{Section: Mean Testing}
We start with the simple problem of one-sample mean testing for which we can concretely deliver our intuition. Suppose that we observe $d$-dimensional i.i.d.~random vectors $X_1,\ldots,X_n$ from a distribution with mean $\mu$ and covariance matrix $\Sigma$. Given this sample, we would like to test whether
\begin{align} \label{Eq: one-sample hypothesis}
	H_0: \mu = 0 \quad \text{versus} \quad H_1:  \mu \neq 0.
\end{align}
This classical problem has been well-studied both in low- and high-dimensional scenarios. In the fixed-dimensional setting, Hotelling's $T^2$ test is one of the most well-known testing procedures with several attractive properties \citep{anderson2003introduction}. However, as highlighted in \cite{bai1996effect}, Hotelling's $T^2$ test performs poorly or is not even applicable when the dimension is comparable to or exceeds the sample size. To address this issue, several tests, which are specifically designed for high-dimensional data, have been proposed in the literature \citep[e.g.][for a review]{hu2016review}. For instance, \cite{chen2010two} introduce a test based on the following U-statistic:
\begin{align} \label{Eq: One-sample Mean U-statistic}
	U_{\mathrm{mean}} := \frac{1}{n(n-1)} \sum_{1 \leq i \neq j \leq n} X_i^\top X_j.
\end{align}
Under a pseudo independence model assumption (see e.g.~Appendix~\ref{Section: Details on the moment condition}), Chen and Qin~\cite{chen2010two} show that $U_{\mathrm{mean}}$ converges to a Gaussian distribution when $d$ increases to infinity with $n$. However, when $d$ is fixed, there is no such guarantee. Indeed classical asymptotic theory for U-statistics reveals that $U_{\mathrm{mean}}$ approximates a weighted sum of chi-square random variables. Moreover, even in high-dimensional settings, the limiting null distribution of $U_{\mathrm{mean}}$ can be far from Gaussian depending on the covariance structure of $X$ \citep{wang2019feasible}. This issue has been partially addressed by \cite{wang2019feasible} where they prove that a bootstrap procedure can consistently estimate the null distribution of $U_{\mathrm{mean}}$ in different asymptotic scenarios. However such asymptotic guarantee hinges on the same model assumption considered in \cite{bai1996effect} and \cite{chen2010two}. In addition the computational cost of bootstrapping is typically prohibitive for large data analysis.

\vskip .8em

\noindent \textbf{A heuristic explanation of why $U_{\mathrm{mean}}$ has different behavior.} Before presenting our approach that alleviates the aforementioned issues, we briefly elaborate on why the asymptotic behavior of $U_{\mathrm{mean}}$ may differ in different regimes. For the sake of simplicity, assume that $X$ is centered and has a diagonal covariance matrix with $d$ elements $\{\sigma_1^2,\ldots,\sigma_d^2 \}$. Then, for a fixed $d$, the asymptotic theory of degenerate U-statistics \citep{lee1990u} shows $nU_{\mathrm{mean}} \approx \sum_{k=1}^d \sigma_k^2 (\xi_k^2 - 1)$, where $\{\xi_i\}_{i=1}^d \overset{i.i.d.}{\sim} N(0,1)$. In other words, $U_{\mathrm{mean}}$ behaves like a sum of independent random variables in large samples. This observation together with the central limit theorem suggests that the distribution of $U_{\mathrm{mean}}$ can be further approximated by a normal distribution as $d \rightarrow \infty$ when $\sigma_k$ are properly bounded (so that each variable does not contribute too much). Indeed this boundedness condition is crucial for the normal approximation. In the extreme case where only one of $\{\sigma_1,\ldots,\sigma_d\}$, let's say $\sigma_1$, is positive, it is clear that the asymptotic behavior of $U_{\mathrm{mean}}$ is dominated by $\sigma_1^2 (\xi_1^2 - 1)$ and the resulting distribution is never close to Gaussian even if $d \rightarrow \infty$. In general the limiting distribution of $U_{\mathrm{mean}}$ is essentially determined by the (unknown) covariance structure of $X$, which is hard to estimate, especially in high dimensional settings. Of course, our explanation here is informal but we hope that it conveys our message that a valid inference based on $U_{\mathrm{mean}}$ is quite challenging in high-dimensions.

\vskip 1em

\subsection{A test statistic by sample-splitting and studentization}
Having described the limitation of the previous approach, we now introduce our test statistic for mean testing and study its asymptotic null distribution. Given positive integers $m_1$ and $m_2 = n-m_1$, we first consider a bilinear statistic defined as
\begin{align*}
	T_{\mathrm{mean}}^{\dagger} := \left(\frac{1}{m_1} \sum_{i=1}^{m_1} X_i \right)^\top \left( \frac{1}{m_2} \sum_{j=m_1+1}^{n} X_{j} \right).
\end{align*}
The definition of $T_{\mathrm{mean}}^{\dagger}$ can be motivated from the two-step approach described in Section~\ref{Section: Some guiding intuition for developing dimension-agnostic tests}. Due to the independence between $\mathcal{X}_{1} := \{X_i\}_{i=1}^{m_1}$ and $\mathcal{X}_{2} := \{X_i\}_{i=m_1+1}^n$, it is clear that the bilinear statistic $T_{\mathrm{mean}}^{\dagger}$ is an unbiased estimator of $\|\mu\|_2^2$ as $U_{\mathrm{mean}}$. However, in sharp contrast to $U_{\mathrm{mean}}$, we claim that the limiting null distribution of a studentized $T_{\mathrm{mean}}^{\dagger}$ is Gaussian regardless of whether we are in a fixed- or a high-dimensional regime. We start with a high-level idea of why it is the case and then present a formal explanation. Let us define a random function $f_{\mathrm{mean}}(\cdot)$ depending on $\mathcal{X}_2$ as
\begin{align} \label{Eq: univariation function 1}
	f_{\mathrm{mean}}(x) := \frac{1}{m_2} \sum_{j=m_1+1}^{n} X_{j}^\top x.
\end{align}
With this notation, $T_{\mathrm{mean}}^{\dagger}$ can be viewed as the sample mean of $f_{\mathrm{mean}}(X_1), \ldots, f_{\mathrm{mean}}(X_{m_1})$, that is $T_{\mathrm{mean}}^{\dagger} =  m_1^{-1} \sum_{i=1}^{m_1} f_{\mathrm{mean}}(X_i)$. Since $f_{\mathrm{mean}}(X_1), \ldots, f_{\mathrm{mean}}(X_{m_1})$ are independent univariate random variables conditional on $\mathcal{X}_2$, one might naturally expect that $T_{\mathrm{mean}}^{\dagger}$ has a Gaussian limiting distribution in both fixed- and high-dimensional regimes.  However we cannot directly apply a standard central limit theorem since the random function $f_{\mathrm{mean}}(x)$ converges to zero almost surely as $m_2 \rightarrow \infty$ (when $\mu=0$), and thus the variance of $f_{\mathrm{mean}}(X)$ can shrink to zero. Therefore it requires some effort to prove this statement rigorously.

\subsection{Asymptotic, dimension-agnostic, Gaussian null distribution}
To this end, let us define our studentized test statistic as
\begin{align} \label{Eq: studentized statistic for mean testing}
	T_{\mathrm{mean}} := \frac{\sqrt{m_1} \overline{f}_{\mathrm{mean}}}{\sqrt{m_1^{-1} \sum_{i=1}^{m_1} \{f_{\mathrm{mean}}(X_i) - \overline{f}_{\mathrm{mean}}\}^2 }},
\end{align}
where we write $\overline{f}_{\mathrm{mean}} = T_{\mathrm{mean}}^{\dagger}$. The next lemma provides a conditional Berry--Esseen bound for the studentized statistic $T_{\mathrm{mean}}$, which is a direct consequence of Theorem 1.1 in \cite{bentkus1996berry}. 
\begin{lemma}[Conditional pointwise Berry--Esseen bound for $T_{\mathrm{mean}}$] \label{Lemma: Conditional BE Bound}
		Suppose that we are under the null where the distribution $P$ has $\mu = 0$ and assume $0 < \mE_P[f_{\mathrm{mean}}^2(X)|\mathcal{X}_2] < \infty$ almost surely (a.s.). Then there exists an absolute constant $C>0$ such that 
	\begin{align} \label{Eq: conditional Berry--Esseen bound}
		\sup_{-\infty < z < \infty} \big| P(T_{\mathrm{mean}}  \leq z | \mathcal{X}_2) - \Phi(z) \big| \leq  C \frac{\mE_P[|f_{\mathrm{mean}}^3(X)| | \mathcal{X}_2]}{\{\mE_P[f_{\mathrm{mean}}^2(X)|\mathcal{X}_2]\}^{3/2} \sqrt{m_1}} \quad \text{a.s.,}
	\end{align}
	where we recall $\mathcal{X}_2 = \{X_i\}_{i=m_1+1}^n$.
\end{lemma}
Due to the above result, our problem is simply transformed into establishing conditions under which the right-hand side of the bound~\eqref{Eq: conditional Berry--Esseen bound} converges to zero. % in probability. Under these conditions, it is clear that $T_{\mathrm{mean}}$ has a Gaussian limiting distribution conditional on $\mathcal{X}_2$. Moreover, since convergence in probability implies convergence in mean for a uniformly integrable random sequence, the conditional convergence statement further guarantees the following unconditional convergence result as well:
%\begin{align*}
%\sup_{-\infty < z < \infty} \big| \mP(T_{\mathrm{mean}}  \leq z) - \Phi(z) \big| \rightarrow 0.
%\end{align*}
For this purpose we consider the following multivariate moment condition.
\begin{assumption} [Lyapunov ratio]\label{assumption: moment condition} \normalfont
	For a given constant $L>1$, assume that $X$ is a random vector having a continuous distribution $P$ in $\mathbb{R}^d$ with mean $\mu$ and satisfying
	\begin{align*}
		\sup_{u \in \mathbb{S}^{d-1}}\frac{\big\{\mE_P \big[|u^\top (X-\mu)|^3 \big] \big\}^{1/3}}{ \big\{ \mE_P \big[ |u^\top (X-\mu)|^2 \big] \big\}^{1/2}} \leq L^{1/3},
	\end{align*}
	 where $\mathbb{S}^{d-1}$ denotes the $(d-1)$-dimensional unit sphere.
\end{assumption}
In our opinion, this moment condition is mild and is satisfied under common multivariate tail assumptions such as multivariate sub-Gaussianity or sub-Exponentiality \citep[e.g.][]{vershynin2018high}. In particular, one can take $L = 2\sqrt{2}/\sqrt{\pi}$ for a Gaussian random vector with an arbitrary positive definite covariance matrix. It is worth pointing out that our moment condition is weaker than the moment condition considered in \cite{bai1996effect,chen2010two,wang2019feasible}; see Appendix~\ref{Section: Details on the moment condition} for details. %We also note that a similar moment condition is imposed by \cite{giessing2020bootstrapping} in the context of studying bootstrap consistency for $\ell_p$-statistics in high-dimensions. 
The continuity of $X$ prevents the case where $f_{\mathrm{mean}}(\cdot)$ becomes zero with a nontrivial probability in finite samples and it is in fact sufficient to assume that one of the components of $X$ has a continuous distribution. 

%In addition, since our test statistic is defined in terms of $X_i^\top X_j$, 

Given Assumption~\ref{assumption: moment condition}, let us denote 
\[
\widehat{\mu}_2:=\frac1{m_2}\sum_{j=m_1+1}^{n}X_j,
\] 
so that we have $f_{\mathrm{mean}}(x) = \widehat{\mu}_2^\top x$. Then the right-hand side of the bound~(\ref{Eq: conditional Berry--Esseen bound}) satisfies
\begin{align*}
	\frac{\mE_P[|f_{\mathrm{mean}}^3(X)| | \mathcal{X}_2]}{\{\mE_P[f_{\mathrm{mean}}^2(X)|\mathcal{X}_2]\}^{3/2} \sqrt{m_1}} = \frac{\mE_P[|\widehat{\mu}_2^\top X|^3 | \mathcal{X}_2]}{\{\mE_P[(\widehat{\mu}_2^\top X)^2|\mathcal{X}_2]\}^{3/2} \sqrt{m_1}} \leq \frac{L}{\sqrt{m_1}} \quad \text{a.s.}
\end{align*}
This observation together with Lemma~\ref{Lemma: Conditional BE Bound} provides the unconditional normality of $T_{\mathrm{mean}}$.
\begin{theorem}[Unconditional uniform Berry--Esseen bound for $T_{\mathrm{mean}}$] \label{Theorem: unconditional BE Bound}
	Let $\mathcal{P}_{0,d}(L)$ be the class of distributions that satisfies Assumption~\ref{assumption: moment condition} with mean $\mu =0$ and some constant $L>1$. Then there exists an absolute constant $C > 0$ such that 
	\begin{align} \label{Eq: Berry--Esseen bound}
		\sup_{d \geq 1}\sup_{P \in \mathcal{P}_{0,d}(L)} \sup_{-\infty < z < \infty} \big| P(T_{\mathrm{mean}}  \leq z) - \Phi(z) \big| \leq \frac{C L}{\sqrt{m_1}}.
	\end{align}
\end{theorem}
Several remarks on Theorem~\ref{Theorem: unconditional BE Bound} are in order.  

\begin{remark} \leavevmode \normalfont
	\begin{itemize}
		\item \emph{Limiting distribution under the alternative.} While our focus is on deriving the limiting null distribution of $T_{\mathrm{mean}}$, the same bound holds after studentizing $\widehat{\mu}_2^\top (X_1 - \mu),\ldots, \widehat{\mu}_2^\top (X_{m_1}-\mu)$, which we use to study the asymptotic power in Theorem~\ref{Theorem: asymptotic power}. \vskip .5em
		\item \emph{Condition on $d$ and $m_2$.} It is interesting to point out that the upper bound for the above Kolmogorov--Smirnov distance solely depends on $m_1$ and $L$. Therefore, by considering $L$ as some fixed constant, $T_{\mathrm{mean}}$ has a Gaussian limiting null distribution as $m_1 \rightarrow \infty$, which is independent of $d$. That means, Theorem~\ref{Theorem: unconditional BE Bound} holds regardless of whether we are under a fixed- or a high-dimensional regime. We also note that the choice of $m_2$ does not affect our result in Theorem~\ref{Theorem: unconditional BE Bound}, i.e.~$T_{\mathrm{mean}}$ converges to $N(0,1)$ for both fixed or increasing $m_2$. However,  Section~\ref{Section: Mean Testing: Asymptotic Power Analysis} shows that the choice of $m_2$ has a significant effect on the asymptotic power of the resulting test. \vskip .5em
		\item \emph{Condition on $L$.} In principle, one need not regard $L$ as a constant, but as an increasing sequence $L \equiv L_n$ satisfying $L/\sqrt{m_1} \rightarrow 0$. This does not affect the asymptotic normality of $T_{\mathrm{mean}}$. However, we treat $L$ as some fixed constant (say $L=2$) for simplicity and write $\mathcal{P}_{0,d}(L)$ as $\mathcal{P}_{0,d}$. 
		\vskip .5em
		 %Next we provide a few remarks on Assumption~\ref{assumption: moment condition} as well as the effect of studentization.% and then move on to a power analysis. 
		\item \emph{Weakening Assumption~\ref{assumption: moment condition}.} While Assumption~\ref{assumption: moment condition} allows us to obtain a \emph{finite-sample} guarantee, the normality of $T_{\mathrm{mean}}$ can be achieved under much weaker conditions than Assumption~\ref{assumption: moment condition}. In fact all we need to show is that the bound~\eqref{Eq: conditional Berry--Esseen bound} (or more generally conditional Lindeberg's condition) converges to zero in probability. Once this result is established, it is clear that $T_{\mathrm{mean}}$ converges to $N(0,1)$ conditional on $\mathcal{X}_2$. Moreover, since convergence in probability implies convergence in mean for a uniformly integrable random sequence, the conditional convergence statement further guarantees the unconditional convergence result as well. We refer the reader to the proof of Theorem~\ref{Theorem: Gaussian Approximation} which uses a similar argument.  \vskip .5em
		\item \emph{Effect of studentization.} We note that the scaling factor used in $T_{\mathrm{mean}}$ is of critical importance for obtaining a Gaussian limiting distribution. To see this, suppose that $X_1,\ldots,X_n$ are univariate random variables with a mean of zero and a variance of $\sigma^2$. For simplicity, we consider the balanced splitting scheme with $m:= m_1 = m_2$. Then the bivariate central limit theorem combined with the continuous mapping theorem yields $m T_{\mathrm{mean}}^{\dagger} \convD \sigma^2 \xi_1 \xi_2$ where we recall $\xi_1,\xi_2 \overset{i.i.d.}{\sim} N(0,1)$. Therefore $T_{\mathrm{mean}}^{\dagger}$ without studentization has asymptotically the same distribution as the product of two independent normal random variables, instead of a single Gaussian. \vskip .5em
		\item \emph{t-distribution.} Since $T_{\mathrm{mean}}$ is essentially the $t$-statistic, it may be worth calibrating the test based on a $t$-distribution when the underlying distribution is close to Gaussian. According to the result of \cite{pinelis2015exact}, the Kolmogorov--Smirnov distance between the standard normal distribution and the $t$-distribution with $n-1$ degrees of freedom is bounded by $C/(n-1)$ where $n \geq 5$ and $C \approx 0.158.$ Therefore, by the triangle inequality, the Berry--Esseen bound~(\ref{Eq: Berry--Esseen bound}) still holds by replacing $\Phi(t)$ with the $t$-distribution with $m_1-1$ degrees of freedom. \vskip .5em
		\item \emph{$p$-value.} Theorem~\ref{Theorem: unconditional BE Bound} implies that $p_{\mathrm{mean}} := 1 - \Phi(T_{\mathrm{mean}})$ is a dimension-agnostic $p$-value in the sense of~\eqref{eq:formal-goal}, as long as $m_1 \to \infty$, no matter how $m_2$ behaves:
		\begin{align*}
		\lim_{m_1 \rightarrow \infty} \sup_{P \in \mathcal{P}_{0,d}} P(p_{\mathrm{mean}} \leq \alpha) = \alpha \text{ for any } \alpha \in [0,1], \text{ regardless of } (d_n)_{n=1}^\infty.
		\end{align*}
		Above, we defined a one-tailed $p$-value as we tend to observe a positive value of $T_{\mathrm{mean}}$ under the alternative. $p$-values for other statistics introduced later in this paper can be similarly defined.
        \item The tests can be inverted to form dimension-agnostic confidence sets for $\mu$ (Appendix~\ref{Section: Dimension-agnostic confidence set}).
	\end{itemize}
\end{remark}
Having characterized the limiting null distribution of $T_{\mathrm{mean}}$, we next turn our attention to power.

\subsection{Asymptotic power analysis} \label{Section: Mean Testing: Asymptotic Power Analysis}
%\subsection{Asymptotic power analysis: optimal up to a $\sqrt 2$ factor} \label{Section: Mean Testing: Asymptotic Power Analysis}
% In the previous subsection we show that $T_{\mathrm{mean}}$ has consistent limiting behavior in different asymptotic regimes. As a result, the asymptotic validity of the resulting test holds independent of $d$ under some mild moment condition. By contrast, the asymptotic distribution of $U_{\mathrm{mean}}$ varies a lot depending on whether $d$ is fixed or not, which raises the question on the right choice of calibration thresholds in practice. 
While it is already clear that our approach is beneficial over $U_{\mathrm{mean}}$ in terms of robust type I error control, we also claim that our test maintains a good power property. In particular we show that the asymptotic power of $T_{\mathrm{mean}}$ is only worse by a factor of $\sqrt{2}$ than that of $U_{\mathrm{mean}}$ under the setting where both test statistics converge to Gaussian. In a more general setting, Theorem~\ref{Theorem: Uniform type I and II error control} proves that our test is minimax rate optimal in terms of the $L_2$ norm, which means that the power of the proposed test cannot be improved beyond a constant factor. Before we proceed, let us describe the following assumptions that facilitate our asymptotic power analysis. 
\begin{assumption} \label{assumption: one-sample mean testing} \normalfont
	Suppose that the following assumptions are satisfied. 
	\begin{enumerate}\setlength{\itemindent}{0.07in}
		\item[(a)] \emph{Gaussianity:} $X$ has a multivariate Gaussian distribution $N(\mu,\Sigma)$. \vskip .5em
		\item[(b)] \emph{Bounded eigenvalues:} there exist constants $\lambda_\ast,\lambda^\ast >0$ s.t.\ $\lambda_\ast \leq \lambda_1(\Sigma) \leq \cdots \leq \lambda_d(\Sigma) \leq \lambda^\ast$. \vskip .5em
		\item[(c)] \emph{Local alternatives:} $\mu^\top \mu = O(\sqrt{d}/n)$. \vskip .5em
		\item[(d)] \emph{Dimension-to-sample size ratio:} $d/n \rightarrow \tau \in (0, \infty)$. \vskip .5em
		\item[(e)] \emph{Sample-splitting ratio:} $m_1/n \rightarrow \kappa \in (0,1)$. 
	\end{enumerate}
\end{assumption}
Of course, these assumptions are not required under the null, and are only needed to derive an explicit power expression under the alternative in Theorem~\ref{Theorem: asymptotic power}. The first Gaussian assumption is not critical and can be relaxed to the pseudo-independence model assumption considered in \cite{bai1996effect}. Under the second condition on eigenvalues, $U_{\mathrm{mean}}$ has a Gaussian limiting distribution and the power expression of our test can be simplified; this assumption can be weakened with more effort, but this takes a little further away from the point of our paper. The last two assumptions on $d/n$ and $m_1/n$ are required for standard technical reasons and can also be found in \cite{kim2016classification}. %The last two assumptions on $d/n$ and $m_1/n$ are required for some technical reasons and can also be found in \cite{kim2016classification}. 

We now present the asymptotic power expression of the test $\phi_{\mathrm{mean}} := \mathds{1}(T_{\mathrm{mean}}  > z_{1-\alpha})$. 
\begin{theorem}[Asymptotic power expression] \label{Theorem: asymptotic power}
	Suppose that Assumption~\ref{assumption: one-sample mean testing} is fulfilled under the alternative. Then it holds that 
	\begin{align*}
		\mE_P[\phi_{\mathrm{mean}}] ~=~  \Phi \left( z_{\alpha} + \frac{n\sqrt{\kappa(1-\kappa)}\mu^\top \mu}{\sqrt{\emph{\tr}(\Sigma^2)}} \right) + o(1).
	\end{align*}
	Therefore the power is asymptotically maximized when $m_1/n \rightarrow \kappa =1/2$. 
\end{theorem}

\noindent It is now apparent that the optimal choice of $\kappa$ equals $1/2$, and the resulting power is
\begin{align} \label{Eq: theoretical power of phi_mean}
	\mE_P[\phi_{\mathrm{mean}}] = \Phi \left( z_{\alpha} + \frac{n \mu^\top \mu}{2\sqrt{\tr(\Sigma^2)}} \right) + o(1).
\end{align}
We note in passing that the same optimal choice of $\kappa$ was observed in \cite{kim2016classification} as well. 
On the other hand, by building on work of \cite{wang2019feasible}, it can be seen that the test based on $U_{\mathrm{mean}}$ has the following asymptotic power under the same assumption in Theorem~\ref{Theorem: asymptotic power}:
\begin{align} \label{Eq: theoretical power of U_1}
	P \left( U_{\mathrm{mean}} > z_{1-\alpha} \sqrt{\frac{2}{n(n-1)} \tr(\Sigma^2)} \right) ~=~  \Phi \left(z_{\alpha}  + \frac{n \mu^\top \mu}{\sqrt{2 \tr(\Sigma^2)}} \right) + o(1).
\end{align}
Thus, the power of our test is lower by a $\sqrt{2}$ factor than the test based on $U_{\mathrm{mean}}$. This  loss of power may be understood as the price to pay for the dimension-agnostic property of $T_{\mathrm{mean}}$ under the null. It is also worth pointing out that the asymptotic expression of the minimax power in the Euclidean norm is equivalent to $\Phi(z_\alpha + n \mu^\top \mu / \sqrt{2d})$ when $\Sigma$ is the identity matrix. See Chapter 1.3 of \cite{ingster2012nonparametric} for details, and Proposition 3.1 of \cite{kim2016classification} that extends this result to the two-sample problem. Thus the asymptotic power expression~\eqref{Eq: theoretical power of phi_mean} matches the minimax power, up to a $\sqrt{2}$ factor, when $\Sigma$ is the identity matrix.

\begin{remark}[Intuition on $\sqrt{2}$ factor]	\normalfont \label{Remark: Intuition of constant factor}
	When a test statistic converges to a Gaussian distribution under the null as well as the alternative, the power of the resulting test is mainly determined by the mean and the variance of the test statistic. Moreover, in the case of local alternatives where a signal shrinks with $n$, the variance of the test statistic tends to be the same under both $H_0$ and $H_1$. With this observation, we first remark that both $U_{\mathrm{mean}}$ and $T_{\mathrm{mean}}^{\dagger}$ are unbiased estimators of $\mu^\top \mu$. On the other hand, the variances of these statistics are
	\begin{align*}
		\mV_P[U_{\mathrm{mean}}] = \frac{2}{n(n-1)}\tr(\Sigma^2) \quad \text{and} \quad \mV_P[T_{\mathrm{mean}}^{\dagger}] = \frac{1}{m_1m_2} \tr(\Sigma^2) \quad \text{under $H_0$. }
	\end{align*}
	Thus, when $m_1$ and $m_2$ are approximately the same, the variance of $T_{\mathrm{mean}}^{\dagger}$ is (approximately) twice larger than that of $U_{\mathrm{mean}}$ under $H_0$, which in turn leads to $\mE_P[T_{\mathrm{mean}}^{\dagger}] / \{\mV_P[T_{\mathrm{mean}}^{\dagger}]\}^{1/2} \approx \mE_P[U_{\mathrm{mean}}] /\{2\mV_P[U_{\mathrm{mean}}]\}^{1/2}$ against local alternatives, roughly explaining the $\sqrt{2}$ in the power~\eqref{Eq: theoretical power of phi_mean}. 
\end{remark}

Next, we relax Assumption~\ref{assumption: one-sample mean testing} and study optimality from a minimax perspective.

\subsection{Minimax rate optimality against Euclidean norm deviations} \label{Section: Minimax optimality}
In the previous subsection, we impose rather strong assumptions on the distribution of $X$ and present an asymptotic power expression, which is precise including constant factors. The aim of this subsection is to relax Assumption~\ref{assumption: one-sample mean testing} and study minimax rate optimality of the proposed test in terms of the $L_2$ distance. To formulate the minimax problem, let $\mathcal{P}_{\mu,d}$ be the set of all possible distributions with mean vector $\mu \in \mathbb{R}^d$ that satisfies Assumption~\ref{assumption: moment condition} with some fixed constant $L$ (say $L=2$). We denote by $\Psi(\alpha)$ the set of all asymptotically level $\alpha$ tests over $\mathcal{P}_{0,d}$ with $\mu=0$. More formally, we let
\begin{align*}
	\Psi(\alpha) := \bigg\{ \phi :  \limsup_{n \rightarrow \infty} \sup_{P \in \mathcal{P}_{0,d}} \mE_{P}[\phi] \leq \alpha \bigg\}.
\end{align*}
Given $\epsilon_n>0$ --- shorthand for a positive sequence $(\epsilon_n)_{n \geq 1}$ --- and $\lambda^\ast>0$, define local alternatives 
\begin{align}
	\mathcal{Q}_{\mathrm{mean},d}(\epsilon_n) := \{P \in \mathcal{P}_{\mu,d} : \| \mu \|_2 \geq \epsilon_n \ \text{and $\lambda_{\mathrm{max}}(\Sigma) < \lambda^\ast$} \}.
\end{align}
We then claim that the proposed mean test $\phi_{\mathrm{mean}}$ is asymptotically level $\alpha$ and the type II error of $\phi_{\mathrm{mean}}$ is uniformly small over $\mathcal{Q}_{\mathrm{mean},d}(\epsilon_n)$ when $\epsilon_n$ is sufficiently larger than $d^{1/4} / \sqrt{n}$. %In view of the optimal choice of $\kappa$~(\ref{Eq: optimal power expression}), the next proposition focuses on the case of $m_1 = \floor{n/2}$ and $m_2=n-m_1$. 
\begin{theorem}[Uniform type I and II error] \label{Theorem: Uniform type I and II error control}
	For fixed $\alpha \in (0,1)$, suppose that $\epsilon_n \geq K_n d^{1/4} n^{-1/2}$ where $K_n$ is a positive sequence increasing to infinity at any rate. Then
	\begin{align*}
		\lim_{n \rightarrow \infty} \sup_{P \in \mathcal{P}_{0,d}} \mE_{P} [\phi_{\mathrm{mean}}] = \alpha \quad \text{and} \quad 	\lim_{n \rightarrow \infty} \sup_{P \in \mathcal{Q}_{\mathrm{mean},d}(\epsilon_n)} \mE_{P} [1 -  \phi_{\mathrm{mean}}] = 0,
	\end{align*}
	where we assume $\kappa_{\ast} \leq m_1 / n \leq \kappa^\ast$ for some universal constants $\kappa_{\ast}, \kappa^\ast >0$. In particular, $1 - \Phi(T_{\mathrm{mean}})$ is a dimension-agnostic $p$-value in the sense of (\ref{eq:formal-goal}).
\end{theorem}
To complement the above result, we recall minimax separation rate for mean testing known in the literature and compare it with our result. Let us define the minimax type II error as 
\begin{align*}
	\textcolor{black}{R(\epsilon_n)} := \inf_{\phi \in \Psi(\alpha)} \textcolor{black}{\limsup_{n \rightarrow \infty}}\sup_{P \in \mathcal{Q}_{\mathrm{mean},d}(\epsilon_n)} \mE_{P}[1 - \phi].
\end{align*}
Given this minimax risk, the minimum separation (or called the critical radius) is given as
\begin{align*}
	\epsilon_n^\dagger := \inf \big\{ \epsilon_n:  \textcolor{black}{R(\epsilon_n)}  \leq 1/2 \big\},
\end{align*}
where the constant $1/2$ can be replaced with any number between $(0,1-\alpha)$. \textcolor{black}{In words, once we commit to a test $\phi$ with asymptotic type I error control at level $\alpha$, the test will be judged based on its asymptotic type II error. For some fast-decaying sequences $\epsilon_n$, the minimax risk~$R(\epsilon_n)$ will remain large close to one, while for slowly decaying sequences $\epsilon_n$, $R(\epsilon_n)$ can be driven to zero. There is typically a sequence of ``critical radii'' where the risk can be brought down to one half, and any slower decay would result in zero risk but any faster decay would cause the risk to exceed one half.} For the Gaussian sequence model, it has been known that the critical radius in the $L_2$ norm is of order $d^{1/4} / \sqrt{n}$ \citep[see e.g.][]{baraud2002non,ingster2012nonparametric,blanchard2018minimax}. In fact, the Gaussian sequence model can be formulated as $X \sim N(\mu, I)$, which satisfies Assumption~\ref{assumption: moment condition} with $L=2\sqrt{2}/\sqrt{\pi}$. Therefore we can conclude that the lower bound matches the upper bound given in Theorem~\ref{Theorem: Uniform type I and II error control}. Consequently the proposed test is minimax rate optimal in terms of the $L_2$ norm.

\begin{remark}[Dimension-to-sample size ratio] \normalfont
	Recall that we assume $d/n \rightarrow \tau \in (0, \infty)$ in Theorem~\ref{Theorem: asymptotic power} under which we can compare the power of $\phi_{\mathrm{mean}}$ with the test using $U_{\mathrm{mean}}$ in an asymptotically precise manner. Without such restriction on $n$ and $d$, the analytic power expressions for both tests are hard to derive and, in fact, the limiting distributions of $U_{\mathrm{mean}}$ and $T_{\mathrm{mean}}$ are not necessarily the same under local alternatives with fixed $d$~(see \cite{wang2019feasible} and Appendix~\ref{Section: Further discussions on power}). In contrast, the uniform consistency result in Theorem~\ref{Theorem: Uniform type I and II error control} holds without the condition~$d/n \rightarrow \tau$.
\end{remark}

\subsection{Non-asymptotic calibration under symmetry}  \label{Section: Non-asymptotic calibration under symmetry}
While our main focus is on asymptotic type I error control and dimension-agnostic property, it is possible to have a finite-sample guarantee using different calibration methods with different assumptions. For example, suppose that $X$ is symmetric about the origin under the null hypothesis. This in turn implies that $f_{\mathrm{mean}}(X)$ and $\varepsilon f_{\mathrm{mean}}(X)$ have the same distribution where $\varepsilon$ is a Rademacher random variable independent of $f_{\mathrm{mean}}(X)$. More generally, for $j=1,\ldots,J$, let $\boldsymbol{\varepsilon}_j:=\{\varepsilon_{1,j},\ldots,\varepsilon_{m_1,j}\}$ be a random vector that consists of i.i.d.~Rademacher random variables independent of the data set. Given mutually independent $\{\boldsymbol{\varepsilon}_j\}_{j=1}^J$, we denote the studentized statistic (\ref{Eq: studentized statistic for mean testing}) computed based on $\{\varepsilon_{1,j} f_{\mathrm{mean}}(X_1),\ldots, \varepsilon_{m_1,j} f_{\mathrm{mean}}(X_{m_1}) \}$ by $T_{\mathrm{mean}}^{\boldsymbol{\varepsilon}_j}$. Under the symmetry assumption, we see that $\{T_{\mathrm{mean}}, T_{\mathrm{mean}}^{\boldsymbol{\varepsilon}_1},\ldots, T_{\mathrm{mean}}^{\boldsymbol{\varepsilon}_J} \}$ are exchangeable and consequently
\begin{align} 
	p_{\mathrm{sym}}^\dagger := \frac{1}{J+1} \Bigg[ 1 + \sum_{j=1}^J \mathds{1} \big\{ T_{\mathrm{mean}} \geq T_{\mathrm{mean}}^{\boldsymbol{\varepsilon}_j} \big\}  \Bigg] \quad \text{and} \quad \phi_{\mathrm{sym}}^{\dagger} :=  \mathds{1} \{ p_{\mathrm{sym}}^\dagger \leq \alpha \} \label{Eq: Monte-Carlo test}
\end{align}
are a valid $p$-value and a level $\alpha$ test, respectively, in finite samples stated as follows. 
\begin{corollary}[Finite sample guarantee of $\phi_{\mathrm{sym}}^{\dagger}$] \label{Corollary: properties of phi_sym_dagger}
	\sloppy Let $\mathcal{P}_{\mathrm{sym},n,d}$ be the set of all symmetric distributions about the origin in $\mathbb{R}^d$. Then, without any moment assumption, it holds that $\sup_{P \in \mathcal{P}_{\mathrm{sym},n,d}}\mE_{P} [\phi_{\mathrm{sym}}^\dagger] \leq \alpha$, and $p_{\mathrm{sym}}^\dagger$ is a dimension-agnostic $p$-value in the sense of (\ref{eq:formal-goal}).
\end{corollary}
Corollary~\ref{Corollary: properties of phi_sym_dagger} follows directly from Lemma 1 of \cite{romano2005exact}. Despite its finite-sample property, we note that $\phi_{\mathrm{sym}}^{\dagger}$ never rejects the null when $\alpha < 1/(J+1)$. Therefore it may be computationally expensive to compute $\phi_{\mathrm{sym}}^{\dagger}$ for small $\alpha$.
Moreover, it is challenging to study the power of $\phi_{\mathrm{sym}}^{\dagger}$ since the critical value is essentially data-dependent. Motivated by these drawbacks, we now derive a computationally efficient but slightly more conservative test building on a classical result of \cite{efron1969student}. More specifically, by imposing Bahadur and Eaton's inequality in \cite{efron1969student}, the self-normalized process $V_n:=\sum_{i=1}^{m_1} f_{\mathrm{mean}}(X_i) / \sqrt{\sum_{i=1}^{m_1}  f_{\mathrm{mean}}^2(X_i)}$ has a sub-Gaussian tail: 
\begin{equation}\label{eq:efron-sym}
	P(V_n \geq t) \leq \exp(-t^2/2)
\end{equation}
under the symmetry assumption. 
Based on this finite-sample bound along with the monotonic relationship $T_{\mathrm{mean}} = \sqrt{m_1} V_n /\sqrt{m_1 - V_n^2}$, it can be seen that the test 
\begin{align*}
	\phi_{\mathrm{sym}} := \mathds{1} \Bigg\{  T_{\mathrm{mean}}  \geq  \sqrt{\frac{2 \log (1/\alpha)}{1 - 2 m_1^{-1}\log(1/\alpha)}} \ \text{and} \  m_1 > 2 \log (1/\alpha)  \Bigg\}
\end{align*}
is equivalent to $\mathds{1}\{V_n > \sqrt{2 \log(1/\alpha)}\}$ and thus $\phi_{\mathrm{sym}}$ is level $\alpha$. 
Moreover, using the fact that $\sqrt{2\log(1/\alpha) } \geq z_{1-\alpha}$ for any $\alpha \in (0,1)$, similar uniform guarantees in Theorem~\ref{Theorem: Uniform type I and II error control} can be established without the symmetry condition. We summarize the properties of $\phi_{\mathrm{sym}}$ in the following corollary. 
\begin{corollary}[Properties of $\phi_{\mathrm{sym}}$] \label{Corollary: properties of phi_sym}
	Recall that $\mathcal{P}_{\mathrm{sym},n,d}$ is the set of all symmetric distributions about the origin in $\mathbb{R}^d$. Then we have $\sup_{P \in \mathcal{P}_{\mathrm{sym},n,d}}\mE_{P} [\phi_{\mathrm{sym}} ] \leq \alpha$ without any moment assumption and $p_{\emph{sym}}:=e^{-V_n^2/2}$ is a dimension-agnostic $p$-value in the sense of (\ref{eq:formal-goal}). Moreover, under the same conditions in Theorem~\ref{Theorem: Uniform type I and II error control}, it holds that
	\begin{align*}
		\lim_{n \rightarrow \infty} \sup_{P \in \mathcal{P}_{0,d}} \mE_{P} [\phi_{\mathrm{sym}}] \leq \alpha \quad \text{and} \quad 	\lim_{n \rightarrow \infty} \sup_{P \in \mathcal{Q}_{\mathrm{mean},d}(\epsilon_n)} \mE_{P} [1 -  \phi_{\mathrm{sym}}] = 0.
	\end{align*}
\end{corollary}
In other words, $\phi_{\mathrm{sym}}$ has a finite-sample guarantee under the symmetry condition and it is minimax rate optimal for one-sample mean testing under the condition in Theorem~\ref{Theorem: Uniform type I and II error control}. While $\phi_{\mathrm{sym}}$ is conservative in general, this may not be a serious issue in ``\emph{large $m_1$ and small $\alpha$}'' scenarios, given that $z_{1-\alpha} / \sqrt{2 \log (1/\alpha)} \rightarrow 1$ as $\alpha \rightarrow 0$.

\subsection{Maximum statistic and its optimality against sparse alternatives} \label{Section: Maximum statistic and its optimality against sparse alternatives}
This section demonstrates our technique using the $L_\infty$-type statistic and proves its minimax rate optimality against sparse alternatives. To this end, recall from \eqref{eq:lp-duality} that the $L_\infty$-norm of $\mu \in \mathbb{R}^d$ has the following variational representation:
\begin{align*}
	\|\mu\|_\infty = \sup_{v \in \mathbb{R}^d: \|v\|_{1} \leq 1} v^\top \mu. 
\end{align*}
As before, let us denote the sample means by $\widehat{\mu}_1$ and $\widehat{\mu}_2$, computed based on the first and second part of the data, respectively. Let $i^\star = \arg\max_{1\leq j \leq d} |\widehat{\mu}_{2}^{(j)}|$, and let $\widehat{v} = (\widehat{v}^{(1)},\ldots,\widehat{v}^{(d)})^\top$ and $\widehat{v}^{(j)} = \text{sign}(\widehat{\mu}_{2}^{(i^\star)})$ for $j=i^\star$ and $\widehat{v}^{(j)} =0$ for $j \neq i^\star$, meaning that $\widehat v$ looks like $(0,0,\dots,0, \pm 1,0,\dots 0)$, with a plus or minus one at position $i^\star$. Then, following the two-step procedure described in Section 1.2, our estimate for $\|\mu\|_\infty$ is $\widehat{v}^\top \widehat{\mu}_1$. We note that our unstudentized statistic $\widehat{v}^\top \widehat{\mu}_1$ is related to the sparse mean statistic studied by \cite{cox1975note} whose focus is on the parametric Gaussian with known variance. In contrast, we do not assume that the variance is known and our focus is on a nonparametric setting. Our statistic is also related to the splitting estimator studied by \cite{rinaldo2019bootstrapping} where the first half of the data is used for model selection and the second half is used for inference.

Similarly as before, let us denote $f_{\mathrm{sparse}}(x) := \widehat{v}^\top x$ and define our studentized statistic as
\begin{align}
	T_{\mathrm{sparse}} := \frac{\sqrt{m_1} \overline{f}_{\mathrm{sparse}}}{\sqrt{m_1^{-1} \sum_{i=1}^{m_1} \{f_{\mathrm{sparse}}(X_i) - \overline{f}_{\mathrm{sparse}}\}^2 }},
\end{align}
where we write $\overline{f}_{\mathrm{sparse}} = m_1^{-1} \sum_{i=1}^{m_1} f_{\mathrm{sparse}}(X_i)$. To study the performance of $T_{\mathrm{sparse}}$, we are concerned with the sparse alternative:
\begin{align*}
	\mathcal{Q}_{\mathrm{sparse},d}(\epsilon_n) :=  \big\{P \in \mathcal{P}_{\mu,d} : \| \mu \|_\infty \geq \epsilon_n \ \text{and $\max_{1\leq j \leq d} \mE_P\big[e^{\lambda (X^{(j)} - \mu^{(j)})}\big] \leq e^{\lambda^2\sigma^2}$ for all $\lambda \in \mathbb{R}$} \big\},
\end{align*}
where the latter condition basically means that the coordinates of $X$ are sub-Gaussian with a common positive parameter $\sigma$. The dependence structure between the coordinates does not matter. We then claim that the sparse mean test $\phi_{\mathrm{sparse}} := \mathds{1} (T_{\mathrm{sparse}} > z_{1-\alpha})$ is asymptotically level $\alpha$ uniformly over $\mathcal{P}_{0,d}$ and the type II error of $\phi_{\mathrm{sparse}}$ is uniformly small over $\mathcal{Q}_{\mathrm{sparse},d}(\epsilon_n)$ when $\epsilon_n$ is sufficiently larger than $\sigma \sqrt{n^{-1}\log d}$. 
\begin{theorem}[Uniform type I and II error control] \label{Theorem: Uniform type I and II error control for sparse alternatives}
	For fixed $\alpha \in (0,1)$, suppose that $\epsilon_n \geq K_n \sigma \sqrt{n^{-1} \log d}$ where $K_n$ is a positive sequence increasing to infinity at any rate. Then
	\begin{align*}
		\lim_{n \rightarrow \infty} \sup_{P \in \mathcal{P}_{0,d}} \mE_{P} [\phi_{\mathrm{sparse}}] = \alpha \quad \text{and} \quad 	\lim_{n \rightarrow \infty} \sup_{P \in \mathcal{Q}_{\mathrm{sparse},d}(\epsilon_n)} \mE_{P} [1 -  \phi_{\mathrm{sparse}}] = 0,
	\end{align*}
	where we assume $\kappa_{\ast} \leq m_1 / n \leq \kappa^\ast$ for some $\kappa_{\ast}, \kappa^\ast >0$. In particular, $1 - \Phi(T_{\mathrm{sparse}})$ is a dimension-agnostic $p$-value in the sense of (\ref{eq:formal-goal}).
\end{theorem}
It is worth pointing out that asymptotic type I error control of $\phi_{\mathrm{sparse}}$ holds regardless of whether $d$ is fixed or not, i.e.~$\phi_{\mathrm{sparse}}$ is dimension-agnostic. It is also important to note that, in view of the minimax lower bound in \cite{cai2014two}, one cannot improve the power of $\phi_{\mathrm{sparse}}$ up to a constant factor, so $\phi_{\mathrm{sparse}}$ is minimax rate optimal against the considered sparse alternative.

Up to now, we have considered the $L_2$ and $L_\infty$ norms of $\mu$ to measure the distance between the null and alternative hypotheses for one-sample mean testing. Indeed, our framework can be readily extended to the $L_p$ norm for any $1\leq p \leq \infty$ using a variational representation. In particular, we can prove that the $L_p$-based test via the two-step procedure preserves dimension-agnostic property over $\mathcal{P}_{0,d}$ as long as the estimated direction vector $\widehat{v}$ is non-zero with probability one. We also expect that the considered test is minimax rate optimal in the $L_p$ norm, but this requires a more delicate analysis. We leave this direction to future work.

\textcolor{black}{We note that our construction of $T_{\mathrm{mean}}$ and $T_{\mathrm{sparse}}$ shares the same principle as the test statistics for mean testing proposed by \cite{huang2015projection,liu2022projection}. In particular, these test statistics are derived by initially identifying an optimal direction, followed by computing a studentized statistic based on the projected univariate random variables along the optimal direction. 
% Nevertheless, their objective is different from ours, leading to Hotelling-type test statistics. 
However, they  appear to \emph{assume} the asymptotic normality of their test statistics (recall Figure~\ref{Figure: CLT under sample splitting}), whereas we offer a rigorous justification for the use of Gaussian calibration in both low- and high-dimensional settings as well as optimal power properties.}

%\begin{figure}[t]
%	\begin{center}
%		\includegraphics[width=1.00\textwidth]{CI.pdf}
%	\end{center}
%	\caption{Illustration of confidence sets $\mathcal{C}_{n,\text{normal}}(\mathcal{X}_1, \mathcal{X}_2;\alpha)$ (left), $\mathcal{C}_{n,\text{cross}}(\mathcal{X}_1, \mathcal{X}_2;\alpha)$ (middle) and $\mathcal{C}_{n,\mathrm{sym}}(\mathcal{X}_1, \mathcal{X}_2;\alpha)$ (right) for $\mu=(\mu^{(1)},\mu^{(2)})^\top$ at $\alpha=0.05$, $0.15$ and $0.25$. These confidence sets are (asymptotically) valid independent of $d$. We focus on the bivariate case for visualization only.} \label{Figure: CI}
%\end{figure}

Next, we switch gears to tackle the problem of one-sample covariance testing. This is mainly in order to demonstrate that our approach was not specialized to just one problem setting, and the methodology and theory generalizes cleanly to more complex settings.

\section{One-sample Gaussian covariance testing} \label{Section: One-sample covariance testing}
We now study the problem of covariance testing under Gaussian assumptions as another application of our approach via sample-splitting and studentization, though we show how to relax the Gaussianity assumption later. Let $X_1,\ldots,X_n$ be i.i.d.~$d$-dimensional random vectors from a multivariate normal distribution $N(0,\Sigma)$. Given this Gaussian sample, we are concerned with testing whether
\begin{align*}
	H_0 : \Sigma = \Sigma_0 \quad \text{versus} \quad H_1: \Sigma \neq \Sigma_0.
\end{align*}
Without loss of generality, we focus on the case of $\Sigma_0 = I$ as one can work with the transformed data $\widetilde{X} = \Sigma_0^{-1/2}X$ such that $\mathrm{cov}[\widetilde{X} ] = I$ whenever $\Sigma_0$ is of full rank. Testing for the equality of covariance matrices is one of the important topics in statistics and it has been actively investigated both in fixed-dimensional settings~\citep[e.g.][]{nagao1973some,anderson2003introduction} and high-dimensional settings~\citep[e.g.][]{birke2005note,chen2010tests,cai2013optimal}. In particular, Cai and Ma~\cite{cai2013optimal} consider the following U-statistic as a test statistic
\begin{align} \label{Eq: CM statistic}
	U_{\mathrm{cov}} := \frac{1}{n(n-1)} \sum_{1 \leq i \neq j \leq n} \tr\big\{ (X_iX_i^\top - I) (X_jX_j^\top - I) \big\},
\end{align}
and show that it converges to a Gaussian distribution in a high-dimensional regime where $d$ can grow independent of $n$. They further prove that their test calibrated by this Gaussian approximation achieves minimax rate optimality in the Frobenius norm. However the validity of their test crucially relies on the assumption that $d$ increases to infinity. In fact, when $d$ is fixed, the theory of degenerate U-statistics reveals that $U_{\mathrm{cov}}$ converges to a sum of chi-square random variables. See also the second row of Figure~\ref{Figure: Covariance} where we provide some simulation results to support this claim, demonstrating that their test is uncalibrated outside of their specific high-dimensional regime. In addition the Berry--Esseen bound for $U_{\mathrm{cov}}$ later in (\ref{Eq: BE bound for U2}) suggests that $U_{\mathrm{cov}}$ has a slow rate of convergence to the normal distribution.

\subsection{A test statistic by sample-splitting and studentization}

Motivated by the aforementioned issue, we propose a sample-splitting analogue of the U-statistic~(\ref{Eq: CM statistic}), which has a Gaussian limit under both fixed- and high-dimensional regimes. First, define 
\begin{align*}
	T_{\mathrm{cov}}^{\dagger} :=  \frac{1}{m_1m_2} \sum_{i=1}^{m_1} \sum_{j=m_1+1}^n \tr\big\{ (X_iX_i^\top - I) (X_jX_j^\top - I) \big\},
\end{align*}
which is an unbiased estimator of the squared Frobenius norm $\| \Sigma - I \|_F^2$. Similarly as before in Section~\ref{Section: Mean Testing}, we define a random function depending solely on $\mathcal{X}_2$ by
\begin{align} \label{Eq: g function}
	f_{\mathrm{cov}}(x) := \frac{1}{m_2} \sum_{j=m_2+1}^n \tr\big\{ (xx^\top - I) (X_jX_j^\top - I) \big\},
\end{align}
so that $T_{\mathrm{cov}}^{\dagger}$ simply equals $\overline{f}_{\mathrm{cov}}:=m_1^{-1} \sum_{i=1}^{m_1} f_{\mathrm{cov}}(X_i)$. We then consider its studentized version:
% of $T_{\mathrm{cov}}^{\dagger}$ defined as
\begin{align}
	T_{\mathrm{cov}} := \frac{\sqrt{m_1} \overline{f}_{\mathrm{cov}}}{\sqrt{m_1^{-1} \sum_{i=1}^{m_1} \big\{ f_{\mathrm{cov}}(X_i) - \overline{f}_{\mathrm{cov}} \big\}^2}},
\end{align}
and provide analogous results in Section~\ref{Section: Mean Testing}. 

\subsection{Asymptotic, dimension-agnostic, Gaussian null distribution} 
We start by exploring the asymptotic null distribution of $T_{\mathrm{cov}}$. Following a similar strategy used in Section~\ref{Section: Mean Testing}, we consider a conditional Berry--Esseen bound \citep[Theorem 1.1 of][]{bentkus1996berry} as
\begin{align} \label{Eq: BE bound for covariance}
	\sup_{-\infty < z < \infty} \big| P_{0}(T_{\mathrm{cov}}  \leq z | \mathcal{X}_{2}) - \Phi(z) \big| \leq  C \frac{\mE_{P_{0}}[|f_{\mathrm{cov}}(X)|^3|\mathcal{X}_{2}]}{\{\mE_{P_{0}}[|f_{\mathrm{cov}}(X)|^2 | \mathcal{X}_{2}] \}^{3/2} \sqrt{m_1}},
\end{align}
where $C>0$ is an absolute constant and we write $P_{0} = N(0,I)$ throughout this section. Therefore $T_{\mathrm{cov}}$ approximates the normal distribution as long as the right-hand side of the bound~\eqref{Eq: BE bound for covariance} decreases to zero as $m_1 \rightarrow \infty$. In this application, it is more convenient to deal with the fourth moment rather than the third moment and thus we work with 
\begin{align} \label{Eq: BE bound for covariance 2}
	\frac{\mE_{P_{0}}[|f_{\mathrm{cov}}(X)|^3|\mathcal{X}_{2}]}{\{\mE_{P_{0}}[|f_{\mathrm{cov}}(X)|^2 | \mathcal{X}_{2}] \}^{3/2} \sqrt{m_1}} \leq \left( \frac{\mE_{P_{0}}[|f_{\mathrm{cov}}(X)|^4|\mathcal{X}_{2}]}{\{ \mE_{P_{0}}[|f_{\mathrm{cov}}(X)|^2 | \mathcal{X}_{2}] \}^2 m_1} \right)^\frac{1}{2},
\end{align}
which is a consequence of the Cauchy--Schwarz inequality. After a careful moment analysis under the Gaussian assumption detailed in the proof of Proposition~\ref{Proposition: BE bound for T2}, we can obtain a simple upper bound for the scaled fourth moment~(\ref{Eq: BE bound for covariance 2}):
\begin{align} \label{Eq: Sufficient condition for Normality}
	\frac{\mE_{P_{0}}[|f_{\mathrm{cov}}(X)|^4|\mathcal{X}_{2}]}{\{\mE_{P_{0}}[|f_{\mathrm{cov}}(X)|^2 | \mathcal{X}_{2}] \}^2 m_1}  \leq \frac{15}{m_1} \quad \text{a.s.}
\end{align}
\textcolor{black}{In fact, the above bound holds over a broad class of distributions, including a Gaussian distribution as a special case, but with a larger constant factor. This, along with the Berry--Esseen bound (\ref{Eq: BE bound for covariance}), implies that $T_{\mathrm{cov}}$ has a Gaussian limiting distribution uniformly over the class for which bound~(\ref{Eq: Sufficient condition for Normality}) holds up to a constant factor. We formally state the result below.}
\begin{proposition}[Berry--Esseen bound for $T_{\mathrm{cov}}$] \label{Proposition: BE bound for T2}
	Consider the class of distributions $\mathcal{P}_{0,d,\mathrm{cov}}$ defined in Definition~\ref{Definition: covariance model}. Then there exists a constant $C > 0$ such that 
	\begin{align} \label{Eq: BE bound for T2}
		\sup_{d \geq 1}\sup_{P \in \mathcal{P}_{0,d,\mathrm{cov}}}\sup_{-\infty < z < \infty} \big| P(T_{\mathrm{cov}}  \leq z) - \Phi(z) \big| \leq \frac{C}{\sqrt{m_1}}.
	\end{align}	
\end{proposition}
The above result later yields dimension-agnostic inference in Proposition~\ref{Proposition: Uniform type I and II error control of the covariance test}.
Proposition~\ref{Proposition: BE bound for T2} shows that the studentized statistic~$T_{\mathrm{cov}}$ has a Gaussian limiting distribution as $m_1 \rightarrow \infty$. Similar to Theorem~\ref{Theorem: unconditional BE Bound}, this asymptotic guarantee is independent of the assumption on $d$ and $m_2$. Furthermore its convergence rate is of the same order as the usual sum of i.i.d.~random variables. In contrast, we note that the Berry--Esseen bound for $U_{\mathrm{cov}}$ \citep{cai2013optimal}:
\begin{align} \label{Eq: BE bound for U2}
	\sup_{-\infty < z < \infty} \bigg| P_{0} \bigg( \frac{U_{\mathrm{cov}}}{\sqrt{\mV[U_{\mathrm{cov}}]} } \leq z \bigg) - \Phi(z) \bigg| \leq  C \left[ \frac{1}{n} + \frac{1}{d} \right]^{1/5}
\end{align}
converges to zero much slower than our bound~(\ref{Eq: BE bound for T2}) and also depends on $d$. At this point, it is unknown whether bound~(\ref{Eq: BE bound for U2}) can be improved by a refined analysis. Nevertheless, it is crucial to assume $d \rightarrow \infty$ for the asymptotic normality of $U_{\mathrm{cov}}$. Our approach, in contrast, does not put any assumption on $d$, and thus it is valid under a more general asymptotic regime of $(n,d)$ than  \cite{cai2013optimal}. We also note that the class of distributions $\mathcal{P}_{0,d,\mathrm{cov}}$ is similarly considered in \cite{chen2010tests} for testing covariance matrices and it includes $P_{0} \equiv N(0,1)$ as an example.

\subsection{Asymptotic power analysis} 
Next we derive the asymptotic power expression of the test based on $T_{\mathrm{cov}}$, denoted by $\phi_{\mathrm{cov}} := \mathds{1}(T_{\mathrm{cov}} > z_{1-\alpha})$. In view of Proposition~\ref{Proposition: BE bound for T2}, it is clear that $\phi_{\mathrm{cov}}$ asymptotically controls the type I error rate both in fixed- and high-dimensional regimes. To provide an explicit power expression, we make similar assumptions to Assumption~\ref{assumption: one-sample mean testing} and consider alternatives where $\| \Sigma - I \|_F = O(\sqrt{d/n})$. Under these settings, we prove that $\phi_{\mathrm{cov}}$ has comparable power to the test by \cite{cai2013optimal}.

\begin{theorem}[Asymptotic power expression] \label{Theorem: asymptotic power of covariance test}
	Suppose that Assumption~\ref{assumption: one-sample mean testing} (a) with $\mu=0$, (b), (d) and (e) are fulfilled under the alternative where $\| \Sigma - I \|_F = O(\sqrt{d/n})$. Then it holds that 
	\begin{align*}
		\mE_P[\phi_{\mathrm{cov}}] ~=~  \Phi \left( z_{\alpha} + \frac{n\sqrt{\kappa(1-\kappa)}\|\Sigma - I\|_F^2}{\sqrt{2}\emph{\tr}(\Sigma^2)} \right) + o(1).
	\end{align*}
	Therefore the power is asymptotically maximized when $m_1/n \rightarrow \kappa =1/2$. 
\end{theorem}

By choosing the optimal choice of $\kappa = 1/2$, the power function of $\phi_{\mathrm{cov}}$ simplifies to
\begin{align} \label{Eq: theoretical power of phi_cov}
	\mE_P[\phi_{\mathrm{cov}}] ~=~  \Phi \left( z_{\alpha} + \frac{n\|\Sigma - I\|_F^2}{2\sqrt{2}\tr(\Sigma^2)} \right) + o(1).
\end{align}
On the other hand, the power of the test by \cite{cai2013optimal} is given by
\begin{align} \label{Eq: theoretical power of U_2}
	P \left( \! U_{\mathrm{cov}} >  2 z_{1-\alpha} \sqrt{\frac{d(d+1)}{n(n-1)}} \right) ~=~ \Phi \left( \! z_{\alpha}   +  \frac{ n \|\Sigma - I \|_F^2}{2\tr(\Sigma^2)} \right) + o(1).
\end{align}
Again, similar to the case of mean testing, we see that the power of $\phi_{\mathrm{cov}}$ is worse by a $\sqrt{2}$ factor. This specific constant factor can be explained in a similar fashion to Remark~\ref{Remark: Intuition of constant factor}. 

\subsection{Minimax rate optimality against Frobenius-norm deviations} 
Lastly we demonstrate that $\phi_{\mathrm{cov}}$ has nontrivial power when $\| \Sigma - I\|_F$ is of order $\sqrt{d/n}$ without making assumptions on the ratio of $d/n$ and the eigenvalues of $\Sigma$. To formally state the result, for a sequence of $\epsilon_n >0$, let us define a class of normal distributions such that 
\begin{align*}
	\mathcal{Q}_{\mathrm{cov},d}(\epsilon_n) := \big\{ P =N(0,\Sigma) : \| \Sigma - I \|_F \geq \epsilon_n  \big\}.
\end{align*}
The following proposition shows that $\phi_{\mathrm{cov}}$ controls the type I error rate in any asymptotic regime and its type II error is uniformly small when $\epsilon_n$ is sufficiently larger than $\sqrt{d/n}$. 
\begin{proposition}[Type I and II error control] \label{Proposition: Uniform type I and II error control of the covariance test}
	For fixed $\alpha \in (0,1)$, suppose that $\epsilon_n \geq K_n \sqrt{d/n}$ where $K_n$ is a positive sequence increasing to infinity as $n \rightarrow \infty$ at any rate. Then we have 
	\begin{align}
		\lim_{n \rightarrow \infty} \mE_{P_{0}} [\phi_{\mathrm{cov}}] = \alpha \quad \text{and} \quad \lim_{n \rightarrow \infty} \sup_{P \in \mathcal{Q}_{\mathrm{cov},d}(\epsilon_n)} \mE_{P}[1 -  \phi_{\mathrm{cov}}] =0,
	\end{align}
	where we assume $\kappa_{\ast} \leq m_1 / n \leq \kappa^\ast$ for some constants $\kappa_{\ast}, \kappa^\ast >0$. In particular, $1 - \Phi(T_{\mathrm{cov}})$ is a dimension-agnostic $p$-value in the sense of (\ref{eq:formal-goal}).
\end{proposition}
\textcolor{black}{We note that, while we state the type I error result under the Gaussian assumption, it is clear from Proposition~\ref{Proposition: BE bound for T2} that the same result holds uniformly over $\mathcal{P}_{0,d,\mathrm{cov}}$.} We also remark that \cite{cai2013optimal} provide a lower bound for the separation boundary $\epsilon_n$ in the Frobenius norm. In detail, they show that if $\epsilon_n$ is much smaller than $\sqrt{d/n}$, no test has uniform nontrivial power against $\mathcal{Q}_{\mathrm{cov},d}(\epsilon_n)$. In this sense our covariance test $\phi_{\mathrm{cov}}$ achieves minimax rate optimality in terms of the Frobenius norm. The same rate optimality is achieved by the test $U_{\mathrm{cov}} >  2 z_{1-\alpha} \sqrt{d(d+1)/n(n-1)}$, as proved by \cite{cai2013optimal}, but we should emphasize that this test is only valid when $d$ goes to infinity with $n$ (see Figure~\ref{Figure: Covariance}).

Since we have a dimension-agnostic test for any positive definite $\Sigma$, the duality of tests and confidence sets allows us to construct a dimension-agnostic confidence set for $\Sigma$ as in Section~\ref{Section: Dimension-agnostic confidence set}. We omit the details due to space limit. In principle, this confidence set can be used to test for composite nulls such as $H_0: \|\Sigma\|_F \leq c_1$ or $H_0: c_2 \leq \|\Sigma\|_F \leq c_1$ for some constants $c_1,c_2>0$.

\section{General results for degenerate \textcolor{black}{cross} U-statistics} \label{Section: General results for degenerate U-statistics}
%So far we have focused on relatively simple problems for which it is possible to develop concrete results such as the dimension-agnostic property. Moving one step further, the goal of this section is to present general conditions under which a sample-splitting analogue of a degenerate U-statistic has a Gaussian limiting distribution. We demonstrate these conditions in the context of goodness-of-fit testing based on Gaussian MMD~\citep{gretton2012kernel,kellner2019one} and also study its minimax rate optimality in the $L_2$ distance. %building on the recent result of \cite{li2019optimality}. 
This section presents general conditions under which a sample splitting analogue of a degenerate U-statistic has a Gaussian limit. We demonstrate these conditions in the context of goodness-of-fit testing based on Gaussian MMD~\citep{gretton2012kernel,kellner2019one} and study its minimax rate optimality in the $L_2$ distance.

\subsection{Cross U-statistics via sample-splitting and studentization}
First, we briefly review the asymptotic theory for a degenerate U-statistic and introduce the proposed statistic via sample splitting and studentization. To formalize the results, consider a symmetric function $h(x,y)$ that satisfies $\mE_P[h(X,y)] = 0$ almost surely. Then the corresponding U-statistic with $h(x,y)$, which is a degenerate U-statistic of order one, is given by
\begin{align} \label{Eq: U-statistic}
	U_h := \frac{1}{n(n-1)}\sum_{1 \leq i \neq j \leq n} h(X_i,X_j),
\end{align}
where $\{X_i\}_{i=1}^n \overset{i.i.d.}{\sim} P$. By choosing $h(x,y) = x^\top y$ or $h(x,y) = \tr\{ (xx^\top - I)(yy^\top - I) \}$, $U_h$ becomes $U_{\mathrm{mean}}$ or $U_{\mathrm{cov}}$ introduced in (\ref{Eq: One-sample Mean U-statistic}) and (\ref{Eq: CM statistic}), respectively. Other examples of degenerate U-statistics include the one-sample version of kernel MMD~\citep{gretton2012kernel,kellner2019one}, energy distance for multivariate Gaussians~\citep{szekely2005new}, kernelized Stein discrepancy \citep{chwialkowski2016kernel,liu2016kernelized} and kernel estimators of a quadratic functional~\citep{hall1987estimation}. The limiting distribution of $U_h$ is typically intractable, having an infinite number of unknown parameters. Specifically, by the spectral decomposition of $h$ under finite second moment, 
\begin{align} \label{Eq: Another representation of h}
	h(x,y) = \sum_{k=1}^\infty \lambda_k \phi_k(x) \phi_k(y),
\end{align} 
where $\{\lambda_k\}_{k=1}^\infty$ and $\{\phi_k(\cdot)\}_{k=1}^\infty$ are eigenvalues and eigenfunctions of the integral equation given as $\mE_P[h(X,y)f(X)]  = \lambda f(y)$. Under a conventional setting where $h(x,y)$ is fixed in the sample size $n$ with finite second moment, the U-statistic converges to a Gaussian chaos as $nU_h \convD \sum_{k=1}^\infty \lambda_k (\xi_k^2-1)$ \citep[see e.g.][for details]{lee1990u}. On the other hand, suppose $h(x,y)$ changes with $n$ and satisfies
\begin{subequations}
	\begin{align}  \label{Eq: normality condition 1}
		& \frac{\mE_P[\{ \mE_P[h(X_1,X_2)h(X_1,X_3)|X_2,X_3] \}^2]}{\{ \mE_P [h^2(X_1,X_2)] \}^2} \rightarrow 0  \quad \text{and}  \\[.5em]  \label{Eq: normality condition 2}
		& \frac{\mE_P[h^4(X_1,X_2)]n^{-1} + \mE_P[h^2(X_1,X_3)h^2(X_2,X_3)] }{n\{ \mE_P [h^2(X_1,X_2)] \}^2} \rightarrow 0 \quad \text{as $n \rightarrow \infty$.}
	\end{align}
\end{subequations}
%\begin{equation}
%\begin{aligned}
%& \frac{\mE[\{ \mE[h(X_1,X_2)h(X_1,X_3)|X_2,X_3] \}^2]}{\{ \mE [h^2(X_1,X_2)] \}^2} \rightarrow 0 \quad \text{and} \\[.5em]
%& \frac{\mE[h^4(X_1,X_2)]n^{-1} + \mE[h^2(X_1,X_3)h^2(X_2,X_3)] }{n\{ \mE [h^2(X_1,X_2)] \}^2} \rightarrow 0 \quad \text{as $n \rightarrow \infty$,}
%\end{aligned}
%\end{equation}
Then, as shown in, for example, \cite{zhang2018conditional}, the scaled $U_h$ converges to a normal distribution as
\begin{align*}
	\frac{nU_h}{\sqrt{\mE_P[h^2(X_1,X_2)]/2}} \convD N(0,1).
\end{align*}
\textcolor{black}{It is worth pointing out that the asymptotic normality of $U_h$ was originally established by \cite{hall1984central} under slightly stronger conditions. In particular, \cite{hall1984central} proves the same result under condition~(\ref{Eq: normality condition 1}) and 
	\begin{align} \label{Eq: normality condition 3}
		\frac{\mE_P[h^4(X_1,X_2)]}{n\{ \mE_P [h^2(X_1,X_2)] \}^2} \rightarrow 0, 
	\end{align}
	which implies condition~(\ref{Eq: normality condition 2}) via the Cauchy--Schwarz inequality.} The conditions in (\ref{Eq: normality condition 1}), (\ref{Eq: normality condition 2}) and (\ref{Eq: normality condition 3}) frequently appear in the literature, studying an asymptotic normality of U-statistics  \citep[e.g.][]{wang2015conditional,li2019optimality,kim2020multinomial}. We also refer to \cite{robins2016asymptotic} who provide different conditions that lead to the asymptotic normality of $U_h$. This complicated asymptotic behavior naturally raises a critical question of how to calibrate a test based on $U_h$. One possible remedy for this issue is to use data-driven calibration methods such as bootstrapping or subsampling, but these methods are computationally expensive to implement.

Motivated by the above challenge associated with $U_h$, we introduce our statistic: 
\begin{align} \label{Eq: proposed statistic}
	T_h^\dagger := \frac{1}{m(n-m)} \sum_{i=1}^m \sum_{j=m+1}^n h(X_i, X_j).
\end{align}
We name this statistic as cross U-statistic as it only considers cross-terms in a kernel matrix~(Figure~\ref{Figure: visualization}). Here we focus on $m=\floor{n/2}$ motivated by its maximum asymptotic power property in Theorem~\ref{Theorem: asymptotic power} and Theorem~\ref{Theorem: asymptotic power of covariance test}. Intuitively, conditional on the second half of the data, one can expect that $T_h^\dagger$ behaves like a Gaussian since $T_h^\dagger$ can be viewed as a linear statistic. 

Our goal below is to make the aformentioned intuition rigorous and show that a studentized $T_h^\dagger$ has a simple Gaussian distribution regardless of whether the kernel $h$ depends on the sample size $n$ or not. In particular, the driving force of different behavior of $U_h$ is the condition~(\ref{Eq: normality condition 1}), which is equivalent to $\lambda_1^2/\sum_{k=1}^\infty \lambda_{k}^2$ $\rightarrow 0$ where $\lambda_1$ is the maximum eigenvalue. We prove that the asymptotic normality for our statistic does not require the condition~(\ref{Eq: normality condition 1}).

%We note in passing that the relationship between the U-statistic and our sample-splitting analogue resembles that between stability selection~\citep{meinshausen2010stability} and complementary pairs stability selection~\citep{shah2013variable}; the former builds on all possible subsamples of the same size, whereas the latter uses complementary independent subsamples in order to weaken the assumptions of original stability selection.

%We note in passing that our purpose of sample-splitting resembles \cite{shah2013variable}, which weakens the assumptions of original stability selection \citep{meinshausen2010stability} by choosing complementary independent subsamples.

\begin{remark}[Intrinsic dimension] \normalfont \label{Remark: intrinsic dimension}
	Since the domain of $h(x,y)$ can be highly abstract, it is unclear how to define ``dimension''. To avoid this ambiguity, we define the intrinsic dimension of the given data associated with $h(x,y)$ as the sum of the squared eigenvalues as 
	\begin{align} \label{Eq: intrinsic dimension}
		D_{n}:=\sum_{k=1}^\infty \lambda_{k,n}^2,
	\end{align}
	which is the same as $\mE[h^2(X_1,X_2)]$. In the case of the bivariate kernel $h(x,y) = x^\top y$ with bounded eigenvalues, the intrinsic dimension equals the original dimension of the data up to a constant factor, meaning that $D_n \propto d_n$. As noted before, the limiting behavior of $U_h$ differs depending on the intrinsic dimension $D_{n}$; indeed in several settings, $d_n$ can stay constant, but $D_n$ can stay finite or approach infinity, the latter behavior driving the limiting distribution. When $D_n$ is fixed and finite, $U_h$ is asymptotically close to a Gaussian chaos, whereas when $D_n$ increases to infinity in a way that the condition~(\ref{Eq: normality condition 1}) holds, the scaled $U_h$ converges to a Gaussian. In contrast, as we shall see, the limiting distribution of a studentized $T_h^\dagger$ is invariant to $D_n$, and this is a nontrivial result in our opinion. 
\end{remark}

\noindent In order to study the limiting distribution of $T_h^\dagger$, define a function $f_h$ as
\begin{align} \label{Eq: general linear function}
	f_h(x) := \frac{1}{n-m} \sum_{j=m+1}^n h(x,X_j).
\end{align}
Given this univariate function, $T_h^\dagger$ can be written as the sample mean of $f_h(X_1),\ldots,f_h(X_m)$, denoted by $\overline{f}_h = m^{-1} \sum_{i=1}^m f_h(X_i)$, and its studentized version can be defined as
\begin{align} \label{Eq: studentized statistic with h}
	T_{h} := \frac{\sqrt{m} \overline{f}_h}{\sqrt{m^{-1} \sum_{i=1}^m \{ f_h(X_i) - \overline{f}_h \}^2 }}.
\end{align}
Generalizing techniques from previous sections, we can now calculate its null distribution.

\subsection{A dimension-agnostic, asymptotic Gaussian limit} 

Since $f_h(X_1),\ldots,f_h(X_m)$ are i.i.d.~random variables centered at zero conditional on $\mathcal{X}_2$, we may apply the conditional Berry--Esseen bound in Lemma~\ref{Lemma: Conditional BE Bound} to obtain a Gaussian limiting distribution of $T_{h}$. However it is worth noting that there might be a nontrivial chance that the conditional variance of $f_h$ becomes zero in finite samples depending on situations. In such a case, the Berry--Esseen bound may not be directly applicable without any additional assumption (e.g.~continuous distributions in Assumption~\ref{assumption: moment condition}). Due to this issue, we focus on an asymptotic scenario and prove that $T_{h}$ converges to the standard normal distribution under some moment conditions.

\begin{theorem}[Gaussian Approximation] \label{Theorem: Gaussian Approximation}
	Let $\mathcal{P}_{n,h}$ be a set of distributions associated with kernel $h$ that can potentially change with $n$. Suppose that any sequence of distributions $P_n \in \mathcal{P}_{n,h}$ satisfies (i)~$\mE_{P_n}[h(X,y)]=0$ for each $n$ and (ii)~the condition~(\ref{Eq: normality condition 2}). % (iii)~$\lambda_{1,n}^2 / \sum_{k=1}^\infty \lambda_{k,n}^2$ has a limit where $\{\lambda_{k,n}\}_{k=1}^\infty$ are eigenvalues defined in (\ref{Eq: Another representation of h}). 
	Then
	\begin{align*}
		\sup_{P \in \mathcal{P}_{n,h}} \sup_{-\infty < z < \infty} \big| P(T_{h}  \leq  z) - \Phi(z) \big| \rightarrow 0 \quad \text{as $n \rightarrow \infty$.}
	\end{align*}
	In particular, $1 - \Phi(T_{h})$ is a dimension-agnostic $p$-value in the sense of (\ref{eq:formal-goal}) where the null sequence is given as $H_{0,n}: \mathcal{P}_{n,h}$ and $d_n$ is replaced by the intrinsic dimension $D_n$~(\ref{Eq: intrinsic dimension}). 
\end{theorem}

A few remarks are in order. Note that our purpose is to provide unified conditions for the normality of $T_{h}$ in both fixed and increasing (intrinsic) dimensions. Depending on the regime of interest and the type of kernel $h$, the given conditions can be sharpened. We also note that the  condition~(\ref{Eq: normality condition 2}) is required for $U_h$ as well. However, in contrast to the normality conditions for $U_h$, our result does not require (\ref{Eq: normality condition 1}), which can only happen when the intrinsic dimension goes to infinity. When the kernel $h$ does not depend on $n$ (i.e.~fixed intrinsic dimension), the condition~(\ref{Eq: normality condition 2}) is easily satisfied given that \[0 < \mE_P[h^2(X_1,X_2)] \leq \sqrt{\mE_P[h^4(X_1,X_2)]} < \infty.\] The condition on the first eigenvalue is mainly for a technical reason, which is always true when the kernel $h$ is fixed in $n$. In summary, $T_{h}$ has a Gaussian limiting distribution regardless of whether the intrinsic dimension is fixed or not under the given conditions in Theorem~\ref{Theorem: Gaussian Approximation}.

%\begin{figure}[h!]
%	\centering
%	% \begin{subfigure}[b]{0.31\textwidth}
%	\includegraphics[width=0.31\textwidth]{SSS2I.pdf}
%	\includegraphics[width=0.31\textwidth]{SSS100I.pdf}
%	\includegraphics[width=0.31\textwidth]{SSS1000I.pdf}
%	% \end{subfigure}
%	% \vskip 1em
%	% \begin{subfigure}[b]{0.31\textwidth}
%	\includegraphics[width=0.31\textwidth]{CZ2I.pdf}
%	% \end{subfigure}
%	% ~ %add desired spacing between images, e. g. ~, \quad, \qquad, \hfill etc. 
%	% %(or a blank line to force the subfigure onto a new line)
%	% \begin{subfigure}[b]{0.31\textwidth}
%	\includegraphics[width=0.31\textwidth]{CZ100I.pdf}
%	% \end{subfigure}
%	% ~ %add desired spacing between images, e. g. ~, \quad, \qquad, \hfill etc. 
%	% %(or a blank line to force the subfigure onto a new line)
%	% \begin{subfigure}[b]{0.31\textwidth}
%	\includegraphics[width=0.31\textwidth]{CZ1000I.pdf}
%	% \end{subfigure}
%	\caption{\small QQ plots of $T_{\mathrm{cov}}$ (first row) and $U_{\mathrm{cov}}$ scaled by its standard deviation $2\sqrt{d(d+1)/n(n-1)}$ (second row) under scenario I. The distribution of $T_{\mathrm{cov}}$ consistently approximates $N(0,1)$ in all dimensions, whereas the distribution of $U_{\mathrm{cov}}$ is far from $N(0,1)$ when $d=2$. The straight line $y=x$ is added as a reference point.} \label{Figure: Covariance}
%\end{figure}

\section{Conclusion} \label{Section: Discussion}

This paper defined the phrase ``dimension-agnostic inference'' and made a case for it being a goal worthy of pursuit. We then developed a new dimension-agnostic approach to multivariate testing and estimation problems based on sample splitting and studentization, resulting in ``cross U-statistics''. The proposed test statistics have a simple Gaussian limiting null distribution irrespective of the scaling of the dimensionality. This key property allows us to calibrate the test easily and reliably without any extra computation. Through the paper, we  examined the power and type I error of the proposed tests applied to several problems, demonstrating their competitive performance. 

With some effort, we have recently extended the results in this paper to nonparametric two-sample and independence testing~\cite{shekhar2022permutation,shekhar2022permutation2}. Two important open directions include:
\begin{itemize}
	\item \emph{U-statistics with high order kernels.} While our focus was on U-statistics with a second order kernel, our framework can be possibly generalized to higher-order kernels, which have many applications. A few examples include independence testing~\citep{gretton2008kernel,huo2016fast}, testing for regression coefficients ~\citep{zhong2011tests}, empirical risk estimation~\citep{papa2015sgd} and ensemble methods~\citep{mentch2016quantifying,wager2018estimation}. 
    % The limiting distribution of high-order U-statistics is nontrivial, especially when their kernel is degenerate (at least asymptotically), our techniques can be beneficial for obtaining a simple limiting distribution in these applications. 
    We expect our approach will lead to the dimension-agnostic property and also minimax optimal power.  
 \vskip .5em
	% \item \emph{Two-sample testing.} With some effort, our framework can potentially be extended to nonparametric two-sample testing. Common test statistics, such as the Gaussian MMD~\citep{gretton2008kernel}, energy distance~\citep{baringhaus2004new,szekely2013energy}, and variants of the t-test~\citep{chen2010two}, are based on degenerate kernels and thus their limiting null distributions are sensitive to the assumption on the dimensionality. \vskip .5em
	% \item \emph{Dimension reduction.} As discussed in Section~\ref{Section: Some guiding intuition for developing dimension-agnostic tests}, we used sample splitting as a tool for dimensionality reduction, by mapping the original space to the real line that maximizes an objective function based on one part of the data. We then project the other part of the data onto the real line following the map constructed from the first stage. Lastly we perform a univariate test to see whether there is evidence towards the alternative hypothesis. As mentioned before in Section~\ref{Section: Some guiding intuition for developing dimension-agnostic tests}, all of the examples considered in this paper maximize an integral probability metric defined over a class of functions in a reproducing kernel Hilbert space. From this perspective, our framework can be further extended to other applications using the integral probability measure with a different function class \citep[e.g.][]{sriperumbudur2012empirical} and more generally to applications where the underlying quantity has a variational representation.  \vskip .5em
 
	\item \emph{Multiple sample splitting.} 
 % Recall that our test statistics are based on a single split. 
 To improve statistical efficiency and reduce the randomness from the single split, it is natural to think about multiple sample splitting where we repeat the procedure over different splits and then aggregate the resulting test statistics in some way. However, we may lose the dimension-agnostic property by doing so. 
 % Interestingly, this is true even for a simple cross-fitted statistic. 
 In particular, we show in Proposition~\ref{Proposition: limiting distribution of cross-splitting} that even a simple cross-fitted statistic converges to the sum of two Gaussians. We also discuss different ways to non-asymptotically calibrate an aggregated test statistic for one-sample mean testing in Appendix~\ref{Section: Multiple sample-splitting}. We do not know in general if one can calibrate many exchangeable tests from different splits with provably higher power than a single split. Developing such methods may improve other procedures that use single sample splitting \citep[e.g.][]{decrouez2014split,tansey2018holdout,lei2018distribution,jankova2020goodness,katsevich2020theoretical,kim2016classification}. 
\end{itemize}

%%%%%%%%%%%%%%%%%%%%%%%%%%%%%%%%%%%%%%%%%%%%%%
%% Example with single Appendix:            %%
%%%%%%%%%%%%%%%%%%%%%%%%%%%%%%%%%%%%%%%%%%%%%%
%\begin{appendix}
%\section*{Title}\label{appn} %% if no title is needed, leave empty \section*{}.
%Appendices should be provided in \verb|{appendix}| environment,
%before Acknowledgements.
%
%If there is only one appendix,
%then please refer to it in text as \ldots\ in the \hyperref[appn]{Appendix}.
%\end{appendix}
%%%%%%%%%%%%%%%%%%%%%%%%%%%%%%%%%%%%%%%%%%%%%%%
%%% Example with multiple Appendixes:        %%
%%%%%%%%%%%%%%%%%%%%%%%%%%%%%%%%%%%%%%%%%%%%%%%
%\begin{appendix}
%\section{Title of the first appendix}\label{appA}
%If there are more than one appendix, then please refer to it
%as \ldots\ in Appendix \ref{appA}, Appendix \ref{appB}, etc.
%
%\section{Title of the second appendix}\label{appB}
%\subsection{First subsection of Appendix \protect\ref{appB}}
%
%Use the standard \LaTeX\ commands for headings in \verb|{appendix}|.
%Headings and other objects will be numbered automatically.
%\begin{equation}
%\mathcal{P}=(j_{k,1},j_{k,2},\dots,j_{k,m(k)}). \label{path}
%\end{equation}
%
%Sample of cross-reference to the formula (\ref{path}) in Appendix \ref{appB}.
%\end{appendix}

%%%%%%%%%%%%%%%%%%%%%%%%%%%%%%%%%%%%%%%%%%%%%%
%% Support information, if any,             %%
%% should be provided in the                %%
%% Acknowledgements section.                %%
% %%%%%%%%%%%%%%%%%%%%%%%%%%%%%%%%%%%%%%%%%%%%%%
\begin{acks}[Acknowledgments]
We thank Diego Martinez Taboada for the insightful comment that helped in removing an unnecessary eigenvalue condition in Theorem~\ref{Theorem: Gaussian Approximation} of the previous manuscript. 
We also thank the referees for their constructive comments that significantly improved this paper. Ilmun Kim acknowledges support from the Yonsei University Research Fund of 2022-22-0289 as well as support from the Basic Science Research Program through the National Research Foundation of Korea (NRF) funded by the Ministry of Education (2022R1A4A1033384) and the Korea government (MSIT) RS-2023-00211073.

 \end{acks}

%%%%%%%%%%%%%%%%%%%%%%%%%%%%%%%%%%%%%%%%%%%%%%
%% Funding information, if any,             %%
%% should be provided in the                %%
%% funding section.                         %%
%%%%%%%%%%%%%%%%%%%%%%%%%%%%%%%%%%%%%%%%%%%%%%
%\begin{funding}
%\end{funding}

%%%%%%%%%%%%%%%%%%%%%%%%%%%%%%%%%%%%%%%%%%%%%%
%% Supplementary Material, including data   %%
%% sets and code, should be provided in     %%
%% {supplement} environment with title      %%
%% and short description. It cannot be      %%
%% available exclusively as external link.  %%
%% All Supplementary Material must be       %%
%% available to the reader on Project       %%
%% Euclid with the published article.       %%
%%%%%%%%%%%%%%%%%%%%%%%%%%%%%%%%%%%%%%%%%%%%%%
%\begin{supplement}
%\stitle{Supplementary Material}
%\sdescription{Enter details.}
%\end{supplement}

\begin{supplement}
   Additional results are provided in the supplementary material. Appendix~\ref{Section: Multiple sample-splitting} discusses multiple sample-splitting, while Appendix~\ref{Section: Further discussions on power} describes a general strategy for studying the asymptotic power of the proposed test. Appendix~\ref{Section: Asymptotic analysis with an increasing conditioning set} gives an example that demonstrates the non-triviality of asymptotic normality under sample splitting. In Appendix~\ref{Section: Dimension-agnostic confidence set}, we give an illustration of how to construct dimension-agnostic confidence sets by inverting dimension-agnostic $p$-values. Appendix~\ref{Section: One-sample MMD and minimax optimality} illustrates our main results using Gaussian MMD and studies minimax power against $L_2$-alternatives. We support our theoretical findings with simulations in Appendix~\ref{Section: Simulations}, and all proofs are provided in Appendix~\ref{Section: Proofs}.
\end{supplement}
%\begin{supplement}
%\stitle{Title of Supplement B}
%\sdescription{Short description of Supplement B.}
%\end{supplement}

%%  imsart-number.bst  will be used to                     %%
%%  MR numbers will be added by VTeX.                      %%
%%                                                         %%
%%  Use \cite{...} to cite references in text.             %%
%%                                                         %%
%%%%%%%%%%%%%%%%%%%%%%%%%%%%%%%%%%%%%%%%%%%%%%%%%%%%%%%%%%%%%

%% if your bibliography is in bibtex format, uncomment commands:
\bibliographystyle{imsart-number} % Style BST file
\bibliography{reference}       % Bibliography file (usually '*.bib')

%\bibliographystyle{apalike}
%\bibliography{reference}

\clearpage 

\appendix

\section{Multiple sample-splitting}  \label{Section: Multiple sample-splitting}
%Recall that our test statistic is based on single splitting where one portion of the data is used to construct a univariate function $f_{\mathrm{mean}}(x)$ and the second portion of the data is used to produce univariate random variables $f_{\mathrm{mean}}(X_1),\ldots, f_{\mathrm{mean}}(X_{m_1})$. To improve statistical efficiency as well as to reduce the randomness from the single split, it is natural to think about multiple sample-splitting where we repeat the previous procedure over different splits and then combine the resulting test statistics in an appropriate way. 

\noindent \textbf{Limiting distribution of a cross-fitted statistic.} From our discussion in Section~\ref{Section: Discussion}, it would be beneficial to consider multiple splits to improve statistical efficiency and also reduce the randomness from a single split. While this is a natural extension of our work, we show that test statistics via multiple sample-splitting do not necessarily have a simple Gaussian limiting distribution, i.e.~they may lose the dimension-agnostic property. To fix ideas, let us focus on one-sample mean testing and  write $\widehat{\mu}_1 = m^{-1} \sum_{i=1}^m X_i$ and $\widehat{\mu}_2 = m^{-1} \sum_{i=m+1}^n X_i$ where $m :=m_1=m_2$. We then define a test statistic via cross-fitting as 
\begin{align*}
	T_{\mathrm{mean}}^{\text{cross}} := \frac{\sqrt{m} \widehat{\mu}_1^\top \widehat{\mu}_2}{\sqrt{m^{-1} \sum_{i=1}^m \big\{ \widehat{\mu}_2^\top (X_i - \widehat{\mu}_1) \big\}^2 }} + \frac{\sqrt{m} \widehat{\mu}_2^\top  \widehat{\mu}_1}{\sqrt{m^{-1} \sum_{i=m+1}^n \big\{ \widehat{\mu}_1^\top (X_i - \widehat{\mu}_2) \big\}^2 }}.
\end{align*}
The limiting null distribution of $T_{\mathrm{mean}}^{\text{cross}}$ under a fixed dimensional setting can be derived based on the multivariate central limit theorem together with the continuous mapping theorem as follows. Throughout the rest of this paper, we say $X_n \convD X$ when $X_n$ converges to $X$ in distribution. 
\begin{proposition}[Limiting null distribution of $T_{\mathrm{mean}}^{\text{cross}}$] \label{Proposition: limiting distribution of cross-splitting}
	Suppose that $X$ has a multivariate distribution in $\mathbb{R}^d$ with $\mu=0$ and a positive definite covariance matrix $\Sigma$. Suppose that the distribution of $X$ does not change with $n$, i.e. the dimension is fixed. Then
	\begin{align*}
		T_{\mathrm{mean}}^{\mathrm{cross}} \convD \frac{\boldsymbol{\xi}_1^\top \Sigma \boldsymbol{\xi}_2}{\sqrt{\boldsymbol{\xi}_1^\top \Sigma^2 \boldsymbol{\xi}_1}} + \frac{\boldsymbol{\xi}_1^\top \Sigma \boldsymbol{\xi}_2}{\sqrt{\boldsymbol{\xi}_2^\top \Sigma^2 \boldsymbol{\xi}_2}},
	\end{align*}
	where $\boldsymbol{\xi}_1, \boldsymbol{\xi}_2 \overset{i.i.d}{\sim} N(0,I)$. 
\end{proposition}
It is worth pointing out that while both $\boldsymbol{\xi}_1^\top \Sigma \boldsymbol{\xi}_2 / \{\boldsymbol{\xi}_1^\top \Sigma^2 \boldsymbol{\xi}_1\}^{1/2}$ and $\boldsymbol{\xi}_1^\top \Sigma \boldsymbol{\xi}_2 / \{\boldsymbol{\xi}_2^\top \Sigma^2 \boldsymbol{\xi}_2\}^{1/2}$ have the standard normal distribution, the distribution of the sum of these two random variables is not Gaussian. We only expect that considering many splits makes the limiting distribution more complicated, which demonstrates challenges associated with multiple splitting in terms of calibration. Alternatively, one can combine different $p$-values from different single splits using Bonferroni correction \citep[see, e.g.][for other methods of combining $p$-values]{vovk2018combining}. However the resulting test is expected to be conservative as it does not use the dependence structure between $p$-values.

\vskip 1em

\noindent \textbf{Calibration under symmetry condition.} 
On the other hand, it is relatively straightforward to calibrate the test statistic using many splits under the symmetry condition. One key observation is that, when $X$ and $-X$ have the same distribution, any test statistic, let's say $T(X_1,\ldots,X_n)$, has the same distribution as $T(\varepsilon_1 X_1,\ldots, \varepsilon_n X_n)$ where $\{\varepsilon_i\}_{i=1}^n$ are i.i.d.~Rademacher random variables independent of the data. We are then able to use exchangeability among different test statistics associated with different Rademacher random variables and produce a valid level $\alpha$ test as in (\ref{Eq: Monte-Carlo test}). We note that, similar to the permutation procedure for two-sample testing or independence testing, this flip-sign procedure works as a wrapper around any statistic including the average of statistics from multiple splits and provides a valid $p$-value in finite samples. %We note that, unlike $\phi_{\mathrm{sym}}^{\dagger}$ where we multiplied Randemacher random variables to the half of the data, the sign-flipping procedure for multiple splits should be applied to the entire observations to guarantee the validity. 

Another way, which is more computationally efficient but more conservative, is to use the concentration bound as before in Corollary~\ref{Corollary: properties of phi_sym}. Let $T_{\mathrm{mean},(1)},\ldots,T_{\mathrm{mean},(K)}$ be the test statistic in (\ref{Eq: studentized statistic for mean testing}) computed over different splits. We denote $V_{(i)}:= \sqrt{m_1} T_{\mathrm{mean},(i)} / \big\{ m_1 + T^2_{\mathrm{mean},(i)} \big\}^{1/2}$ for $i=1,\ldots,K$ and define a test
\begin{align*}
	\phi_{\mathrm{sym},K} := \mathds{1} \Bigg\{ \frac{1}{K} \sum_{i=1}^K e^{\lambda_i V_{(i)} - \lambda_i^2/2} > \frac{1}{\alpha}\Bigg\},
\end{align*}
for some $\lambda_i \in \mathbb{R}$ and $i=1,\ldots,K$. Then Markov's inequality shows that $\sup_{P \in \mathcal{P}_{\mathrm{sym},n,d}}\mE_{P} [\phi_{\mathrm{sym},K} ] \leq \alpha$ since $\mE_P[\exp(\lambda_i V_{(i)} - \lambda_i^2/2)] \leq 1$ for each $i=1,\ldots,K$ under the symmetry assumption. We note that a similar idea appears in \cite{wasserman2019universal} where they study likelihood ratio tests based on sample-splitting. %Although it is intuitively clear that 

%\paragraph{Practical recommendation.}

\section{Further discussions on power} \label{Section: Further discussions on power}

This section provides a couple of results on the power of the proposed test. In Appendix~\ref{Section: Sufficient conditions for uniform power}, we present sufficient conditions that guarantee that the proposed test is uniformly consistent. In Appendix~\ref{Eq: Asymptotic power expression}, we derive the asymptotic power expression of the proposed test in a fixed-dimensional setting and compare it with that of the corresponding U-statistic.

\subsection{Sufficient conditions for uniform power} \label{Section: Sufficient conditions for uniform power}
We start by identifying sufficient conditions under which the test based on the proposed cross U-statistic is asymptotically consistent. All of the uniform consistency results in this paper (except Theorem~\ref{Theorem: Uniform type I and II error control for sparse alternatives} on the maximum-type statistic) can be proved by verifying these sufficient conditions.

Let us first recall a general setting considered in Section~\ref{Section: General results for degenerate U-statistics}. Given $\{X_i\}_{i=1}^n \overset{i.i.d.}{\sim} P$ and a symmetric kernel $h(x,y)$, we reject the null hypothesis of $\mE_P[h(X_1,X_2)] = 0$ if 
\begin{align} \label{Eq: general test}
	\overline{f}_h > z_{1-\alpha} \sqrt{ \frac{1}{m_1^2} \sum_{i=1}^{m_1} \{ f_h(X_i) - \overline{f}_h \}^2},
\end{align}
where $f_h(x)$ is defined in (\ref{Eq: general linear function}) computed over $\{X_{m_1+1},\ldots,X_n\}$ and $\overline{f}_h = m_1^{-1}\sum_{i=1}^{m_1} f_h(X_i)$. Here, $m_1$ and $m_2$ are some positive integers such that $m_1+m_2=n$.

Consider a positive sequence $\psi_n$ that goes to infinity at any rate as $\min\{m_1,m_2\} \rightarrow \infty$. Given the sequence $\psi_n$, we denote by $\mathcal{Q}_1:= \mathcal{Q}_1(\psi_n)$ a class of alternative distributions such that for each $P \in \mathcal{Q}_1$ and each $m_1,m_2 \in \mathbb{N}$,
\begin{equation}
\begin{aligned}  \label{Eq: Sufficient Condition for Power}
	\mE_P[h(X_1,X_2)] \geq \psi_n \max\{z_{1-\alpha}, 1\} \max \Bigg\{ & \sqrt{\frac{ \mE_P[\mV_P[h(X_1,X_2)|X_1]]}{m_1m_2}}, \\
	& \sqrt{\biggl( \frac{1}{m_1} + \frac{1}{m_2} \biggr) \mV_P[\mE_P[h(X_1,X_2)|X_1]]} \Bigg\}.
\end{aligned}
\end{equation}
The next lemma shows that our test is uniformly consistent against the given class of alternatives $\mathcal{Q}_1$.

\begin{lemma}[Sufficient conditions for uniform type II error control] \label{Lemma: Sufficient conditions for uniform type II error control}
	For $\alpha \in (0,1)$, suppose that the condition~(\ref{Eq: Sufficient Condition for Power}) holds uniformly over $P \in \mathcal{Q}_1$. Then the uniform type II error of the test~(\ref{Eq: general test}) converges to zero:
	\begin{align*}
		\lim_{\min\{m_1,m_2\} \rightarrow \infty }\sup_{P \in \mathcal{Q}_1} P \Biggl( \overline{f}_h \leq  z_{1-\alpha} \sqrt{ \frac{1}{m_1^2} \sum_{i=1}^{m_1} \{ f_h(X_i) - \overline{f}_h \}^2} \Biggr) = 0.
	\end{align*}
\end{lemma}

\vskip 1em

\noindent \textbf{Remark.}
\begin{itemize}
	\item In Lemma~\ref{Lemma: Sufficient conditions for uniform type II error control}, we allow $\alpha \in (0,1)$ to change with $n$. Indeed $\alpha$ can decrease to zero as long as the sufficient condition~(\ref{Eq: Sufficient Condition for Power}) holds.
	\item By Jensen's inequality, it holds that 
	\begin{align*}
		\mE_P[\mV_P[h(X_1,X_2)|X_1]] \leq \mE_P[h^2(X_1,X_2)],
	\end{align*}
	based on which the sufficient condition~(\ref{Eq: Sufficient Condition for Power}) can be further simplified as
	\begin{align}  \label{Eq: Sufficient Condition for Power II}
			\mE_P[h(X_1,X_2)] ~\geq~  & \psi_n \max\{z_{1-\alpha}, 1\} \\ \nonumber
			& \times  \max \Bigg\{ \sqrt{\frac{ \mE_P[h^2(X_1,X_2)]}{m_1m_2}}, 
			\sqrt{\biggl( \frac{1}{m_1} + \frac{1}{m_2} \biggr) \mV_P[\mE_P[h(X_1,X_2)|X_1]]} \Bigg\}.
	\end{align}
	\item The proof of this result can be found in Appendix~\ref{Section: Proof of Lemma: Sufficient conditions for uniform type II error control}.
\end{itemize}

\subsection{Asymptotic power expression when $d=1$} \label{Eq: Asymptotic power expression}
In Theorem~\ref{Theorem: asymptotic power} and Theorem~\ref{Theorem: asymptotic power of covariance test}, we demonstrate that the asymptotic power of the proposed test has the form of a normal cumulative distribution function in the high-dimensional regime where $d/n \rightarrow \tau \in (0,\infty)$. These results allow us to show that the asymptotic power of the proposed test is only worse by a factor of $\sqrt{2}$ than that of the corresponding U-statistic for mean and covariance testing. However, outside of this high-dimensional regime, the limiting distributions of the proposed statistic as well as the original U-statistic are not necessarily Gaussian under alternatives, and thus it is difficult to figure out the relative efficiency in power without relying on numerical methods.

In this subsection, we focus on one-sample mean testing and demonstrate this difficulty. In particular, we derive the asymptotic power expressions of both $T_{\mathrm{mean}}$ and $U_{\mathrm{mean}}$ for the one-dimensional case, i.e.~$d=1$ so that $d/n \rightarrow 0$. For simplicity, assume that $X_1,\ldots,X_{2n} \overset{i.i.d.}{\sim} N(hn^{-1/2},1)$ for some fixed $h \neq 0$. Let us write 
\begin{align*}
	\widehat{\mu}_1 = \frac{1}{n} \sum_{i=1}^n X_i \quad \text{and} \quad \widehat{\mu}_2 = \frac{1}{n} \sum_{i=n+1}^{2n} X_i.
\end{align*} 
Then the proposed test statistic has the form of 
\begin{align*}
	T_{\mathrm{mean}} = \mathrm{sign}(\widehat{\mu}_2) \frac{\sqrt{n} \widehat{\mu}_1}{\sqrt{\frac{1}{n}\sum_{i=1}^n(X_i - \widehat{\mu}_1)^2}},
\end{align*}
where $\mathrm{sign}(x)$ denotes the sign function of $x$. Using this expression, the resulting test $T_{\mathrm{mean}} > z_{1-\alpha}$ can be shown to have the asymptotic power against this local alternative, i.e.~$\mu_n = hn^{-1/2}$, as
\begin{align*}
	P(T_{\mathrm{mean}} > z_{1-\alpha}) ~=~ & P(\mathrm{sign}(\widehat{\mu}_2) = +1) \times P\Biggl( \frac{\sqrt{n} \widehat{\mu}_1}{\sqrt{\frac{1}{n}\sum_{i=1}^n(X_i - \widehat{\mu}_1)^2}} > z_{1-\alpha} \Biggr) \\[.5em]
	 +~ & P(\mathrm{sign}(\widehat{\mu}_2) = -1) \times P\Biggl( \frac{\sqrt{n} \widehat{\mu}_1}{\sqrt{\frac{1}{n}\sum_{i=1}^n(X_i - \widehat{\mu}_1)^2}} < - z_{1-\alpha} \Biggr) \\[.5em]
	 = ~ & \Phi(h) \times \Phi(z_\alpha + h) +  \Phi(-h) \times \Phi(z_\alpha - h) +o(1).
\end{align*}
On the other hand, $nU_{\mathrm{mean}} + 1$ converges to a chi-square distribution with one degree of freedom under the null, and it converges to a non-central chi-square distribution with one degree of freedom and non-centrality parameter $h^2$ under the local alternative. In other words, the resulting power functions are defined in terms of different distribution functions, and they need to be compared numerically.

\section{A failure of standard CLT under sample splitting} \label{Section: Asymptotic analysis with an increasing conditioning set}

Suppose that $X_1,\ldots,X_n$ are uniformly distributed over $[0,1]$, and we let $\widehat{g}^\ast(x) = \mathds{1}(x \leq m_2^{-1})$ where $m_2 = n - m_1$. Then $\widehat{g}^\ast(X_1),\ldots,\widehat{g}^\ast(X_{m_1})$ are i.i.d.~Bernoulli random variables with success probability $m_2^{-1}$. Suppose further that $\widehat{g}^\ast$ is fixed (thereby $m_2$ is fixed as well). Under this setting, the standard central limit theorem along with Slutsky's theorem shows that 
\begin{align*}
	T := \frac{\sqrt{m_1} \bigl( \frac{1}{m_1} \sum_{i=1}^{m_1} \widehat{g}^\ast(X_i) - \mE[\widehat{g}^\ast(X_i) \,|\, \widehat{g}^\ast] \bigr)}{\sqrt{ \frac{1}{m_1} \sum_{i=1}^{m_1} \bigl\{\widehat{g}^\ast(X_i) - \frac{1}{m_1} \sum_{j=1}^{m_1} \widehat{g}^\ast(X_j) \bigr\}^2}}  \convD N(0,1),
\end{align*}
as $n \rightarrow \infty$. However, if $m_2$ increases with $n$, the asymptotic normality is not guaranteed. 

To fix ideas, we assume that $m_1 = m_2$. Under this scenario, the probability of observing $\widehat{g}^\ast(X_1) = \ldots = \widehat{g}^\ast(X_{m_1}) = 0$ is non-vanishing in the limit. We simply let $T = 0$ whenever this event occurs. Hence $T$ can be written as
\begin{align*}
	T ~=~ &  \frac{\sum_{i=1}^{m_1} \widehat{g}^\ast(X_i)  - 1}{\sqrt{\sum_{i=1}^{m_1} \widehat{g}^\ast(X_i)}} \times \mathds{1}\biggl( \sum_{i=1}^{m_1} \widehat{g}^\ast(X_i)  \neq 0 \biggr) \times \sqrt{\frac{1}{1 - m_1^{-1} \sum_{i=1}^{m_1} \widehat{g}^\ast(X_i)}}   \\[.5em]
	  =~&  \frac{\sum_{i=1}^{m_1} \widehat{g}^\ast(X_i)  - 1}{\sqrt{\sum_{i=1}^{m_1} \widehat{g}^\ast(X_i)}} \times \mathds{1}\biggl( \sum_{i=1}^{m_1} \widehat{g}^\ast(X_i)  \neq 0 \biggr) \times \{1 + o_P(1)\},
\end{align*}
where the second equality follows by the weak law of large numbers and the continuous mapping theorem. By defining a continuous function $h$ as 
\begin{align*}
	h(x) = \begin{cases}
		\frac{x-1}{\sqrt{x}} & \text{if $x \geq 1$,} \\
		0 & \text{if $x < 1$,}
	\end{cases}
\end{align*}
we have 
\begin{align*}
	h\biggl( \sum_{i=1}^{m_1} \widehat{g}^\ast(X_i)  \biggr)  = \frac{\sum_{i=1}^{m_1} \widehat{g}^\ast(X_i) - 1}{\sqrt{\sum_{i=1}^{m_1} \widehat{g}^\ast(X_i)}} \times \mathds{1}\biggl( \sum_{i=1}^{m_1} \widehat{g}^\ast(X_i)  \neq 0 \biggr). 
\end{align*}
Given this identity and letting $W \sim \mathrm{Poisson}(1)$, the Poisson limit theorem along with the continuous mapping theorem yields
\begin{align*}
	h\biggl( \sum_{i=1}^{m_1} \widehat{g}^\ast(X_i)  \biggr) \convD h(W).
\end{align*}
Finally, Slutsky's theorem establishes that
\begin{align*}
	T \convD \frac{W - 1}{\sqrt{W}} \mathds{1}(W \neq 0),
\end{align*}
and thus the limiting distribution is not Gaussian.

While the above example is not practically oriented given that $\widehat{g}^\ast$ only uses the cardinality of $\mathcal{X}_2$, it clearly reveals that the standard central limit theorem is not directly applicable with an increasing conditioning set. 

\section{Dimension-agnostic confidence sets} \label{Section: Dimension-agnostic confidence set}

%So far we have proposed dimension-agnostic tests for one-sample mean testing. 
This section leverages the previous results in Section~\ref{Section: Mean Testing} and constructs dimension-agnostic confidence sets for $\mu \in \mathbb{R}^d$ in the sense of (\ref{eq:formal-goal2}). We achieve this goal by the duality of hypothesis tests and confidence sets. 

\begin{figure}[t!]
	\begin{center}
		\includegraphics[width=1.00\textwidth]{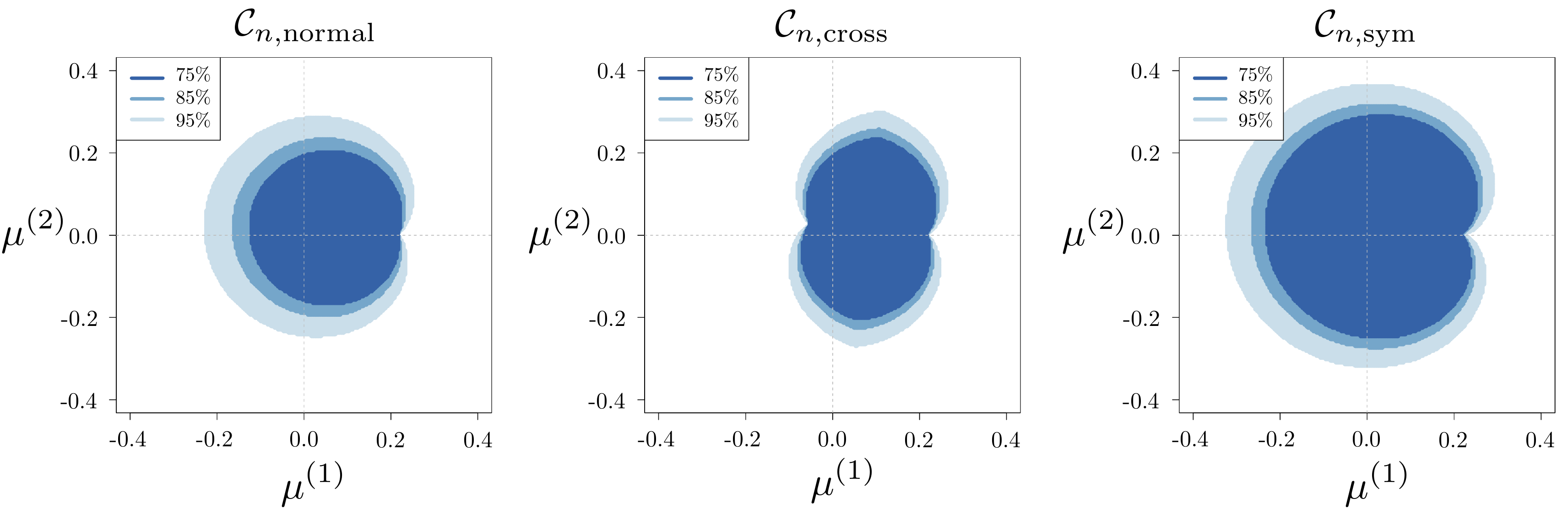}
	\end{center}
	\caption{Illustration of confidence sets $\mathcal{C}_{n,\text{normal}}(\mathcal{X}_1, \mathcal{X}_2;\alpha)$ (left), $\mathcal{C}_{n,\text{cross}}(\mathcal{X}_1, \mathcal{X}_2;\alpha)$ (middle) and $\mathcal{C}_{n,\mathrm{sym}}(\mathcal{X}_1, \mathcal{X}_2;\alpha)$ (right) for $\mu=(\mu^{(1)},\mu^{(2)})^\top$ at $\alpha=0.05$, $0.15$ and $0.25$. These confidence sets are (asymptotically) valid independent of $d$. We focus on the bivariate case for visualization only.} \label{Figure: CI}
\end{figure}

Given i.i.d.~observations $\mathcal{X}_1 = \{X_i\}_{i=1}^{m_1}$ and $\mathcal{X}_2= \{X_i\}_{i=m_1+1}^n$ from a distribution with mean $\mu \in \mathbb{R}^d$, we denote by $T_{\mathrm{mean},\mu}$ the studentized statistic~(\ref{Eq: studentized statistic for mean testing}) computed based on the centered samples $X_1 -\mu,\ldots,X_n - \mu$. Then Theorem~\ref{Theorem: unconditional BE Bound} verifies that $T_{\mathrm{mean},\mu}$ has a Gaussian limiting distribution under Assumption~\ref{assumption: moment condition}. This is true for any $\mu \in \mathbb{R}^d$, which implies the confidence set defined below is asymptotically valid:
\begin{align*}
	\mathcal{C}_{n,\text{normal}}(\mathcal{X}_1, \mathcal{X}_2;\alpha) := \big\{\mu \in \mathbb{R}^d : T_{\mathrm{mean},\mu} < z_{1-\alpha} \big\}.
\end{align*}
By the nature of the test statistic, the shape of $\mathcal{C}_{n,\text{normal}}(\mathcal{X}_1, \mathcal{X}_2;\alpha)$ tends to be asymmetric (see Figure~\ref{Figure: CI}). To fix this issue, we also consider a cross-fitted confidence set with the Bonferoni correction defined as
\begin{align*}
	\mathcal{C}_{n,\text{cross}}(\mathcal{X}_1, \mathcal{X}_2;\alpha) := \mathcal{C}_{n,\text{normal}}(\mathcal{X}_1, \mathcal{X}_2;\alpha/2) \cap \mathcal{C}_{n,\text{normal}}(\mathcal{X}_2, \mathcal{X}_1;\alpha/2). 
\end{align*}
In addition, motivated by our discussion in Section~\ref{Section: Non-asymptotic calibration under symmetry}, we can construct a confidence set
\begin{align*}
	\mathcal{C}_{n,\mathrm{sym}}(\mathcal{X}_1, \mathcal{X}_2;\alpha) := \Bigg\{  \mu \in \mathbb{R}^d : T_{\mathrm{mean},\mu}  <  \sqrt{\frac{2 \log (1/\alpha)}{1 - 2 m_1^{-1}\log(1/\alpha)}} \Bigg\},
\end{align*}
which has a finite sample coverage guarantee when $m_1 > 2 \log (1/\alpha)$ and $X - \mu$ is symmetric at the origin. We formalize these results in the following proposition. 

\begin{proposition}[Confidence sets]
	Let $\mathcal{P}_{\mu,d}$ be the class of distributions that satisfies Assumption~\ref{assumption: moment condition}. Then the confidence sets $\mathcal{C}_{n,\mathrm{normal}}$ and  $\mathcal{C}_{n,\mathrm{cross}}$ are dimension-agnostic in the sense of (\ref{eq:formal-goal2}). In particular, we have 
	\begin{align*}
		\lim_{m_1 \rightarrow \infty} \inf_{P \in \mathcal{P}_{\mu,d}} P\{ \mu \in \mathcal{C}_{n,\mathrm{normal}}(\mathcal{X}_1, \mathcal{X}_2;\alpha)\}& = 1 - \alpha,\\
		\text{and ~}~ \lim_{m_1,m_2 \rightarrow \infty} \inf_{P \in \mathcal{P}_{\mu,d}} P\{ \mu \in \mathcal{C}_{n,\mathrm{cross}}(\mathcal{X}_1, \mathcal{X}_2;\alpha)\} & \geq 1 - \alpha.
	\end{align*} 
	Let $\mathcal{P}_{\mathrm{sym},\mu,d}$ be the collection of distributions that are symmetric at $\mu$. Then, provided that $m_1 > 2 \log (1/\alpha)$, the confidence set $\mathcal{C}_{n,\mathrm{sym}}(\mathcal{X}_1, \mathcal{X}_2;\alpha)$ has a finite sample guarantee as 
	\[
	\inf_{P \in \mathcal{P}_{\mathrm{sym},\mu,d}} P\{ \mu \in \mathcal{C}_{n,\mathrm{normal}}(\mathcal{X}_1, \mathcal{X}_2;\alpha)\} \geq 1 - \alpha.
	\] 
\end{proposition}
The result follows directly by the duality of hypothesis tests and confidence sets. Hence we omit the proof. In Figure~\ref{Figure: CI}, we visually demonstrate the constructed confidence sets at $\alpha = 0.05, 0.15$ and $0.25$ by generating $n=100$ independent samples from $N(\mu,\Sigma)$ where $\mu = (0,0)^\top$ and $\Sigma = I$. They have all reasonable shapes containing the true mean vector in this particular example and we observed similar results using different random seeds as well. From Figure~\ref{Figure: CI}, we note that neither $\mathcal{C}_{n,\mathrm{normal}}$ nor $\mathcal{C}_{n,\mathrm{cross}}$ dominates the other, and not surprisingly, $\mathcal{C}_{n,\mathrm{sym}}(\mathcal{X}_1, \mathcal{X}_2;\alpha)$ based on a non-asymptotic concentration bound has the largest area among the considered confidence sets. Next, building on the results in this section, we provide a remark on testing for composite nulls. 
\begin{remark}[Composite null in $\mu$] \normalfont
	Once we have a valid confidence set for $\mu$, we can use it to test a general composite null $H_0: \mu \in S$, where $S$ is an arbitrary set in $\mathbb{R}^d$. In particular, $p_n:=\inf\{ \alpha: \mathcal{C}_{n,\text{normal}}(\mathcal{X}_1,\mathcal{X}_2;\alpha) \cap S = \emptyset \}$ is a dimension-agnostic $p$-value in the sense of (\ref{eq:formal-goal}) and, as a result, $\mathds{1}(p_n \leq \alpha)$ is a dimension-agnostic level-alpha test for $H_0: \mu \in S$. Similarly, we can also use the other derived confidence sets to test for the same composite null. 
	We are not aware of any other tests or confidence sets with such a dimension-agnostic property.
\end{remark}

\section{One-sample maximum mean discrepancy (MMD)}  \label{Section: One-sample MMD and minimax optimality}
Here we illustrate Theorem~\ref{Theorem: Gaussian Approximation} in the context of nonparametric goodness-of-fit testing. Suppose we are interested in testing whether $H_0: P=Q$ versus $H_1: P\neq Q$ where $P$ is a data generating distribution and $Q$ is some theoretical distribution. One method for this problem is based on the Gaussian MMD. We denote a Gaussian kernel with a scaling parameter $\nu_n >0$ by 
\begin{align} \label{Eq: Gaussian kernel0}
	G_{\nu_n}(x,y) := \exp(-\nu_n \|x - y\|^2_d), \quad \forall x, y \in \mathbb{R}^d.
\end{align}
Then the Gaussian MMD is defined as
\begin{align*}
	\text{MMD}(P,Q) := \mE_P[G_{\nu_n}(X_1,X_2)] + \mE_{Q}[G_{\nu_n}(Y_1,Y_2)] - 2  \mE_{P,Q}[G_{\nu_n}(X_1,Y_1)], 
\end{align*}
where $\{Y_1,Y_2\} \overset{i.i.d.}{\sim} Q$. Given a sample $\{X_i\}_{i=1}^n \overset{i.i.d.}{\sim} P$, $\text{MMD}(P,Q)$ can be estimated using the U-statistic~(\ref{Eq: U-statistic}) with the kernel
\begin{align} \label{Eq: Gaussian kernel}
	h_{\text{Gau},\nu_n}(x,y) := G_{\nu_n}(x,y) - \mE_Q[G_{\nu_n}(Y_1,y)] - \mE_Q[G_{\nu_n}(x,Y_2)] + \mE_Q[G_{\nu_n}(Y_1,Y_2)].
\end{align}
When $Q$ is a multivariate Gaussian, a closed-form expression of $h_{\text{Gau},\nu_n}$ is known \citep{kellner2019one,makigusa2020asymptotics}. However the resulting U-statistic is degenerate under the null and thus consistent estimation of its asymptotic distribution can be challenging. On the other hand our test statistic~(\ref{Eq: studentized statistic with h}) with kernel $h_{\text{Gau},\nu_n}$ has a simple Gaussian limiting distribution under the conditions in Theorem~\ref{Theorem: Gaussian Approximation}. 

%\subsection{Minimax rate optimality against $L_2$ alternatives}

In addition to the dimension-agnostic null for our modified one-sample MMD, we can build on the recent result of \cite{li2019optimality} to show that the resulting test is minimax rate optimal in the $L_2$ distance and hence its power cannot be improved beyond a constant factor. For the purposes of examining power, it is standard in nonparametric testing (such as the aforementioned work) to assume that the dimension of the data $d$ is fixed in the sample size but allow the intrinsic dimension $D_n$~(\ref{Eq: intrinsic dimension}) to increase to infinity, recalling that it is the latter that determines limiting behavior. To describe the result, let us denote our test statistic~(\ref{Eq: studentized statistic with h}) with kernel $h_{\text{Gau},\nu_n}$ by $T_{\text{MMD}}$ and define a test $$\phi_{\mathrm{MMD}} := \mathds{1}(T_{\text{MMD}} > z_{1-\alpha}).$$ In addition, following the notation in \cite{li2019optimality}, we denote the $s$th order Sobolev space in $\mathbb{R}^d$ by
\begin{align*}
	\mathcal{W}_{d}^{s,2} := \Big\{ f : \mathbb{R}^d \mapsto \mathbb{R} ~ \big| ~ \text{$f$ is continuous and } \int (1 + \|\omega \|^2)^{s/2} \| \mathcal{F}(g)(\omega)\|^2 d\omega < \infty  \Big\},
\end{align*}
where $\|\cdot\|$ denotes the Euclidean norm in $\mathbb{R}^d$ and $\mathcal{F}(g)$ is the Fourier transform of $g$:
\begin{align*}
	\mathcal{F}(g)(\omega) := \frac{1}{(2\pi)^{d/2}} \int_{\mathbb{R}^d} g(x) e^{-i x^\top \omega} dx.
\end{align*}
We denote the Sobolev norm of $g$ by $\| g \|^2_{\mathcal{W}^{s,2}} =  \int_{\mathbb{R}^d}  (1 + \|\omega \|^2)^{s/2} \| \mathcal{F}(g)(\omega)\|^2 d\omega$ and define
\begin{align*}
	\mathcal{W}^{s,2} (M) := \big\{g \in \mathcal{W}^{s,2}: \| g \|_{\mathcal{W}^{s,2}} < M  \big\}.
\end{align*}
Finally, for a given sequence of $\epsilon_n > 0$, we define a class of alternatives 
\begin{align*}
	\mathcal{Q}_{\text{density}}(\epsilon_n):= \big\{ p \in \mathcal{W}^{s,2} (M) : \| p - q\|_{L_2} \geq \epsilon_n \big\},
\end{align*}
where $p$ and $q$ are density functions of $P$ and $Q$, respectively, and $\|g\|_{L_2}^2 = \int_{\mathbb{R}^d} g^2(x) dx$. We show that the uniform type II error of $\phi_{\mathrm{MMD}}$ goes to zero as long as $\epsilon_n n^{2s/(d+4s)} \rightarrow \infty$. 
\begin{proposition}[Type I and II error control] \label{Proposition: MMD test}
	For fixed $\alpha \in (0,1)$, suppose that $\epsilon_n \geq K_n n^{-2s/(d+4s)}$ where $K_n$ is a positive sequence increasing to infinity at any rate. Take $\nu_n \asymp n^{4/(d+4s)}$, which results in $D_n \rightarrow \infty$. Then for $P,Q \in \mathcal{W}^{s,2} (M)$, 	
	we have 
	\begin{align*}
		\lim_{n \rightarrow \infty} \mE_{Q} [\phi_{\mathrm{MMD}}] = \alpha \quad \text{and} \quad \lim_{n \rightarrow \infty} \sup_{P \in \mathcal{Q}_{\mathrm{density}}(\epsilon_n)} \mE_{P} [1 -  \phi_{\mathrm{MMD}}] =0.
	\end{align*}
\end{proposition} 
%We note that the choice of scaling parameter $\nu_n \asymp n^{4/(d+4s)}$ is crucial for controlling the type II error but not the type I error. 
Note that the smoothness parameter $s$ is assumed to be known in the above result. We believe that it is possible to remove this assumption at the expense of a $\log\log n$ factor, in view of \cite{ingster2000adaptive} and \cite{li2019optimality}, and we leave this direction to future work. As proved in \cite{li2019optimality}, no test can be uniformly powerful when $\epsilon_n$ is of order smaller than $n^{-2s/(d+4s)}$. In fact, this is a standard minimax rate for nonparametric testing in the $L_2$ distance~\cite[e.g.][]{ingster1982minimax}. This means that $\phi_{\mathrm{MMD}}$ is minimax rate optimal in the $L_2$ distance under the considered Sobolev smoothness. Thus, our test has a good power property even in a nonparametric setting, while possessing a simple Gaussian null distribution irrespective of the intrinsic dimension of the data.

\section{Numerical simulations} \label{Section: Simulations}

In this section, we present numerical experiments that validate our theoretical findings. In the first experiment, we compare the empirical power of $\phi_{\mathrm{mean}}$ and $\phi_{\mathrm{cov}}$ with their theoretical asymptotic power derived in Theorem~\ref{Theorem: asymptotic power} and Theorem~\ref{Theorem: asymptotic power of covariance test}, respectively. In the second experiment, we illustrate the dimension-agnostic property of the proposed methods and contrast them with the corresponding U-statistics.

\subsection{Power comparisons: U-statistics versus sample-splitting}
\begin{figure}[t!]
	\begin{center}		
		\begin{minipage}[b]{0.46\textwidth}
			\includegraphics[width=\textwidth]{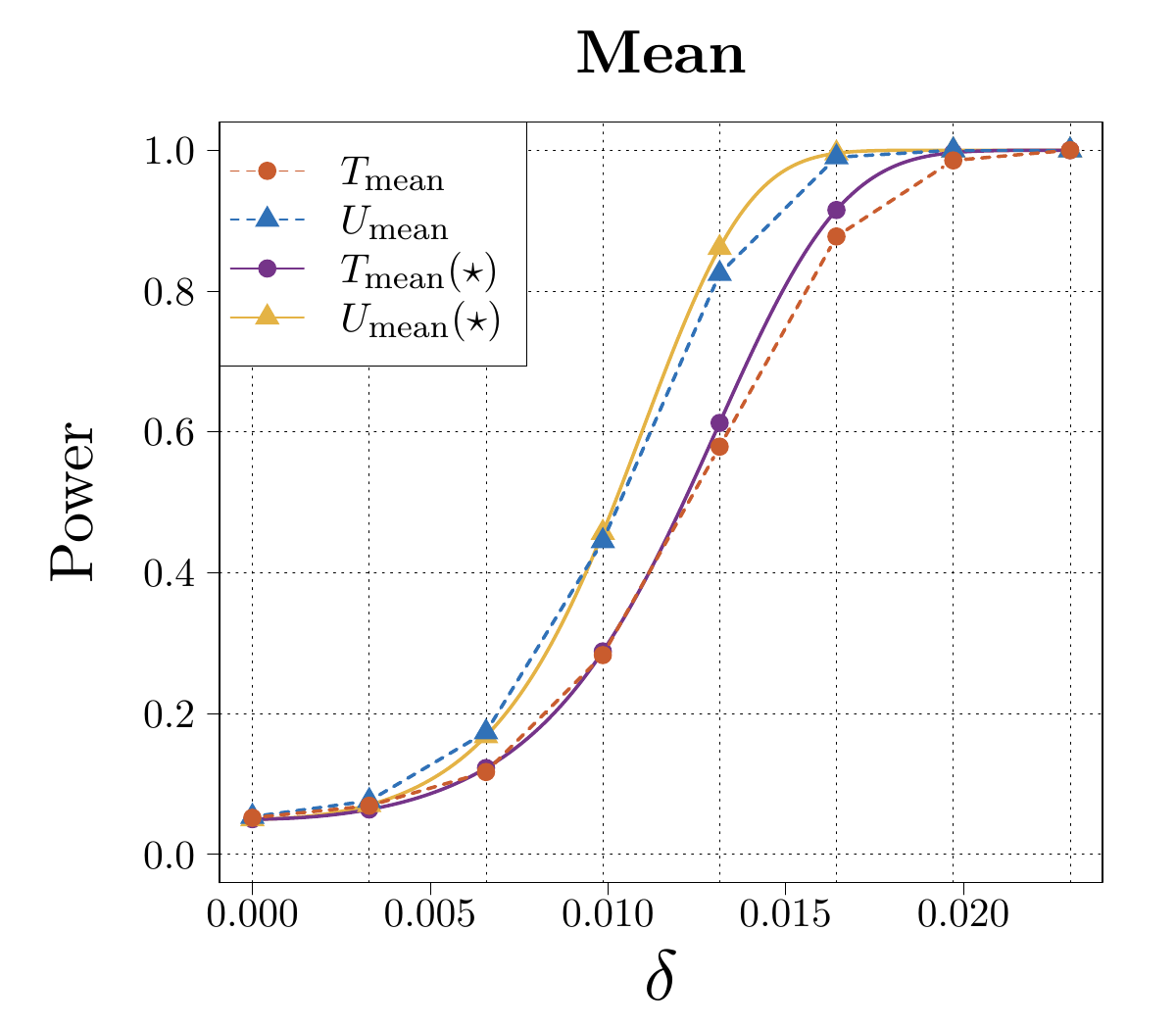}
		\end{minipage} 
		\hskip 2em
		\begin{minipage}[b]{0.46\textwidth}
			\includegraphics[width=\textwidth]{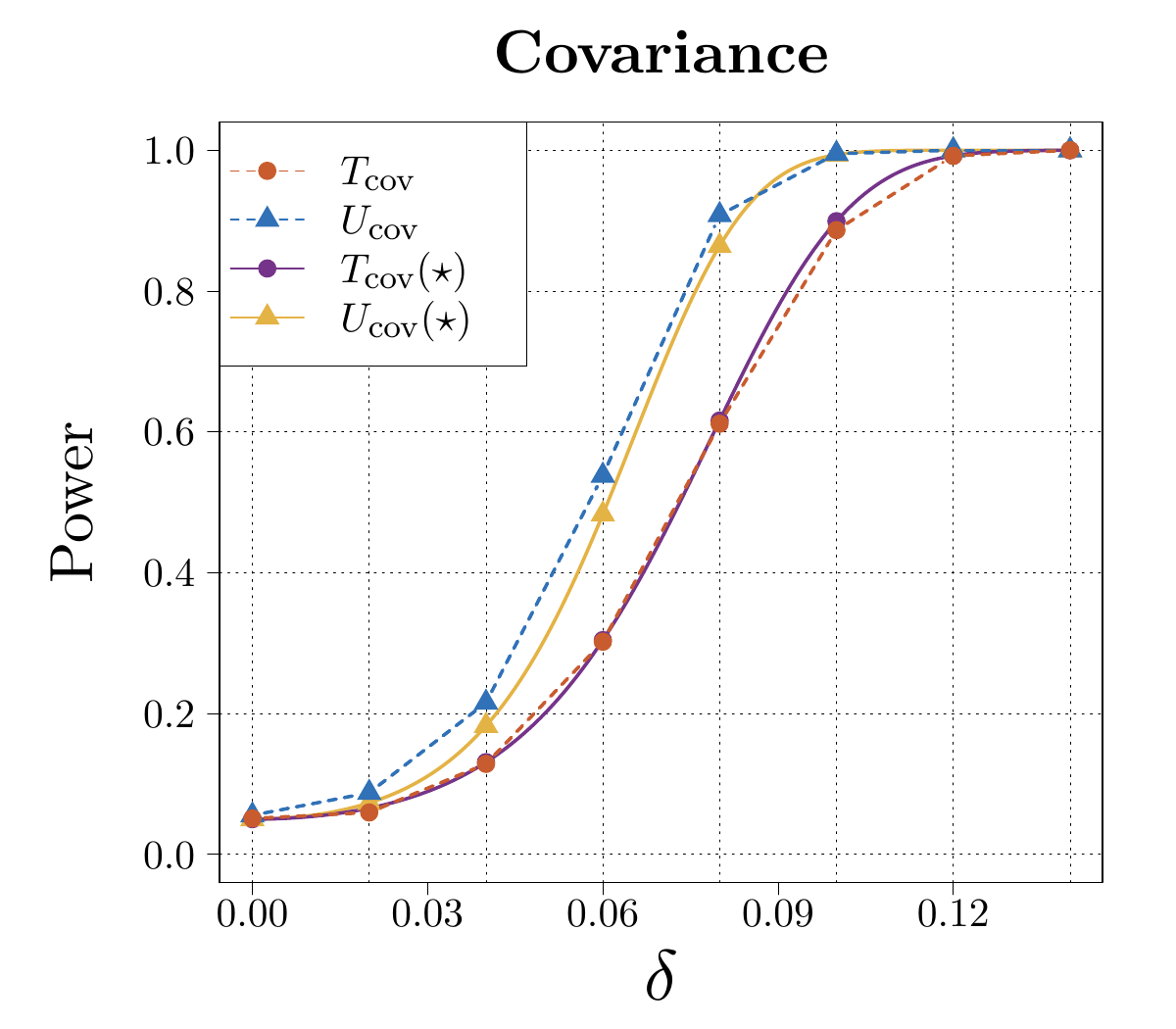}
		\end{minipage}
		\caption{\small Comparisons of the empirical power to the theoretical (asymptotic) power of the considered tests where $(\star)$ denotes the theoretical power of the corresponding test. Left panel: empirical power of $\phi_{\mathrm{mean}}$ closely tracks (\ref{Eq: theoretical power of phi_mean}). Right panel: empirical power of $\phi_{\mathrm{cov}}$ closely tracks (\ref{Eq: theoretical power of phi_cov}). Each plot also contains the empirical and theoretical power of $U_{\mathrm{mean}}$ and $U_{\mathrm{cov}}$, presented in (\ref{Eq: theoretical power of U_1}) and (\ref{Eq: theoretical power of U_2}) respectively.} \label{Figure: Power}
	\end{center}
\end{figure}

The purpose of this subsection is to confirm the validity of the theoretical power expressions for $\phi_{\mathrm{mean}}$ and $\phi_{\mathrm{cov}}$ by simulations. For both tests, we set the sample splitting ratio to be $1/2$ that achieves the asymptotic maximum power presented in (\ref{Eq: theoretical power of phi_mean}) and (\ref{Eq: theoretical power of phi_cov}), respectively. In this numerical experiment, we draw a sample of size $n=1000$ from a multivariate Gaussian with the mean vector $\mu$ and the covariance matrix $\Sigma$ where the dimension is chosen to be $d=500$. For mean testing, we let the mean vector be $\mu = (\delta,\ldots,\delta)^\top$ with $\delta$ ranging from $0$ to $0.023$ and the covariance matrix $\Sigma$ be the identity matrix. Similarly, for covariance testing, we let $\mu = (0,\ldots,0)^\top$ and $\Sigma = \text{diag}\{1+\delta,\ldots,1+\delta\}$ with $\delta$ ranging from $0$ to $0.14$. Note that for both cases, $\delta=0$ corresponds to the null where we expect the power becomes the nominal level $\alpha$. The finite-sample power of each test is approximated based on 2000 Monte Carlo replicates at the significance level $\alpha = 0.05$ and the results can be found in Figure~\ref{Figure: Power}.

From the results, we find that the empirical power is indeed well approximated by the power expression predicted by our theory, especially when $\delta$ is relatively small. A minor gap in large $\delta$ is expected since our theory holds under local alternatives where $\mu^\top \mu$ is of order smaller than or comparable to $\sqrt{d}/n$. We also plot the empirical power of the tests based on $U_{\mathrm{mean}}$ and $U_{\mathrm{cov}}$ along with their asymptotic power expressions given in (\ref{Eq: theoretical power of U_1}) and (\ref{Eq: theoretical power of U_2}) respectively. Again we see that their empirical power closely follows the corresponding asymptotic power. By comparing these results with U-statistics, we confirm that the tests via sample splitting lose a factor of $\sqrt{2}$ in their power as a price for being dimension-agnostic.

\subsection{Tightness of the Gaussian limiting null distribution}
\begin{figure}[h!]
	\centering
	\begin{subfigure}[b]{0.31\textwidth}
		\includegraphics[width=\textwidth]{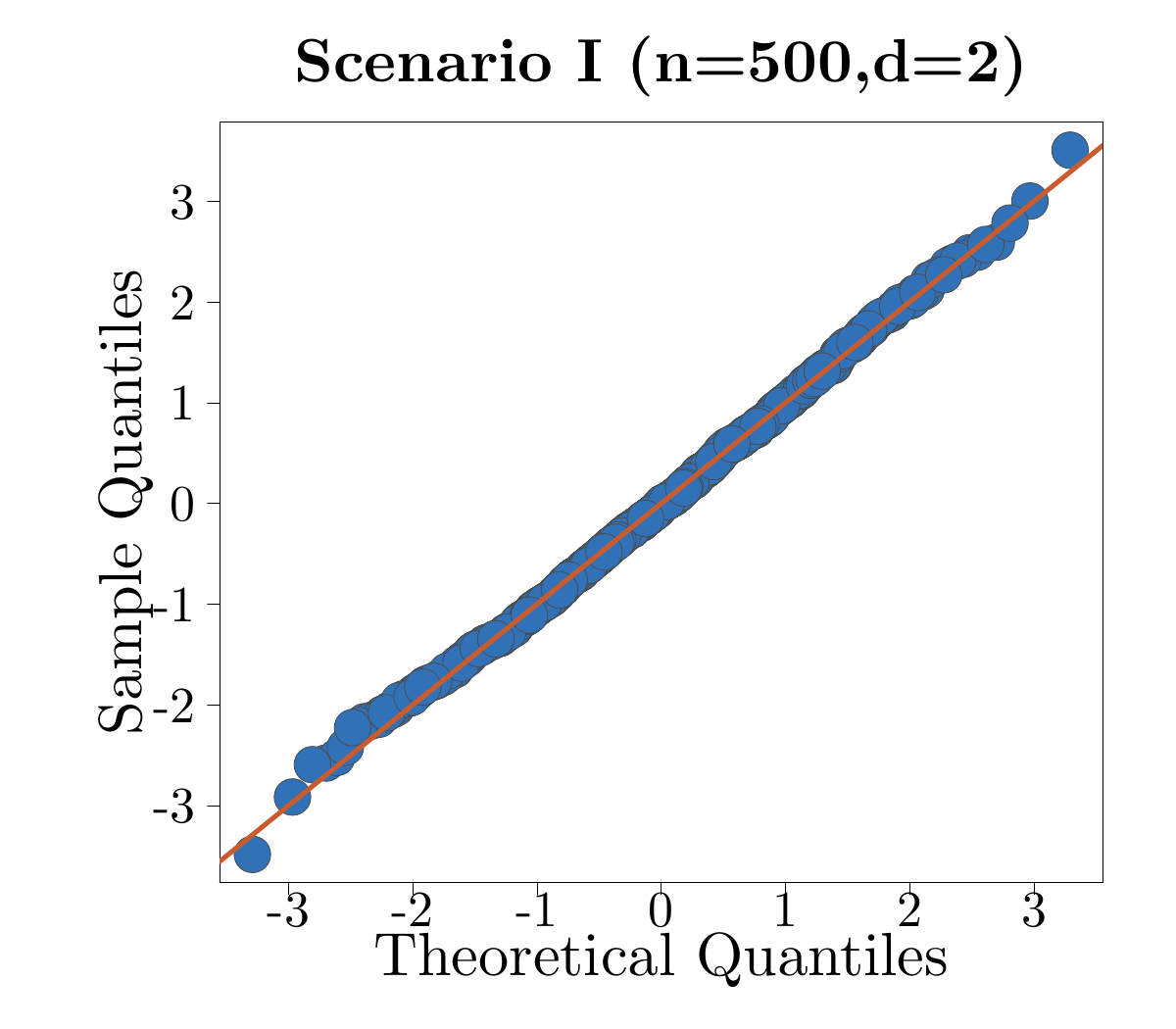}
	\end{subfigure}
	~ %add desired spacing between images, e. g. ~, \quad, \qquad, \hfill etc. 
	\begin{subfigure}[b]{0.31\textwidth}
		\includegraphics[width=\textwidth]{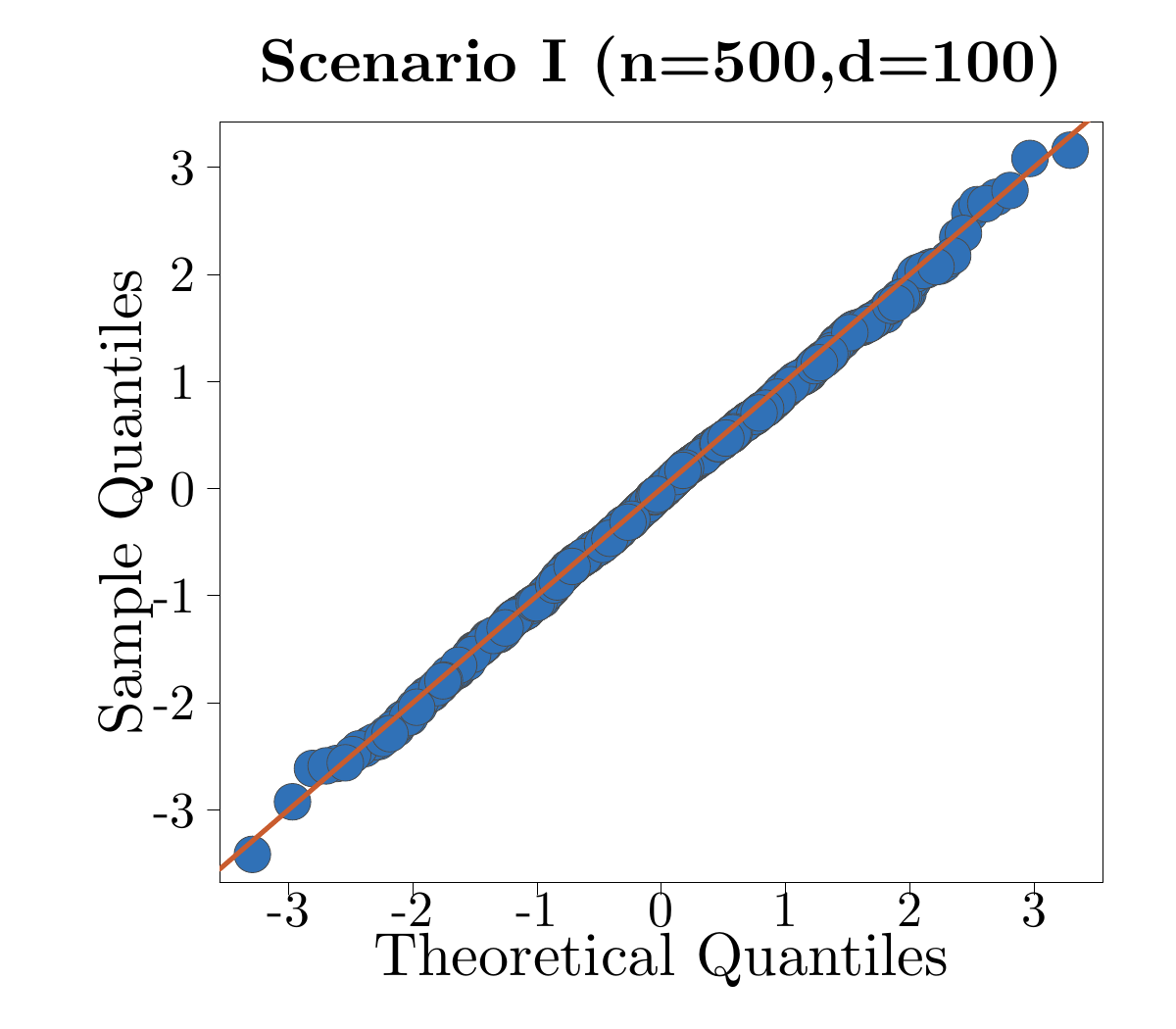}
	\end{subfigure}
	~ %add desired spacing between images, e. g. ~, \quad, \qquad, \hfill etc. 
	\begin{subfigure}[b]{0.31\textwidth}
		\includegraphics[width=\textwidth]{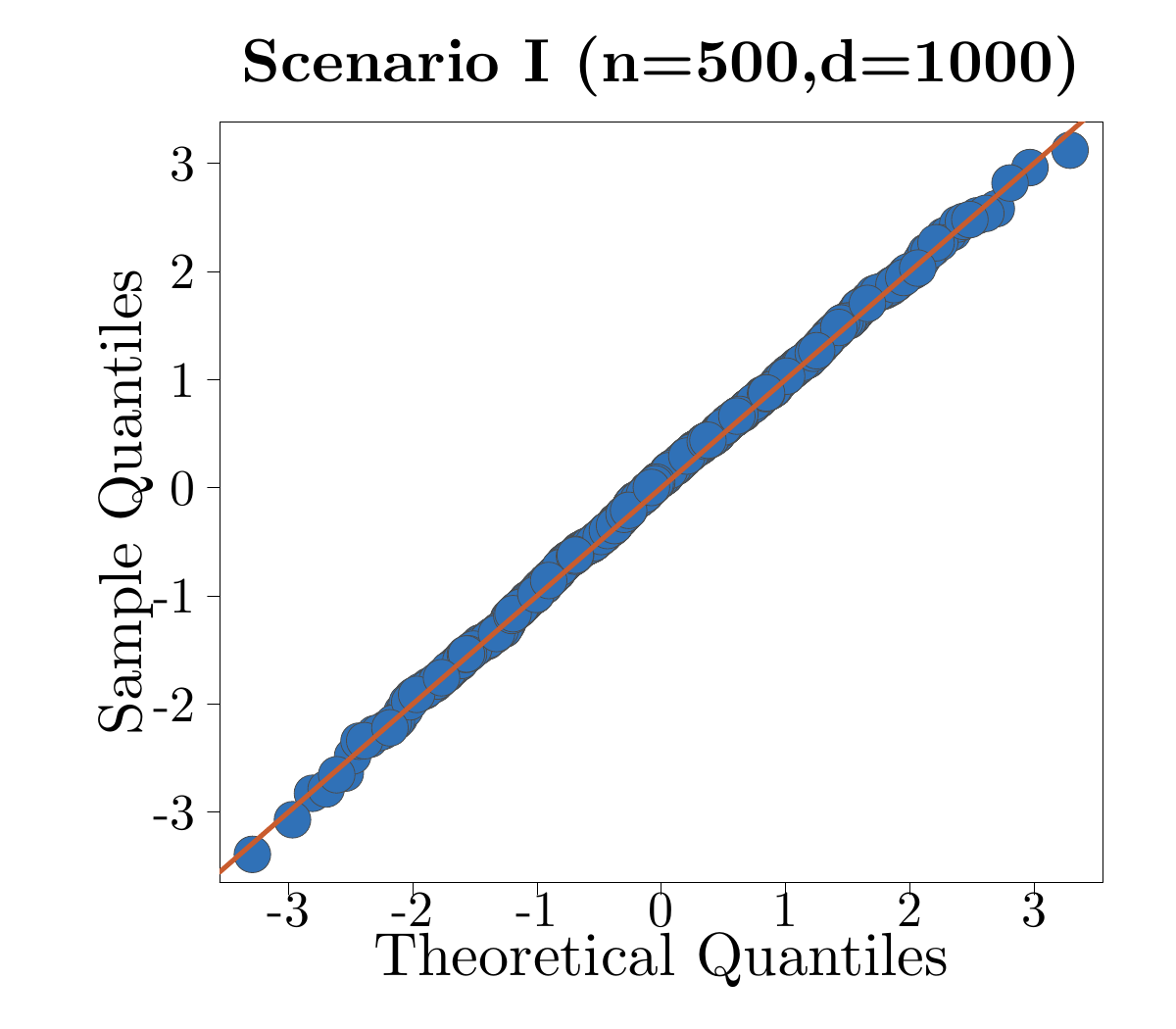}
	\end{subfigure}
	\vskip 1em
	\begin{subfigure}[b]{0.31\textwidth}
		\includegraphics[width=\textwidth]{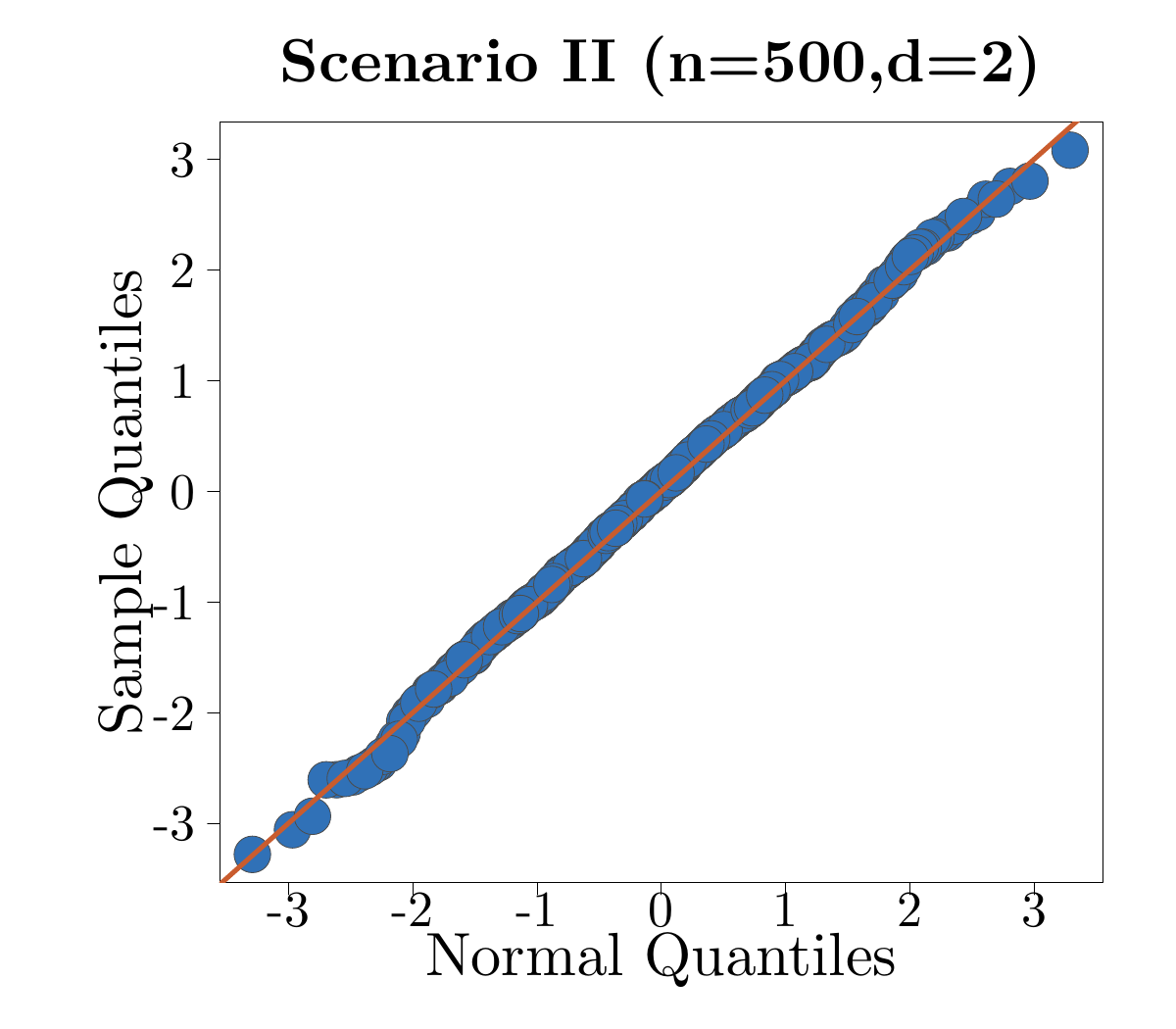}
	\end{subfigure}
	~ %add desired spacing between images, e. g. ~, \quad, \qquad, \hfill etc. 
	\begin{subfigure}[b]{0.31\textwidth}
		\includegraphics[width=\textwidth]{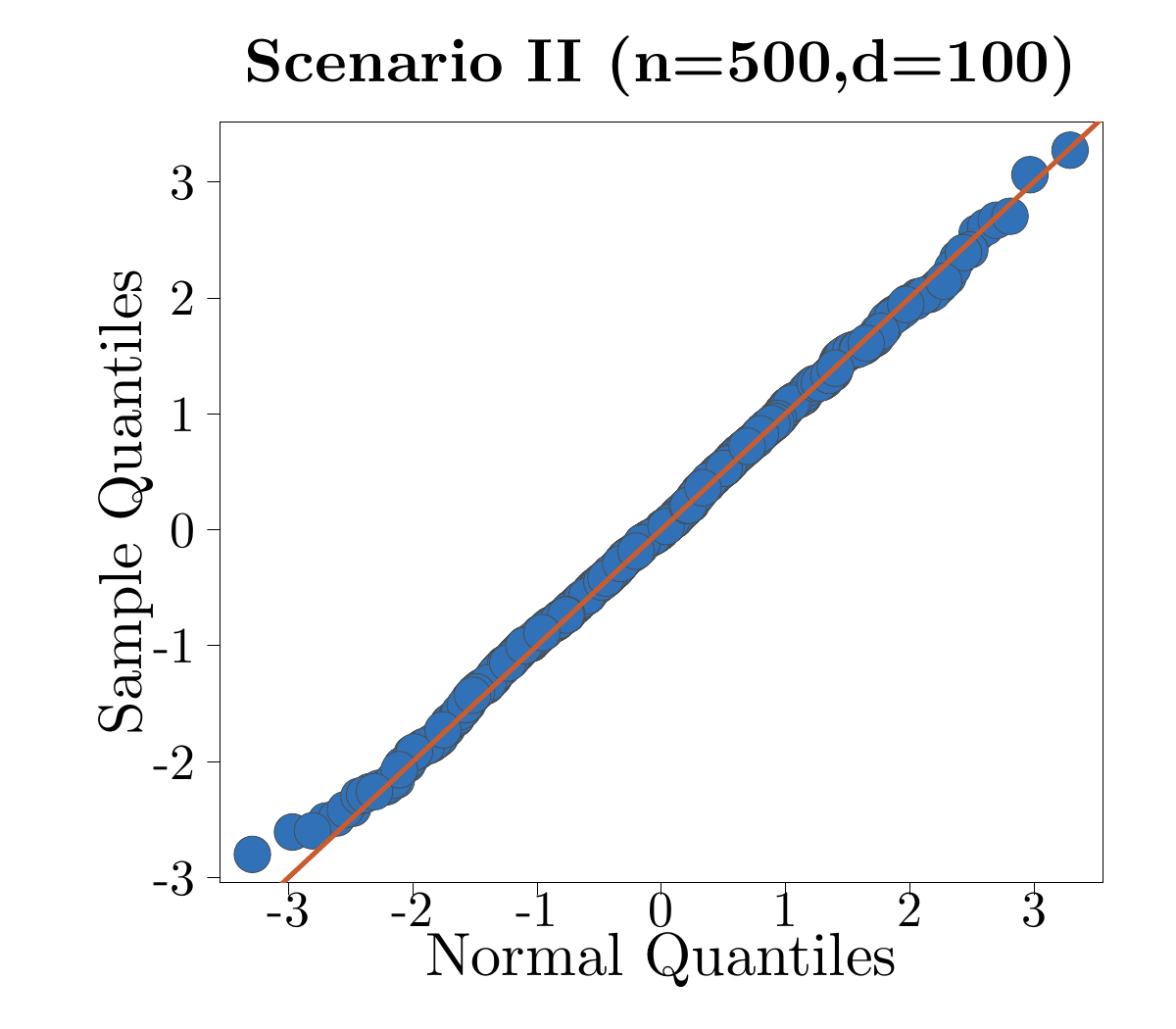}
	\end{subfigure}
	~ %add desired spacing between images, e. g. ~, \quad, \qquad, \hfill etc. 
	\begin{subfigure}[b]{0.31\textwidth}
		\includegraphics[width=\textwidth]{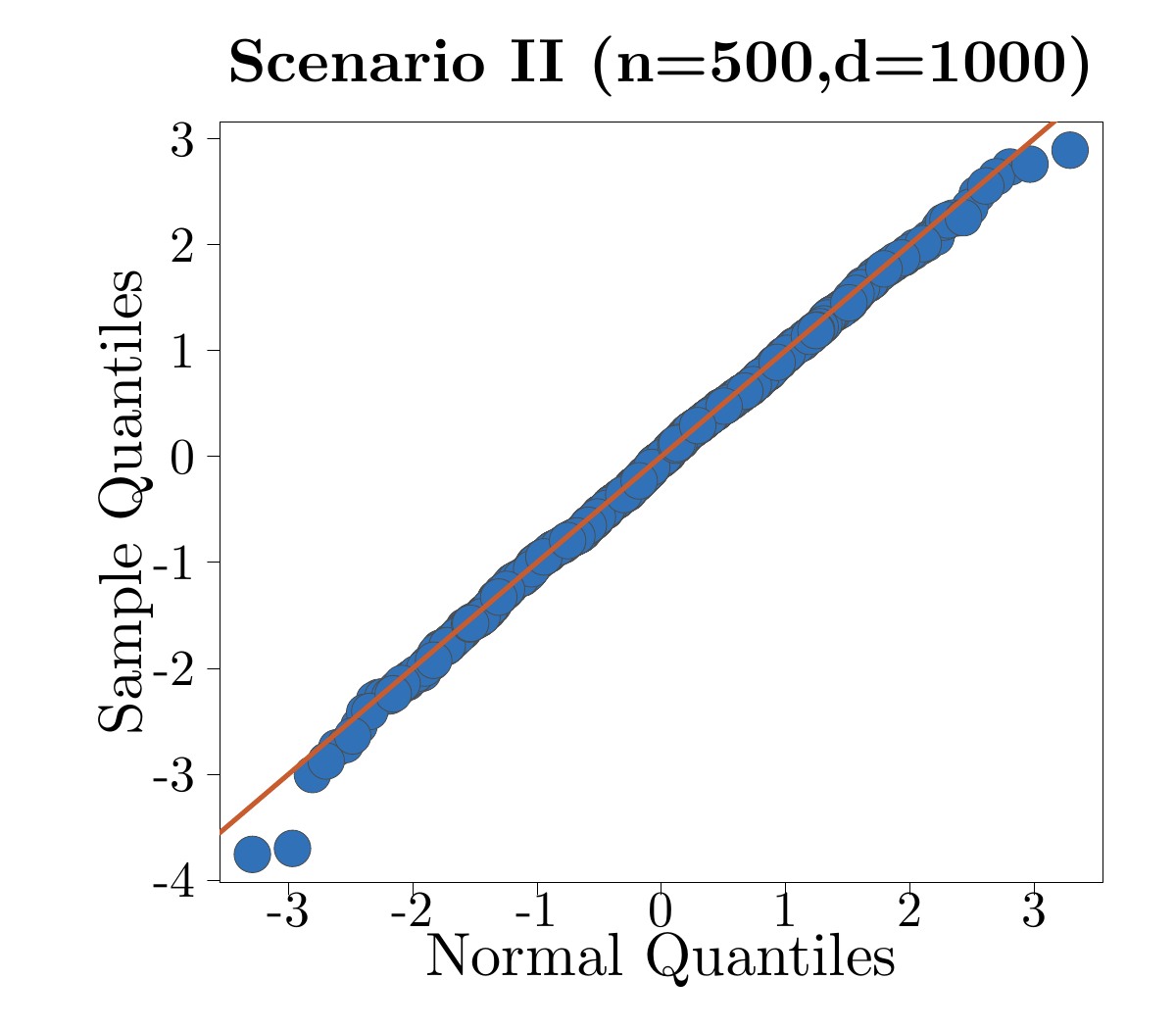}
	\end{subfigure}
	\caption{\small QQ plots of $T_{\mathrm{mean}}$. The results show that the null distribution of $T_{\mathrm{mean}}$ closely follows the standard normal distribution under both scenarios irrespective of the ratio of $n$ and $d$, which coincides with our theory in Theorem~\ref{Theorem: unconditional BE Bound}. The straight line $y=x$ is added as a reference point.}  \label{Figure: Mean 1}
\end{figure}

\begin{figure}[h!]
	\centering
	% \begin{subfigure}[b]{0.31\textwidth}
	\includegraphics[width=0.31\textwidth]{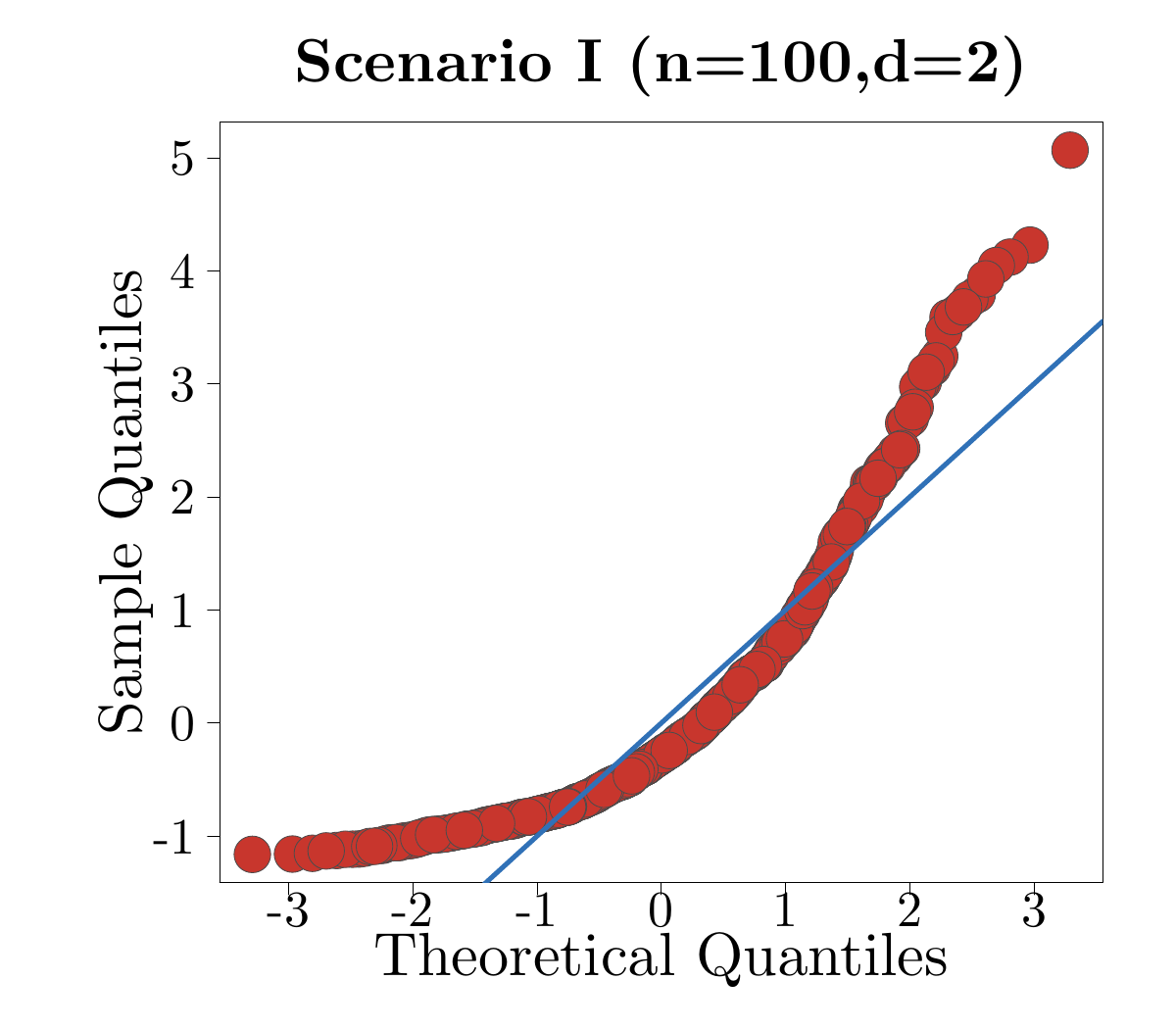}
	% \end{subfigure}
	% ~ %add desired spacing between images, e. g. ~, \quad, \qquad, \hfill etc. 
	% \begin{subfigure}[b]{0.31\textwidth}
	\includegraphics[width=0.31\textwidth]{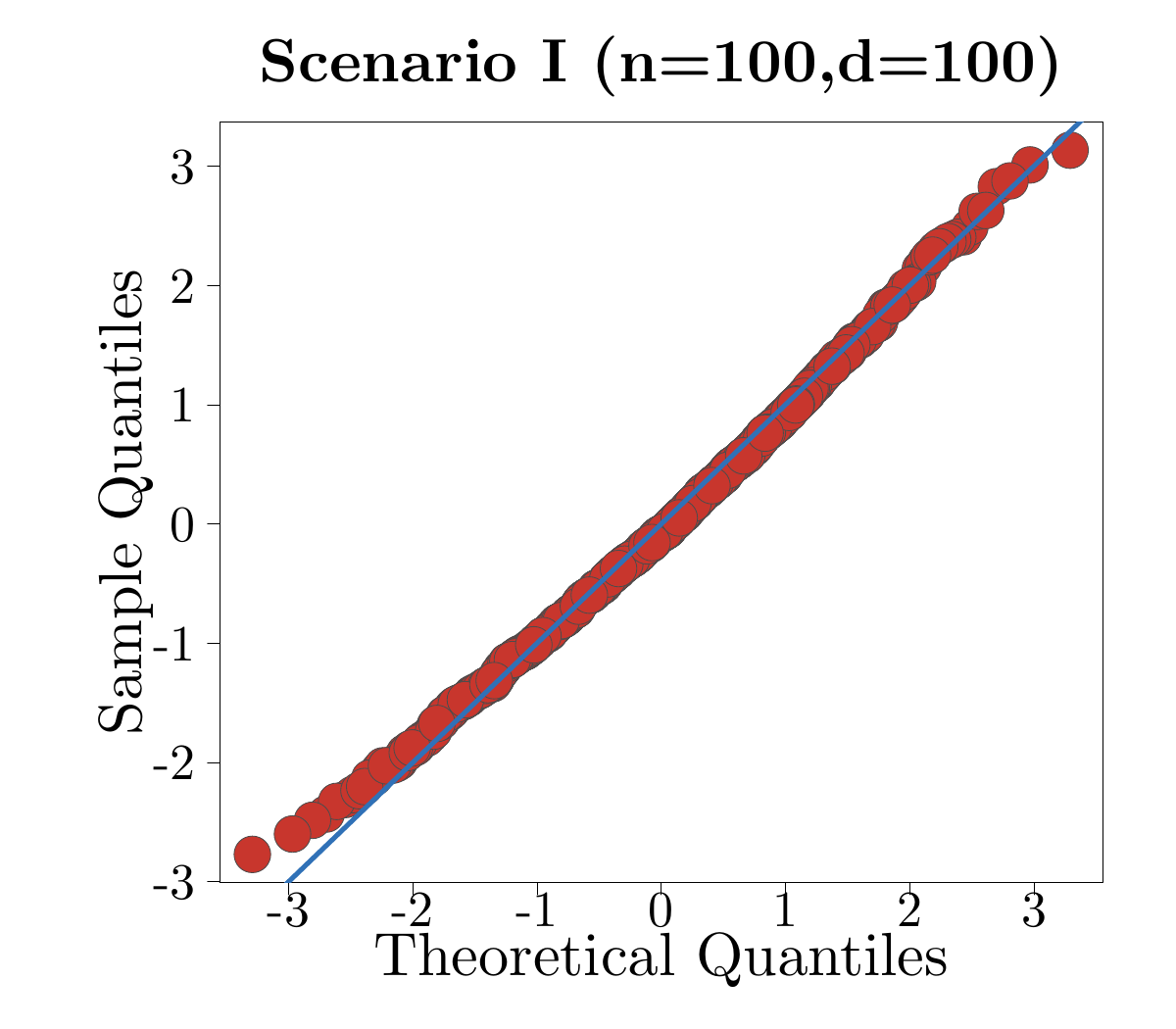}
	% \end{subfigure}
	% ~ %add desired spacing between images, e. g. ~, \quad, \qquad, \hfill etc. 
	% \begin{subfigure}[b]{0.31\textwidth}
	\includegraphics[width=0.31\textwidth]{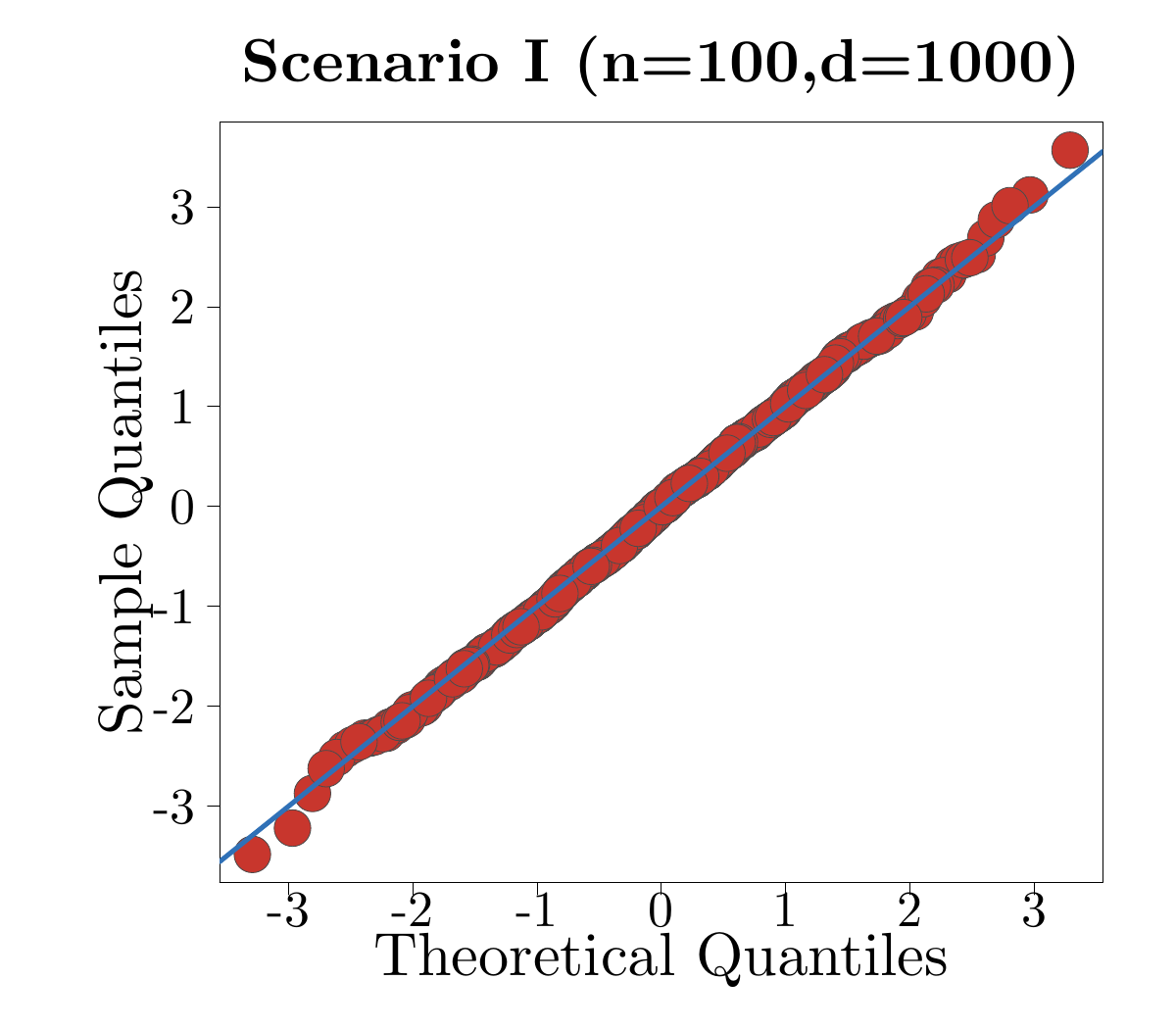}
	% \end{subfigure}
	% \vskip 1em
	% \begin{subfigure}[b]{0.31\textwidth}
	\includegraphics[width=0.31\textwidth]{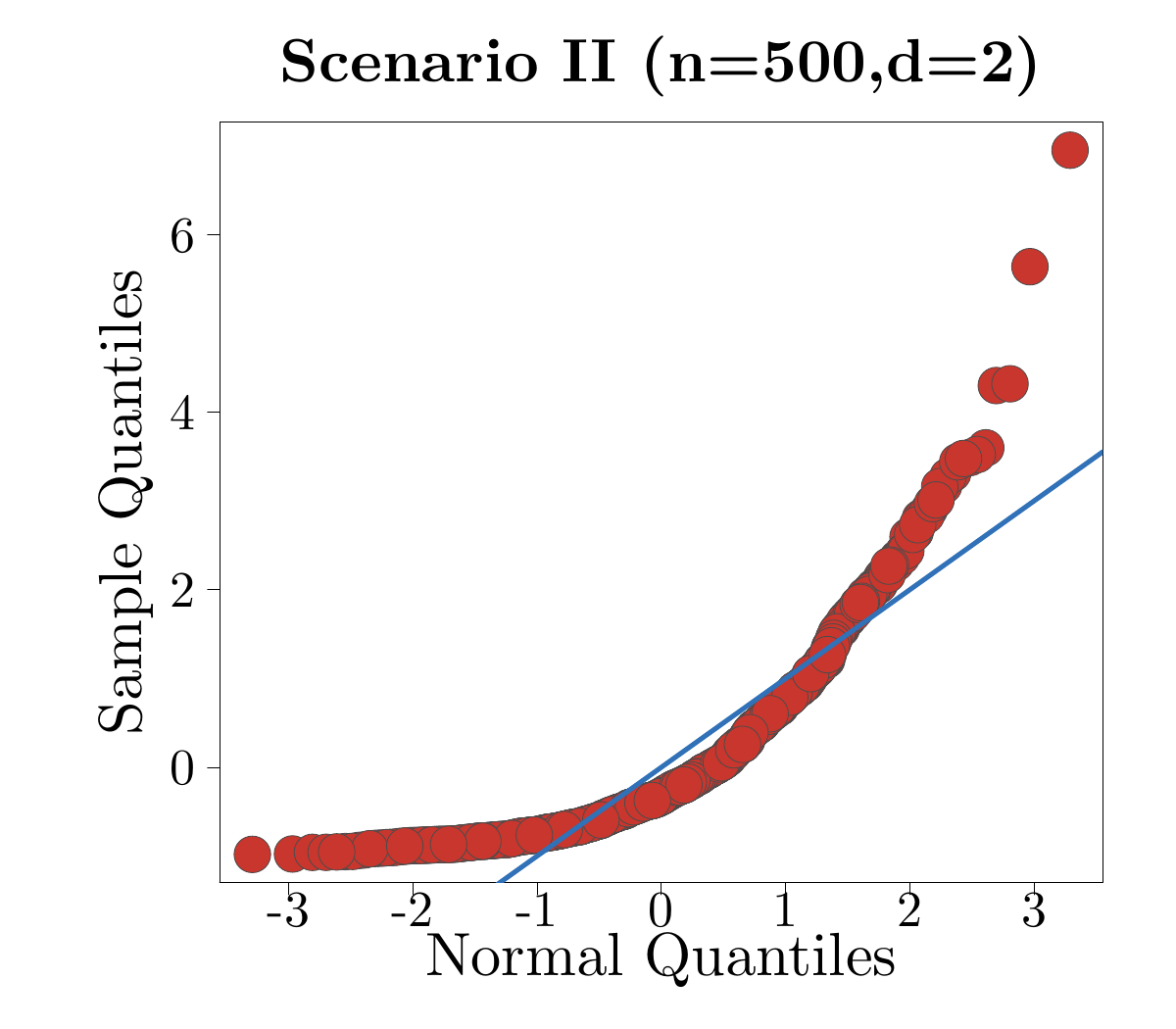}
	% \end{subfigure}
	% ~ %add desired spacing between images, e. g. ~, \quad, \qquad, \hfill etc. 
	% \begin{subfigure}[b]{0.31\textwidth}
	\includegraphics[width=0.31\textwidth]{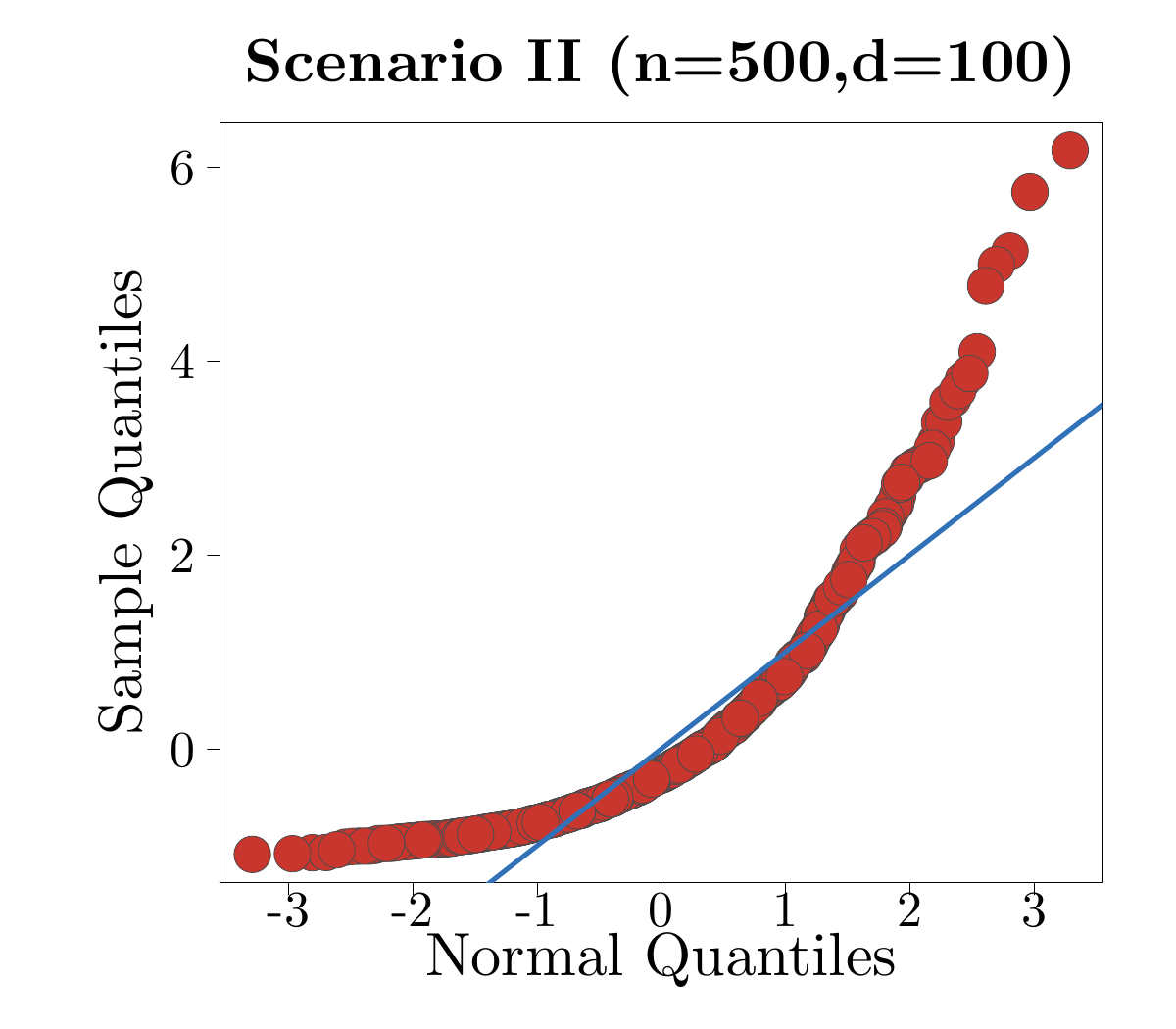}
	% \end{subfigure}
	% ~ %add desired spacing between images, e. g. ~, \quad, \qquad, \hfill etc. 
	% \begin{subfigure}[b]{0.31\textwidth}
	\includegraphics[width=0.31\textwidth]{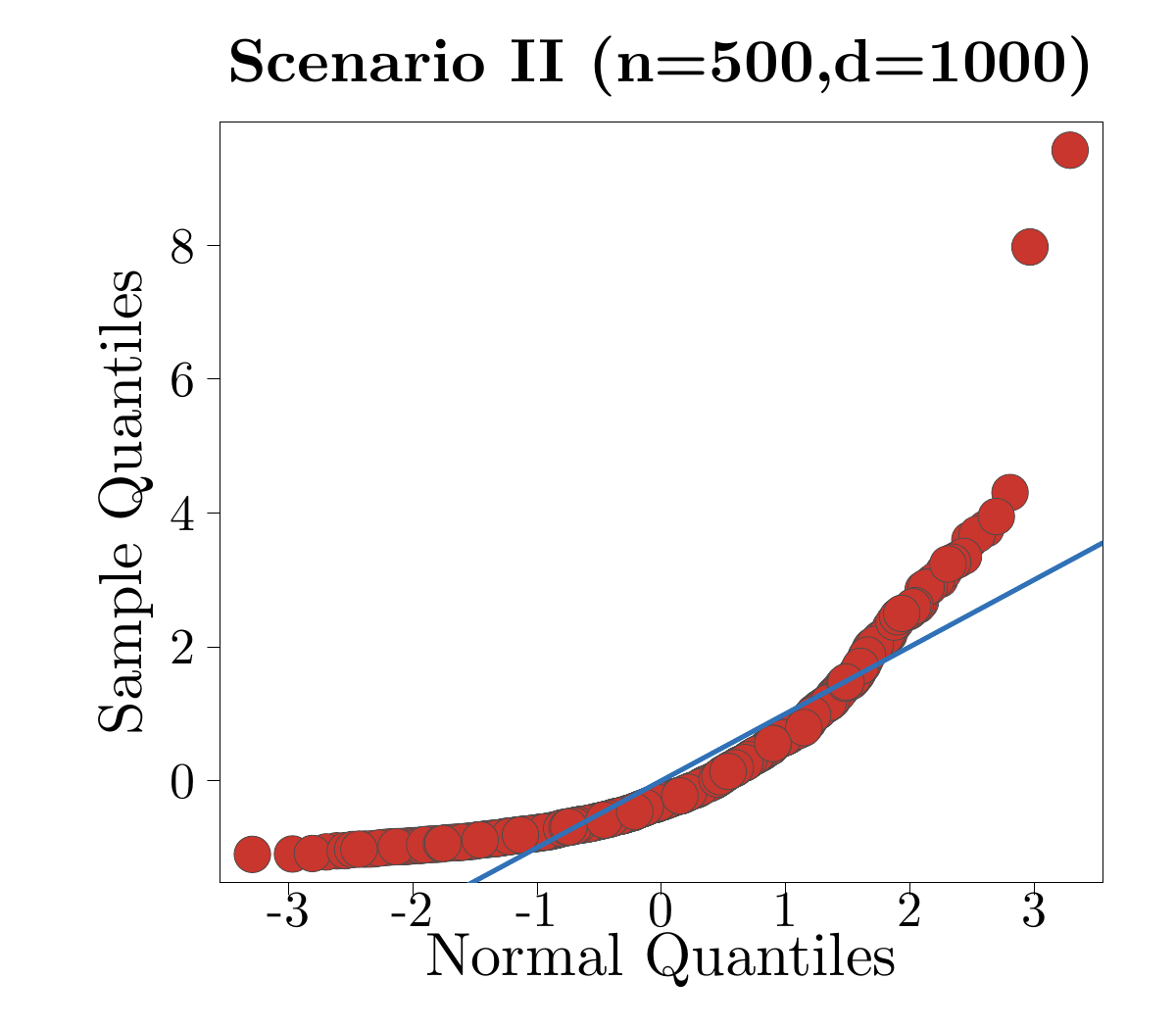}
	% \end{subfigure}
	\caption{\small QQ plots of $U_{\mathrm{mean}}$ scaled by its standard deviation $n^{-1}\sqrt{2\tr(\Sigma^2)}$. Under scenario I, the distribution of $U_{\mathrm{mean}}$ approaches $N(0,1)$ as $d$ increases. Under scenario II, in contrast, the distribution of $U_{\mathrm{mean}}$ is far from $N(0,1)$ regardless of $n,d$. The straight line $y=x$ is added as a reference point.}  \label{Figure: Mean 2}
\end{figure}

\begin{figure}[h!]
	\centering
	% \begin{subfigure}[b]{0.31\textwidth}
	\includegraphics[width=0.31\textwidth]{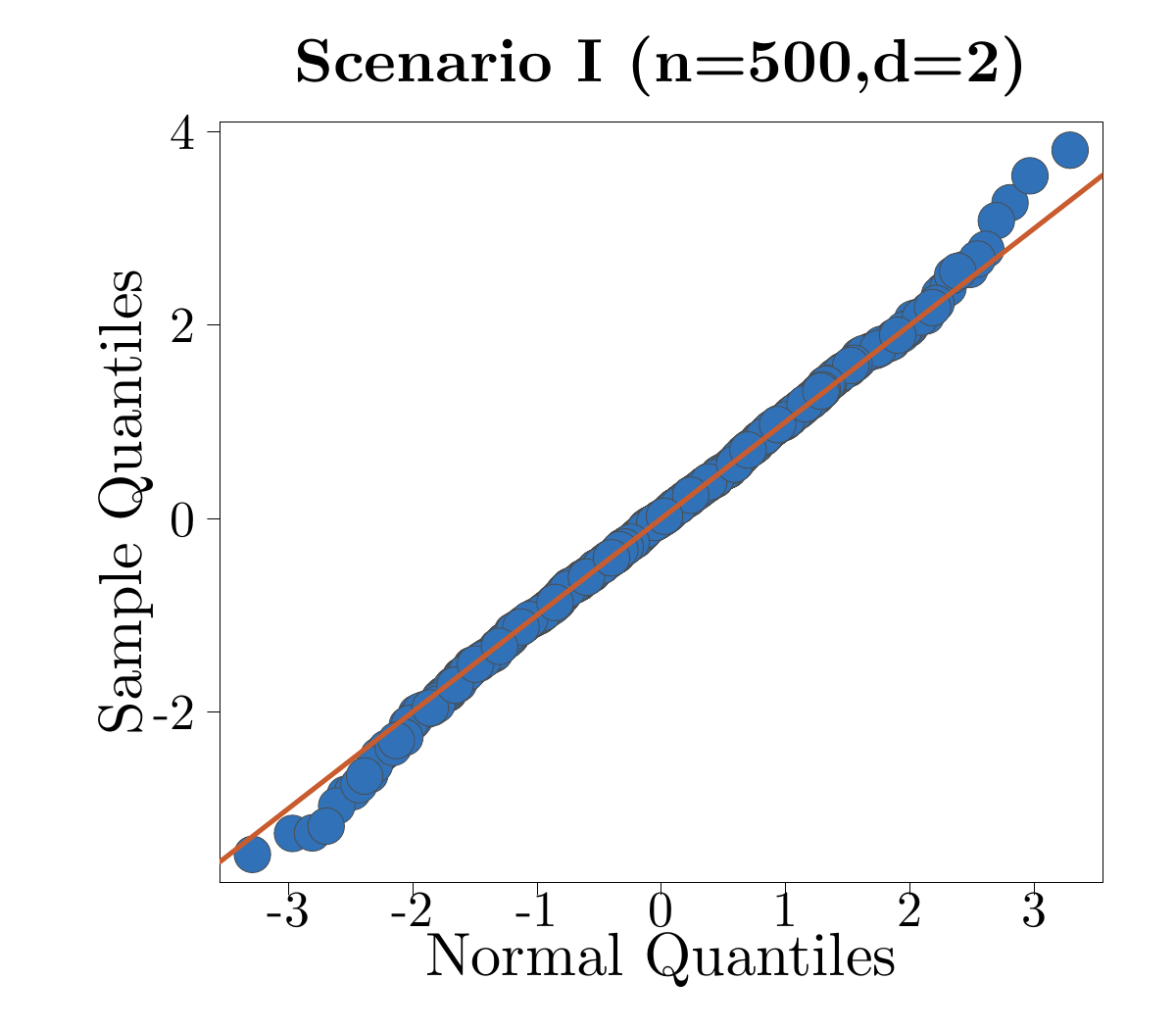}
	\includegraphics[width=0.31\textwidth]{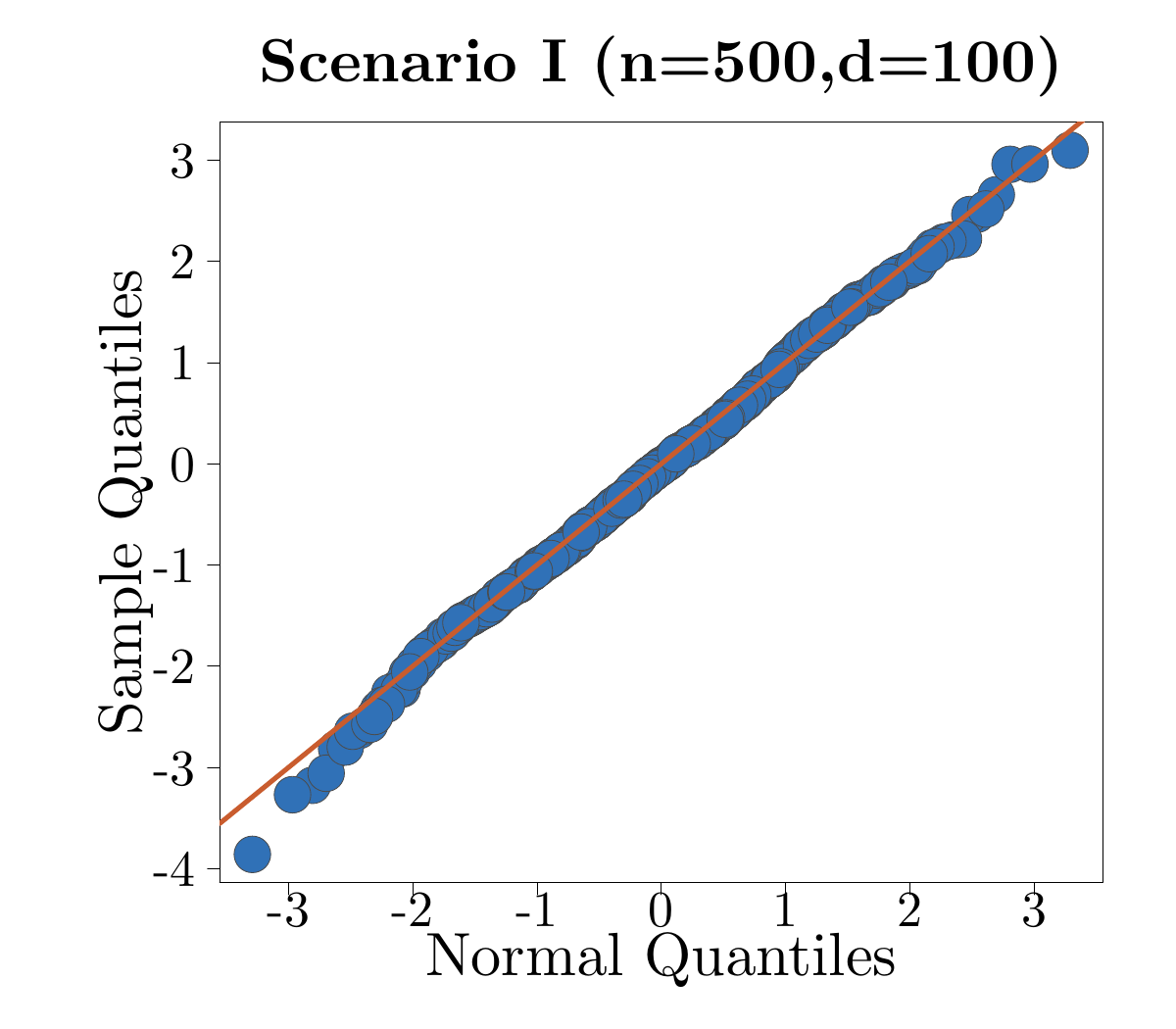}
	\includegraphics[width=0.31\textwidth]{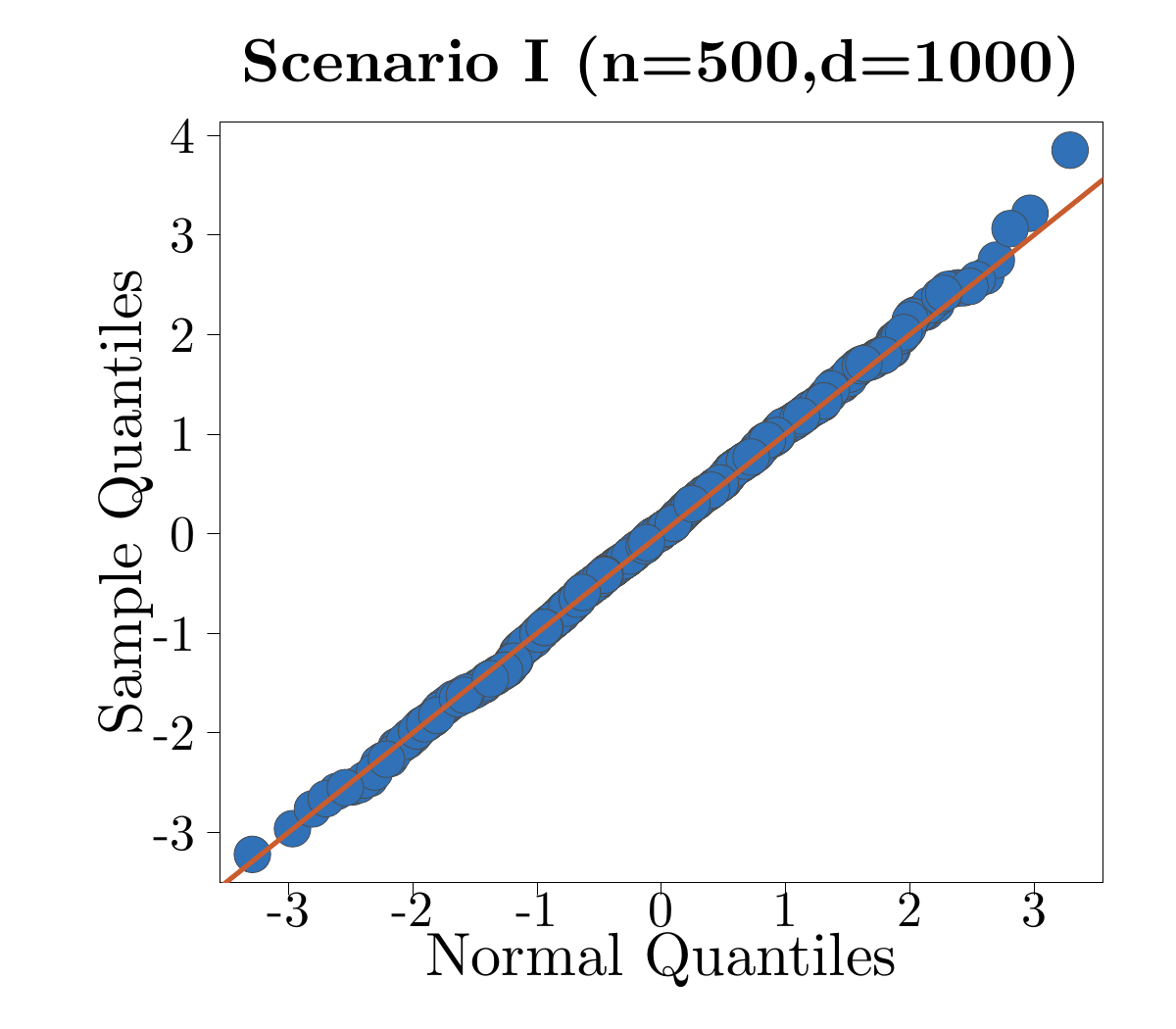}
	% \end{subfigure}
	% \vskip 1em
	% \begin{subfigure}[b]{0.31\textwidth}
	\includegraphics[width=0.31\textwidth]{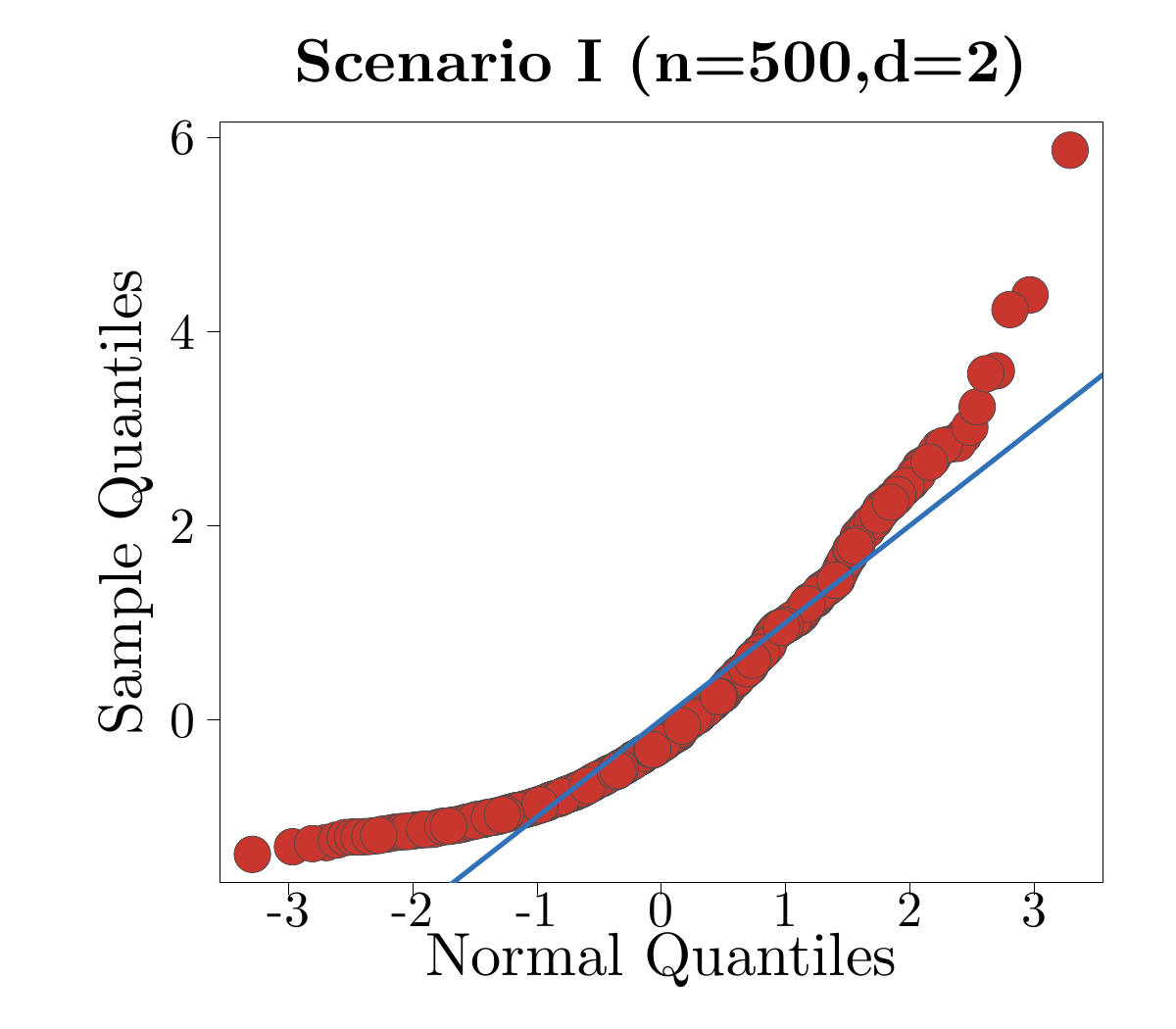}
	% \end{subfigure}
	% ~ %add desired spacing between images, e. g. ~, \quad, \qquad, \hfill etc. 
	% %(or a blank line to force the subfigure onto a new line)
	% \begin{subfigure}[b]{0.31\textwidth}
	\includegraphics[width=0.31\textwidth]{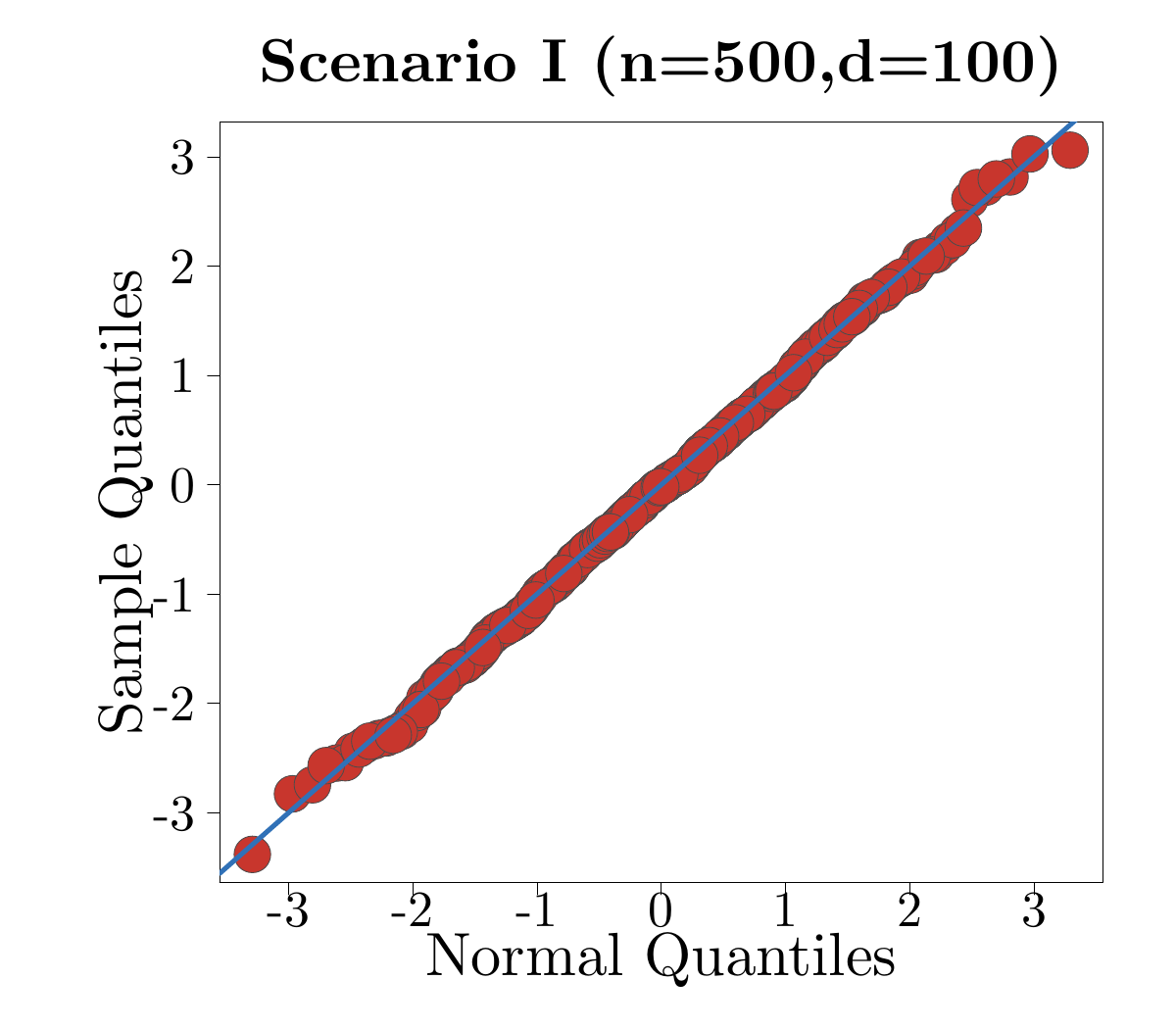}
	% \end{subfigure}
	% ~ %add desired spacing between images, e. g. ~, \quad, \qquad, \hfill etc. 
	% %(or a blank line to force the subfigure onto a new line)
	% \begin{subfigure}[b]{0.31\textwidth}
	\includegraphics[width=0.31\textwidth]{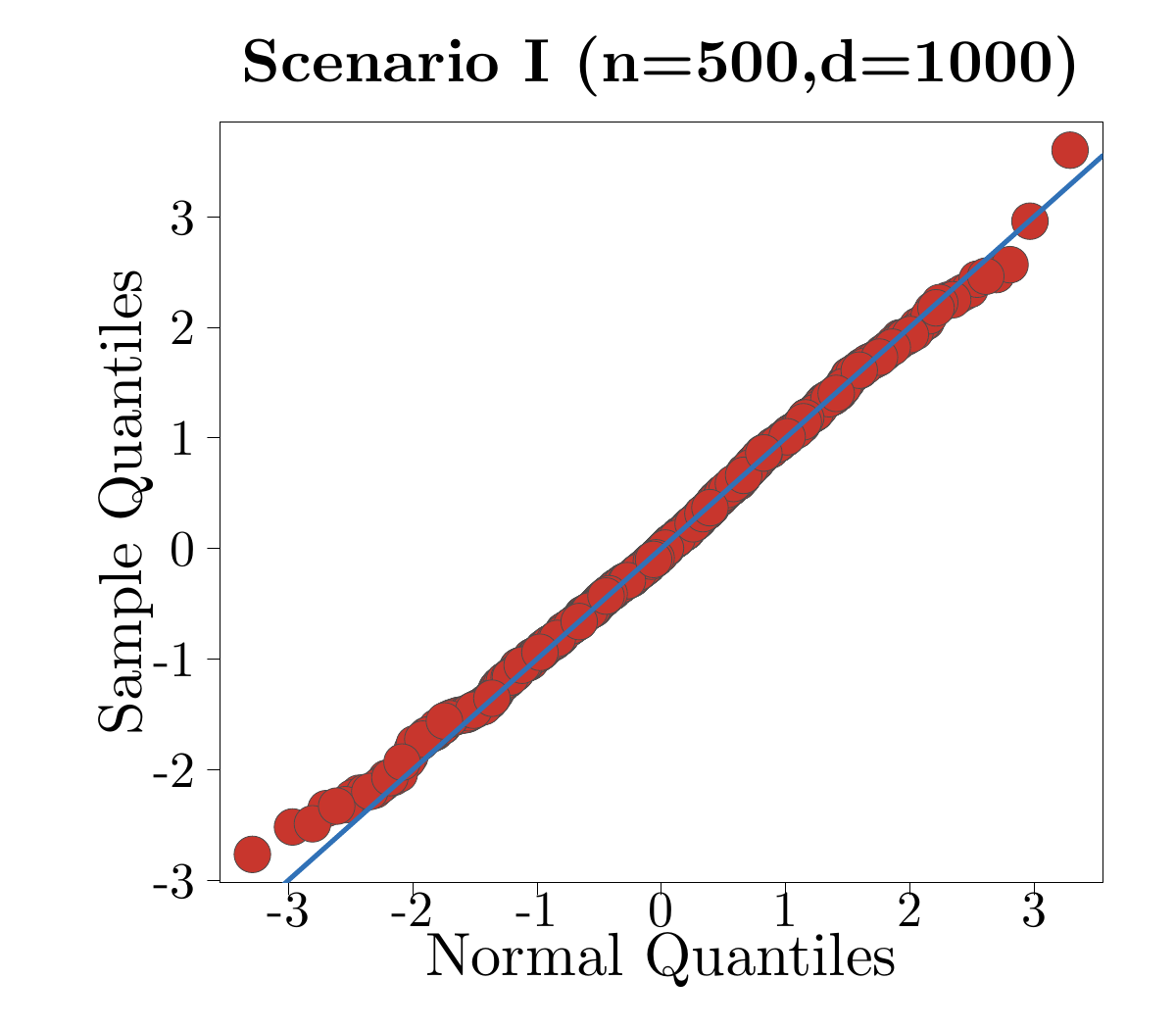}
	% \end{subfigure}
	\caption{\small QQ plots of $T_{\mathrm{cov}}$ (first row) and $U_{\mathrm{cov}}$ scaled by its standard deviation $2\sqrt{d(d+1)/n(n-1)}$ (second row) under scenario I. The distribution of $T_{\mathrm{cov}}$ consistently approximates $N(0,1)$ in all dimensions, whereas the distribution of $U_{\mathrm{cov}}$ is far from $N(0,1)$ when $d=2$. The straight line $y=x$ is added as a reference point.} \label{Figure: Covariance}
\end{figure}

In the second experiment, we illustrate the dimension-agnostic property of the proposed methods using QQ (quantile-quantile) plots. In particular, we demonstrate that the null distributions of our test statistics well approximate Gaussian in both fixed- and high-dimensional scenarios stated as follows.
\begin{itemize}
	\item \emph{Scenario I}: We observe a sample of size $n=500$ from a $d$-dimensional multivariate Gaussian distribution with the mean vector $\mu=(0,\ldots,0)^\top$ and the identity covariance matrix where the dimension is chosen among $d \in \{2,100,1000\}$. \vskip .5em
	\item \emph{Scenario II}: We consider the same scenario as before except that we choose a diagonal covariance matrix with its eigenvalues satisfying $\lambda_i(\Sigma) \propto \exp(i)$ and $\sum_{i=1}^d \lambda_i(\Sigma)=d$. In such setting, only a handful of eigenvalues are significantly greater than zero. 
\end{itemize}
By focusing on the balanced splitting ratio, i.e.~$m_1=m_2$, Figure~\ref{Figure: Mean 1} shows QQ plots that compare the quantiles of $N(0,1)$ with the empirical quantiles of $T_{\mathrm{mean}}$ based on 1000 repetitions. From Figure~\ref{Figure: Mean 1}, we see that most points lie on the straight line $y=x$ under both scenario I and scenario II, which implies that $T_{\mathrm{mean}}$ closes to $N(0,1)$ in distribution. This empirical observation in turn supports our theory in Theorem~\ref{Theorem: unconditional BE Bound}. Similarly, the first row of Figure~\ref{Figure: Covariance} indicates that the null distribution of $T_{\mathrm{cov}}$ approximates $N(0,1)$ well under scenario I, which can be predicted from Proposition~\ref{Proposition: BE bound for T2}. In contrast, the null distributions of the corresponding U-statistics vary a lot depending on the dimension as well as the condition on eigenvalues. In particular, Figure~\ref{Figure: Mean 2} shows that the distribution of the scaled $U_{\mathrm{mean}}$ differs from $N(0,1)$ when $d=2$ under both scenarios. Moreover, as we can see from the second row of Figure~\ref{Figure: Mean 2}, the distribution of $U_{\mathrm{mean}}$ does not approach to Gaussian even when the dimension increases. Similar behavior is observed for $U_{\mathrm{cov}}$ in the second row of Figure~\ref{Figure: Covariance}. 
These empirical results highlight the benefit of the proposed method via sample-splitting that is easier to calibrate and more robust to model assumptions than the corresponding U-statistics. 

\section{Proofs of all results} \label{Section: Proofs}

In this section we collect technical details and proofs of our main results. In addition to the notation introduced in Section~\ref{Section: Notation}, we need to introduce some additional notation used throughout this section.

\vskip 1em
 
\noindent \textbf{Additional notation.} For a sequence of random variables $X_n$ and constants $a_n$, we write $X_n = O_P(a_n)$ if $a_n^{-1} X_n$ is stochastically bounded and $X_n = o_P(a_n)$ if $a_n^{-1} X_n$ converges to zero in probability or equivalently $a_n^{-1} X_n \convP 0$. For real sequences $a_n$ and $b_n$, we write $a_n \lesssim b_n$ or $b_n \gtrsim a_n$ if there exists a constant $C$ such that $a_n \leq C b_n$ for all $n$.

To simplify the presentation, we often suppress dependence on $P$ and write $\mE$ instead of $\mE_{P}$.

\subsection{Details on Assumption~\ref{assumption: moment condition}} \label{Section: Details on the moment condition}
In this subsection, we demonstrate that Assumption~\ref{assumption: moment condition} is implied by the model assumption made in \cite{bai1996effect,chen2010two,wang2019feasible}. We first recall their multivariate model assumption:
\begin{align} \label{Eq: pseudo independence}
	X = \mu + \Gamma Z,
\end{align}
where $\Gamma$ is a $d \times r$ matrix with some $r \geq d$ and thus $\mathrm{cov}(X)=  \Gamma \Gamma^\top$. Here $Z = (z_1,\ldots,z_r)^\top$ satisfies $\mE[Z]=0$, $\mV[Z] = I$, $\mE[z_i^4] = 3 + \Delta < \infty$ and $\mE[z_i^8] < C$ for $i=1,\ldots,r$ where $\Delta$ is a difference between the fourth moment of $z_i$ and $N(0,1)$ and $C$ is some positive constant. Furthermore it is assumed that 
\begin{align*}
	\mE[z_{\ell_1}^{\alpha_1}z_{\ell_2}^{\alpha_2} \cdots z_{\ell_q}^{\alpha_q}] = \mE[z_{\ell_1}^{\alpha_1}] \mE[z_{\ell_2}^{\alpha_2}] \cdots \mE[z_{\ell_q}^{\alpha_q}] 
\end{align*}
for a positive integer $q$ such that $\sum_{\ell=1}^q \alpha_\ell \leq 8$ and $\{\ell_1,\ldots,\ell_q\}$ are distinct indices. 

Now, for a given nonzero $u \in \mathbb{R}^d$ and by assuming $\mu=0$, Lemma 1 of \cite{wang2019feasible} yields
\begin{align*}
	\mE[(u^\top X)^4] = \mE[(Z^\top \Gamma^\top u u^\top \Gamma Z )^2] \leq (3+ \Delta) \tr^2(\Gamma^\top u u^\top \Gamma) = (3 + \Delta) \{ \mE[(u^\top X)^2]\}^2.
\end{align*} 
This result together with Jensen's inequality implies 
\begin{align*}
	\left[ \frac{\{\mE |u^\top X|^3 \}^{1/3} }{\{ \mE |u^\top X|^2 \}^{1/2} } \right]^4 \leq \frac{\mE [(u^\top X)^4] }{\{ \mE[(u^\top X)^2]\}^2} \leq 3 + \Delta.
\end{align*}
Thus, under the model~(\ref{Eq: pseudo independence}), Assumption~\ref{assumption: moment condition} is satisfied with $K = (3+\Delta)^{1/4}$. 

\textcolor{black}{In order to demonstrate the dimension-agnostic property of $T_{\mathrm{cov}}$ under a non-Gaussian setting, we define a class of distributions for which bound~(\ref{Eq: Sufficient condition for Normality}) holds with a potentially different constant. 
	\begin{definition}[Covariance model] \label{Definition: covariance model}
		Define a class of distribution $\mathcal{P}_{0,d,\mathrm{cov}}$ where a random vector $X \sim P \in\mathcal{P}_{0,d,\mathrm{cov}}$ is continuous and satisfies model~(\ref{Eq: pseudo independence}) with $\mu = 0$, $d=r$ and $\Gamma = I$. Moreover, assume that there exists a small constant $c>0$ such that $\Delta + 2 > c$. 
	\end{definition}
}

\subsection{Proof of Theorem~\ref{Theorem: unconditional BE Bound}}
First note that $\mE_P[f_{\mathrm{mean}}^2(X)| \mathcal{X}_2] = \widehat{\mu}_2^\top \Sigma \widehat{\mu}_2$ is strictly greater zero and finite almost surely under Assumption~\ref{assumption: moment condition}. Thus, based on Lemma~\ref{Lemma: Conditional BE Bound}, the studentized statistic has the following conditional guarantee
\begin{align*} 
	\sup_{-\infty < z < \infty} \big| P(T_{\mathrm{mean}}  \leq z | \mathcal{X}_2) - \Phi(z) \big| \leq   \frac{C \mE_P[|f_{\mathrm{mean}}^3(X)| | \mathcal{X}_2]}{\{\mE_P[f_{\mathrm{mean}}^2(X)|\mathcal{X}_2]\}^{3/2} \sqrt{m_1}} \leq  \frac{C L}{\sqrt{m_1}} \quad \text{a.s.,} 
\end{align*}
where the second inequality uses Assumption~\ref{assumption: moment condition}. Now notice that
\begin{align*}
	\sup_{-\infty < z < \infty} \big| P(T_{\mathrm{mean}}  \leq z) - \Phi(z) \big| ~= ~ & \sup_{-\infty < z < \infty} \big| \mE_P[P(T_{\mathrm{mean}}  \leq z | \mathcal{X}_2)] - \Phi(z) \big| \\[.5em]
	\leq~&  \mE_P \bigg[	\sup_{-\infty < z < \infty} \big| P(T_{\mathrm{mean}}  \leq z | \mathcal{X}_2) - \Phi(z) \big| \bigg] \\[.5em]
	\leq ~ & \frac{C L}{\sqrt{m_1}},
\end{align*}
and the claim follows by taking the supremum over $\mathcal{P}_{0,d}(L)$ on both sides. 

\subsection{Proof of Theorem~\ref{Theorem: asymptotic power}}
By recalling $\widehat{\mu}_2 = m_2^{-1} \sum_{i=m_1+1}^n X_i$ and denoting $\widehat{\mu}_1 = m_1^{-1} \sum_{i=1}^{m_1}X_i$, the power function can be written as
\begin{align*}
	& P(T_{\mathrm{mean}} > z_{1-\alpha}) \\[.5em] 
	=~ & P  \left( \frac{\sqrt{m_1} \widehat{\mu}_2^\top (\widehat{\mu}_1  - \mu)}{\sqrt{m_1^{-1} \sum_{i=1}^{m_1} \{\widehat{\mu}_2^\top (X_i - \widehat{\mu}_1 )\}^2}} > z_{1-\alpha}  - \frac{\sqrt{m_1} \widehat{\mu}_2^\top \mu}{\sqrt{{m_1}^{-1} \sum_{i=1}^{m_1} \{\widehat{\mu}_2^\top (X_i - \widehat{\mu}_1 )\}^2}} \right).
\end{align*}
Note that $X-\mu$ satisfies Assumption~\ref{assumption: moment condition} under the Gaussian assumption and so it follows that  
\begin{align} \label{Eq: normal approximation}
	\frac{\sqrt{m_1} \widehat{\mu}_2^\top (\widehat{\mu}_1  - \mu)}{\sqrt{{m_1}^{-1} \sum_{i=1}^{m_1} \{\widehat{\mu}_2^\top (X_i - \widehat{\mu}_1 )\}^2}} \convD N(0,1).
\end{align}
Now we claim that under Assumption~\ref{assumption: one-sample mean testing}:
\begin{enumerate}\setlength{\itemindent}{0.073in}
	\item[\textbf{(C1).}] $\displaystyle  \sqrt{m_1} \widehat{\mu}_2^\top \mu = \sqrt{m_1} \mu^\top \mu + O_P\left( \! \sqrt{\frac{m_1}{m_2} \mu^\top \Sigma \mu}\right)$,
	\item[\textbf{(C2).}] $\displaystyle  \widehat{\mu}_2^\top \Sigma \widehat{\mu}_2 = \frac{1}{m_2} \tr (\Sigma^2) + \mu^\top \Sigma \mu + O_P \left( \! \sqrt{\frac{2}{m_2^2} \tr(\Sigma^4) + \frac{4}{m_2} \mu^\top \Sigma^3 \mu} \right)$,
	\item[\textbf{(C3).}] $\displaystyle  \frac{\widehat{\mu}_2^\top \widehat{\Sigma} \widehat{\mu}_2}{\widehat{\mu}_2^\top \Sigma \widehat{\mu}_2} \convP 1$ where $\displaystyle \widehat{\Sigma} = m_1^{-1} \sum_{i=1}^{m_1} (X_i - \widehat{\mu}_1 )(X_i - \widehat{\mu}_1 )^\top$.
\end{enumerate}
Once these conditions are given, it follows that
\begin{align*}
	\frac{\sqrt{m_1} \widehat{\mu}_2^\top \mu}{\sqrt{{m_1}^{-1} \sum_{i=1}^{m_1} \{\widehat{\mu}_2^\top (X_i - \widehat{\mu}_1 )\}^2}}  ~=~ & \frac{\sqrt{m_1} \widehat{\mu}_2^\top \mu}{\sqrt{\widehat{\mu}_2^\top \Sigma \widehat{\mu}_2}} \times \sqrt{\frac{\widehat{\mu}_2^\top \Sigma \widehat{\mu}_2}{\widehat{\mu}_2^\top \widehat{\Sigma} \widehat{\mu}_2}} \\[.5em]
	\overset{\text{(i)}}{=} ~ & \frac{\sqrt{m_1} \widehat{\mu}_2^\top \mu}{\sqrt{\widehat{\mu}_2^\top \Sigma \widehat{\mu}_2}} \{ 1 + o_P(1) \} \\[.5em]
	\overset{\text{(ii)}}{=} ~ & \frac{ \sqrt{m_1} \mu^\top \mu + O_P (n^{-1/4})}{\sqrt{m_2^{-1}\tr(\Sigma^2) + O_P(n^{-1/2}) }} \{ 1 + o_P(1) \} \\[.5em]
	\overset{\text{(iii)}}{=} ~ & \frac{\sqrt{m_1} \mu^\top \mu}{\sqrt{m_2^{-1} \tr(\Sigma^2)}} + o_P(1), 
\end{align*}
where step~(i) uses (C3) with the continuous mapping theorem and step~(ii) follows by (C1) and (C2) along with Assumption~\ref{assumption: one-sample mean testing}. The last step~(iii) holds since $\sqrt{m_1} \mu^\top \mu = O(1)$ and $m_2^{-1} \tr(\Sigma^2)$ is strictly greater than zero and finite under the eigenvalue condition in Assumption~\ref{assumption: one-sample mean testing}. Combining the last equality~(iii) with the normal approximation~(\ref{Eq: normal approximation}), we arrive at the conclusion that 
\begin{align*}
	P(T_{\mathrm{mean}} > z_{1-\alpha})  ~=~  & \Phi\left( \! z_\alpha + \frac{\sqrt{m_1} \mu^\top \mu}{\sqrt{m_2^{-1} \tr(\Sigma^2)}}  \right) + o(1) \\[.5em]
	= ~ &  \Phi \left( \! z_{\alpha} + \frac{n\sqrt{\kappa(1-\kappa)}\mu^\top \mu}{\sqrt{\tr(\Sigma^2)}} \right) + o(1).
\end{align*} 
To complete the proof, we verify the conditions (C1), (C2) and (C3) in order. 

\vskip 1em

\noindent \textbf{$\bullet$ Verification of (C1) and (C2).} It is clear to see that 
\begin{align*}
	\mE_P[\sqrt{m_1} \widehat{\mu}_2^\top \mu] = \sqrt{m_1} \mu^\top \mu \quad \text{and} \quad \mV_P[\sqrt{m_1} \widehat{\mu}_2^\top \mu] = \frac{m_1}{m_2} \mu^\top \Sigma \mu. 
\end{align*}
The condition (C1) follows since $X_n = \mE_P[X_n] + O_P(\sqrt{\mV_P[X_n]})$. 

For (C2), we first note the following moment calculations
\begin{align*}
	&\mE_P[\xi^\top \Lambda \xi ] ~=~ \tr(\Lambda \Sigma) + \mu^\top \Lambda \mu, \\[.5em]
	&\mV_P[\xi^\top \Lambda \xi ] ~=~ 2 \tr(\Lambda \Sigma \Lambda \Sigma) + 4 \mu^\top \Lambda \Sigma \Lambda \mu,
\end{align*}
where $\xi \sim N(\mu, \Sigma)$. Then, using the fact that $\widehat{\mu}_2 \sim N(\mu, m_2^{-1} \Sigma)$, we can obtain
\begin{align*}
	\mE_P[\widehat{\mu}_2^\top \Sigma \widehat{\mu}_2] ~=~ \frac{1}{m_2} \tr (\Sigma^2) + \mu^\top \Sigma \mu \quad \text{and} \quad  \mV_P[\widehat{\mu}_2^\top \Sigma ] ~=~ \frac{2}{m_2^2} \tr(\Sigma^4) + \frac{4}{m_2} \mu^\top \Sigma^3 \mu.
\end{align*}
Thus the condition (C2) follows similarly as (C1).

\vskip 1em

\noindent \textbf{$\bullet$ Verification of (C3).} We want to prove that for any $\epsilon >0$,
\begin{align*}
	P \bigg( \bigg|  \frac{\widehat{\mu}_2^\top \widehat{\Sigma} \widehat{\mu}_2}{\widehat{\mu}_2^\top \Sigma \widehat{\mu}_2}  -  1 \bigg|  > \epsilon \bigg) \rightarrow 0 \quad \text{as $n \rightarrow \infty$.}
\end{align*}
By the law of total expectation and conditional Chebyshev's inequality,
\begin{align*}
	P \bigg( \bigg|  \frac{\widehat{\mu}_2^\top \widehat{\Sigma} \widehat{\mu}_2}{\widehat{\mu}_2^\top \Sigma \widehat{\mu}_2}  -  1 \bigg|  > \epsilon \bigg) ~\leq~ \frac{1}{\epsilon^2} \mE_P \bigg[ \mE_P \bigg\{ \bigg( \frac{\widehat{\mu}_2^\top \widehat{\Sigma} \widehat{\mu}_2 - \widehat{\mu}_2^\top \Sigma \widehat{\mu}_2}{\widehat{\mu}_2^\top \Sigma \widehat{\mu}_2} \bigg)^2 \bigg| \widehat{\mu}_2 \bigg\} \bigg]
\end{align*} 
and the bias-variance decomposition gives
\begin{align*}
	\mE_P \bigg\{ \bigg( \frac{\widehat{\mu}_2^\top \widehat{\Sigma} \widehat{\mu}_2 - \widehat{\mu}_2^\top \Sigma \widehat{\mu}_2}{\widehat{\mu}_2^\top \Sigma \widehat{\mu}_2} \bigg)^2 \bigg| \widehat{\mu}_2 \bigg\}  ~=~ & \frac{1}{m_1^2} + \mV_P \bigg\{  \frac{\widehat{\mu}_2^\top \widehat{\Sigma} \widehat{\mu}_2}{\widehat{\mu}_2^\top \Sigma \widehat{\mu}_2}\bigg| \widehat{\mu}_2 \bigg\}.
\end{align*}
Using the variance formula of the sample variance~\citep[e.g.][]{lee1990u} and by letting $\widetilde{\mu}_2 = \widehat{\mu}_2 / \sqrt{\widehat{\mu}_2^\top \Sigma \widehat{\mu}_2}$,
\begin{align*}
	\mV_P \bigg\{  \frac{\widehat{\mu}_2^\top \widehat{\Sigma} \widehat{\mu}_2}{\widehat{\mu}_2^\top \Sigma \widehat{\mu}_2}\bigg| \widehat{\mu}_2 \bigg\} ~=~ &   \frac{(m_1-1)^2}{m_1^2} \bigg\{ \frac{\mE_P[ \{ \widetilde{\mu}_2^\top (X - \mu)\}^4 | \widehat{\mu}_2 ]}{m_1}  - \frac{\{\mE_P[(\widetilde{\mu}_2^\top (X - \mu))^2 | \widehat{\mu}_2 ]\}^4(m_1-3)}{m_1(m_1-1)} \bigg\} \\[.5em]
	= ~ & \frac{(m_1-1)^2}{m_1^2} \bigg\{ \frac{3}{m_1} - \frac{(m_1-3)}{m_1(m_1-1)}   \bigg\},
\end{align*}
where we use the following identities:
\begin{align*}
	\mE_P[ \{ \widetilde{\mu}_2^\top (X - \mu)\}^4 | \widehat{\mu}_2 ] ~=~ & \mV_P[ \{ \widetilde{\mu}_2^\top (X - \mu)\}^2 | \widehat{\mu}_2 ]  + \{ \mE_P[ \{ \widetilde{\mu}_2^\top (X - \mu)\}^2 |\widehat{\mu}_2 ] \}^2 \\[.5em]
	= ~ & 2 \tr( \widetilde{\mu}_2 \widetilde{\mu}_2^\top  \Sigma \widetilde{\mu}_2 \widetilde{\mu}_2^\top  \Sigma ) + \tr^2(\widetilde{\mu}_2 \widetilde{\mu}_2^\top  \Sigma) = 3,
\end{align*}
and 
\begin{align*}
	\mE_P[ \{ \widetilde{\mu}_2^\top (X - \mu) \}^2 | \widehat{\mu}_2 ] ~=~  \tr ( \widetilde{\mu}_2^\top \Sigma \widetilde{\mu}_2) = 1.
\end{align*}
Combining the results yields
\begin{align*}
	\mE_P \bigg[ \mE_P \bigg\{ \bigg( \frac{\widehat{\mu}_2^\top \widehat{\Sigma} \widehat{\mu}_2 - \widehat{\mu}_2^\top \Sigma \widehat{\mu}_2}{\widehat{\mu}_2^\top \Sigma \widehat{\mu}_2} \bigg)^2 \bigg| \widehat{\mu}_2 \bigg\} \bigg] ~=~  \frac{1}{m_1^2} +  \frac{2(m_1-1)}{m_1^2} \rightarrow 0 \quad \text{as $m_1 \rightarrow \infty$,}
\end{align*}
which verifies (C3). This completes the proof of Theorem~\ref{Theorem: asymptotic power}.

\subsection{Proof of Theorem~\ref{Theorem: Uniform type I and II error control}}

The type I error result clearly follows by Theorem~\ref{Theorem: unconditional BE Bound}. Therefore we only need to show the type II error result. Note that under Assumption~\ref{assumption: moment condition}, the denominator of $T_{\mathrm{mean}}$ is greater than zero almost surely. Consequently we can reformulate the type II error as 
\begin{align} \label{Eq: type II error}
	\mE_{P}[1 - \phi_{\mathrm{mean}}]  ~=~ P \bigg\{ \widehat{\mu}_2^\top \widehat{\mu}_1  \leq z_{1-\alpha} \times \bigg( m_1^{-2}\sum_{i=1}^{m_1} \{\widehat{\mu}_2^\top (X_i - \widehat{\mu}_1 ) \}^2\bigg)^{1/2} \bigg\},
\end{align}
from which it suffices to check sufficient condition~(\ref{Eq: Sufficient Condition for Power}) in Lemma~\ref{Lemma: Sufficient conditions for uniform type II error control} to ensure type II error control. First note that the kernel considered in Theorem~\ref{Theorem: Uniform type I and II error control} is $h(x,y) =x^\top y$. Given this kernel, let us compute the conditional expectation and the conditional variance in sufficient condition~(\ref{Eq: Sufficient Condition for Power}) as
\begin{align*}
	& \mE_P[\mV_P[X_1^\top X_2|X_2]] = \mE_P[X_2^\top \Sigma X_2] = \mE_P[\tr(X_2X_2^\top \Sigma)] = \tr(\Sigma^2) + \mu^\top \Sigma \mu \quad \text{and} \\[.5em]
	& \mV_P[\mE_P[X_1^\top X_2|X_2]] = \mV_P[\mu^\top X_2] = \mu^\top \Sigma \mu.
\end{align*}
Also note that $\mE_P[h(X_1,X_2)] = \|\mu\|_2^2$. Therefore, the considered test $\phi_{\mathrm{mean}}$ has the uniform type II error converging to zero provided that 
\begin{align}  \label{Eq: intermediate condition}
	\|\mu\|_2^2 \geq \psi_n \max\{z_{1-\alpha}, 1\} \max \Bigg\{ & \sqrt{\frac{\tr(\Sigma^2) + \mu^\top \Sigma \mu}{m_1m_2}}, \ \sqrt{\biggl( \frac{1}{m_1} + \frac{1}{m_2} \biggr) \mu^\top \Sigma \mu} \Bigg\},
\end{align}
where $\psi_n$ is an arbitrary positive sequence that goes to infinity as $n \rightarrow \infty$. Since we assume that (i)~$\alpha$ is fixed, (ii)~the maximum eigenvalue of $\Sigma$ is bounded by $\lambda^\ast$ and (iii)~the splitting ratio satisfies $\kappa_\ast \leq m_1/n \leq \kappa^\ast$, the above condition~(\ref{Eq: intermediate condition}) is implied by
\begin{align*}
	\|\mu\|_2 \geq C(\alpha, \lambda^\ast, \kappa_\ast, \kappa^\ast) \max\{\psi_n, \sqrt{\psi_n}\} \times  \frac{d^{1/4}}{n^{1/2}}, 
\end{align*}
where $C(\alpha, \lambda^\ast, \kappa_\ast, \kappa^\ast)$ is some positive constant which only depends on $\alpha, \lambda^\ast, \kappa_\ast, \kappa^\ast$. Now the desired result follows by choosing $\psi_n$ such that $K_n \geq C(\alpha, \lambda^\ast, \kappa_\ast, \kappa^\ast) \max\{\psi_n, \sqrt{\psi_n}\}$. This completes the proof of Theorem~\ref{Theorem: Uniform type I and II error control}.

\subsection{Proof of Corollary~\ref{Corollary: properties of phi_sym}}
The first claim follows by our discussion in the main text. To show the second claim, note that $1 - \Phi(t) \leq \exp(-t^2/2)$ for all $t \in \mathbb{R}$, which implies $\sqrt{2\log(1/\alpha)} \geq z_{1-\alpha}$ for any $\alpha \in (0,1)$. This means that $\phi_{\mathrm{sym}} \leq \phi_{\mathrm{mean}}$ and thus $\lim_{n \rightarrow \infty} \sup_{P \in \mathcal{P}_{0,d}} \mE_{P} [\phi_{\mathrm{sym}}] \leq \alpha$ due to Theorem~\ref{Theorem: Uniform type I and II error control}. Note that we assume $\alpha$ is fixed and the proof for the uniform type II error of $\phi_{\mathrm{mean}}$ goes through even when $z_{1-\alpha}$ is replaced by some other positive constant. Therefore the second claim also holds.

\subsection{Proof of Theorem~\ref{Theorem: Uniform type I and II error control for sparse alternatives}}

\noindent \textbf{$\bullet$ Type I error.}  Since we assume that $X$ has a continuous distribution, $\widehat{v}$ is nonzero with probability one.  Thus, for any $P \in \mathcal{P}_{0,d}$, we have
\begin{align*}
	\frac{\mE_P[|f_{\mathrm{sparse}}^3(X)| | \mathcal{X}_2]}{\{\mE_P[f_{\mathrm{sparse}}^2(X)|\mathcal{X}_2]\}^{3/2} \sqrt{m_1}} = \frac{\mE_P[|\widehat{v}^\top X|^3 | \mathcal{X}_2]}{\{\mE_P[(\widehat{v}^\top X)^2|\mathcal{X}_2]\}^{3/2} \sqrt{m_1}} \leq \frac{L}{\sqrt{m_1}} \quad \text{a.s.}
\end{align*}
Therefore the same Berry--Esseen bound used earlier gives us
\begin{align*}
	\sup_{P \in \mathcal{P}_{0,d}} \sup_{-\infty < z < \infty} \big| P(T_{\mathrm{sparse}}  \leq z) - \Phi(z) \big| \leq \frac{C L}{\sqrt{m_1}},
\end{align*}
which in turn implies that type I error control:
\begin{align*}
	\sup_{P \in \mathcal{P}_{0,d}} \big| \mE_P[\phi_{\mathrm{sparse}}] -\alpha \big| = O\left(\frac{1}{\sqrt{m_1}}\right) = o(1).
\end{align*}

\vskip 1em

\noindent \textbf{$\bullet$ Type II error.} Turning to type II error control, assume without loss of generality that $|\mu^{(j)}| > \epsilon_n$ for $j \in \mathcal{J}:=\{1,\ldots,J\}$ where $1 \leq J \leq d$. For each $j \in \mathcal{J}$, let us define an event $\mathcal{A}_j := \big\{ \widehat{\mu}_2^{(j)} > \max_{i \in \mathcal{J}^c} |\widehat{\mu}_2^{(i)}| \big\}$ if $\mu^{(j)} >0$ and $\big\{ -\widehat{\mu}_2^{(j)} > \max_{i \in \mathcal{J}^c} |\widehat{\mu}_2^{(i)}| \big\}$ if $\mu^{(j)} <0$. 
%\begin{align*}
%	\mathcal{A}_j := \big\{ |\widehat{\mu}_j| > \max_{1 \leq i \neq j \leq d} |\widehat{\mu}_i| \big\}.
%\end{align*}
%\begin{align*}
%	\mathcal{A}_j := 
%	\begin{cases}
%		\big\{ \widehat{\mu}^{(j)} > \max_{i \in \mathcal{J}^c} |\widehat{\mu}^{(i)}| \big\}, \quad & \text{if $\mu_j >0$} \\[.1em]
%		\big\{ -\widehat{\mu}^{(j)} > \max_{i \in \mathcal{J}^c} |\widehat{\mu}^{(i)}| \big\}, \quad & \text{if $\mu_j<0$.}
%	\end{cases}
%\end{align*}
Then the type II error of our test is bounded above by
\begin{align} \nonumber
	P (T_\mathrm{sparse} \leq z_{1-\alpha})  ~\leq~ & \underbrace{P (T_\mathrm{sparse} \leq z_{1-\alpha}, \cup_{j \in \mathcal{J}} \mathcal{A}_j )}_{(I)} + \underbrace{P(\cap_{j \in \mathcal{J}} \mathcal{A}_j^c)}_{(II)}.
\end{align}
We analyze $(I)$ and $(II)$ separately. To start with the term $(I)$, the event $\cup_{j \in \mathcal{J}} \mathcal{A}_j$ implies that $i^\star = \arg\max_{1\leq j \leq d} |\widehat{\mu}_{2}^{(j)}|$ should belong to $\mathcal{J}$. Thus by letting events $\mathcal{B}_j := \{i^\star =j\}$ for $j=1,\ldots,J$, which are disjoint with probability one, we have 
\begin{align} \label{Eq: conditional event}
	P (T_\mathrm{sparse} \leq z_{1-\alpha}, \cup_{j \in \mathcal{J}} \mathcal{A}_j) ~=~ & \sum_{k \in \mathcal{J}} P (T_\mathrm{sparse} \leq z_{1-\alpha}, \cup_{j \in \mathcal{J}} \mathcal{A}_j, \mathcal{B}_{k}) \\[.5em] \nonumber 
	= ~ &  \sum_{k \in \mathcal{J}} P(\mathcal{B}_{k},\cup_{j \in \mathcal{J}} \mathcal{A}_j) P(T_{\mathrm{sparse}} \leq z_{1-\alpha} ~| \cup_{j \in \mathcal{J}} \mathcal{A}_j, \mathcal{B}_{k}).
\end{align}
Now, by letting $\widehat{s}_{1,k}^2= m_1^{-1} \sum_{i=1}^{m_1} \{f_{\mathrm{sparse}}(X_i) - \overline{f}_{\mathrm{sparse}}\}^2$ when $i^\star = k$ and assuming $\mu^{(k)}>0$ (and the same argument holds for $\mu^{(k)}<0$), we can simply write 
\begin{align*}
	& P(T_{\mathrm{sparse}} \leq z_{1-\alpha} ~| \cup_{j \in \mathcal{J}} \mathcal{A}_j, \mathcal{B}_k) \\[.5em]
	= ~ & P \bigg\{\text{sign}\Big(\widehat{\mu}_{2}^{(k)} \Big) \times  \frac{\sqrt{m_1} \widehat{\mu}_{1}^{(k)}}{\widehat{s}_{1,k}} \leq z_{1-\alpha} ~\Big|  \cup_{j \in \mathcal{J}} \mathcal{A}_j, \mathcal{B}_k \bigg\} \\[.5em]
	= ~ & P \bigg\{ \frac{\sqrt{m_1}( \widehat{\mu}_{1}^{(k)} - \mu^{(k)}) }{\widehat{s}_{1,k}} \leq z_{1-\alpha} - \ \frac{\sqrt{m_1}\mu^{(k)}}{\widehat{s}_{1,k}} ~\Big|  \cup_{j \in \mathcal{J}} \mathcal{A}_j, \mathcal{B}_k \bigg\}.
\end{align*}
Define another event $\mathcal{C}_k := \{ \widehat{s}_{1,k} < \nu^{-1} \sigma \}$ where $\nu>0$ will be specified later. Notice that the sub-Gaussian parameter $\sigma^2$ is lower bounded by the variance of $X^{(j)}$ for all $j=1,\ldots,d$. Then, by Markov's inequality, it holds $P\{\mathcal{C}_k^c ~| \cup_{j \in \mathcal{J}} \mathcal{A}_j, \mathcal{B}_k \} \leq \nu$ and thus we may further bound the above equation by
\begin{align*}
	& P(T_{\mathrm{sparse}} \leq z_{1-\alpha} ~| \cup_{j \in \mathcal{J}} \mathcal{A}_j, \mathcal{B}_k) \\[.5em]
	\leq ~ &  P \bigg\{ \frac{\sqrt{m_1}( \widehat{\mu}_{1}^{(k)} - \mu^{(k)}) }{\widehat{s}_{1,k}} \leq z_{1-\alpha} - \ \frac{\sqrt{m_1}\epsilon_n}{\nu^{-1}\sigma} ~\Big|  \cup_{j \in \mathcal{J}} \mathcal{A}_j, \mathcal{B}_k \bigg\} + \nu \\[.5em]
	\leq ~ & \Phi \bigg\{ z_{1-\alpha} - \ \frac{\sqrt{m_1}\epsilon_n}{\nu^{-1}\sigma}  \bigg\} + \frac{C}{\sqrt{m_1}} + \nu,
\end{align*}
where the first inequality uses our assumption that $\mu^{(k)} \geq \epsilon_n >0$ and the second inequality uses the Berry--Esseen bound for the studentized statistic~\citep{bentkus1996berry}. Here $C$ is a positive constant depending on $L$ given in Assumption~\ref{assumption: moment condition}. Combining this bound with the equation~(\ref{Eq: conditional event}) yields
\begin{align*}
	P (T_\mathrm{sparse} \leq z_{1-\alpha}, \cup_{j \in \mathcal{J}} \mathcal{A}_j) ~\leq~ & \Phi \bigg\{ z_{1-\alpha} - \ \frac{\sqrt{m_1}\epsilon_n}{\nu^{-1}\sigma}  \bigg\} + \frac{C}{\sqrt{m_1}} + \nu \\[.5em]
	\leq ~ & \Phi \bigg\{ z_{1-\alpha} - \ \frac{\sqrt{m_1}}{\nu^{-1}\sigma} \times K_n \sigma\sqrt{\frac{\log d}{n}}  \bigg\} + \frac{C}{\sqrt{m_1}}  + \nu,
\end{align*}
where the last inequality holds due to $\epsilon_n \geq K_n \sigma \sqrt{n^{-1}\log d}$. Note that the upper bound does not depend on $P$. Therefore, by taking $\nu = K_n^{-q}$ where $0 < q < 1$ and assuming $\kappa_\ast < m_1/n < \kappa^\ast$, we see that 
\begin{align} \label{Eq: bound for the first term}
	\sup_{P \in \mathcal{Q}_{\mathrm{sparse},d}(\epsilon_n)} P (T_\mathrm{sparse} \leq z_{1-\alpha}, \cup_{j \in \mathcal{J}} \mathcal{A}_j) \rightarrow 0 \quad \text{as $n \rightarrow \infty$.}
\end{align}

Next turning to the term $(II)$, we first bound $(II)$ by $P(\cap_{j \in \mathcal{J}} \mathcal{A}_j^c) \leq P(\mathcal{A}_1^c)$. For the case when $\mu^{(1)} >0$ (and similarly for $\mu^{(1)}<0$), 
\begin{align*}
	P(\mathcal{A}_1^c) ~=~ & P \big\{ \widehat{\mu}_2^{(1)} \leq \max_{i \in \mathcal{J}^c} |\widehat{\mu}_2^{(i)}| \big\}  = P \big\{ \mu^{(1)}  \leq \max_{i \in \mathcal{J}^c} |\widehat{\mu}_2^{(i)}| - (\widehat{\mu}_2^{(1)} - \mu^{(1)}) \big\} \\[.5em]
	\leq ~ & P \big\{ \mu^{(1)}  \leq \max_{i \in \mathcal{J}^c} |\widehat{\mu}_2^{(i)}| + |\widehat{\mu}_2^{(1)} - \mu^{(1)}| \big\} \\[.5em]
	\leq ~ & P \big\{ \mu^{(1)} \leq 2 \max_{1\leq j \leq d} |\widehat{\mu}_2^{(i)} - \mu^{(i)}| \big\}.
\end{align*}
We note that $\widehat{\mu}_2^{(i)} - \mu^{(i)}$ is sub-Gaussian with parameter $\sigma^2/m_2$. Therefore, together with the union bound, 
\begin{align*}
	P \big\{ \mu^{(1)} \leq 2 \max_{1\leq j \leq d} |\widehat{\mu}_2^{(i)} - \mu^{(i)}| \big\} ~ \leq ~ & \sum_{i=1}^d P \{ \mu^{(1)}/2 \leq |\widehat{\mu}_2^{(i)} - \mu^{(i)}|  \} \\[.5em]
	\leq ~ & 2d \exp \bigg\{ - \frac{m_2 (\mu^{(1)})^2}{8\sigma^2} \bigg\},
\end{align*}
Now by our choice of $\mu^{(1)} \geq \epsilon_n$, the upper bound converges to zero uniformly over $\mathcal{Q}_{\mathrm{sparse},d}(\epsilon_n)$ and so 
\begin{align*}
	\sup_{P \in \mathcal{Q}_{\mathrm{sparse},d}(\epsilon_n)} P(\cap_{j \in \mathcal{J}} \mathcal{A}_j^c) \rightarrow 0 \quad \text{as $n \rightarrow \infty$.}
\end{align*}
Combining this with (\ref{Eq: bound for the first term}) completes the proof of Theorem~\ref{Theorem: Uniform type I and II error control for sparse alternatives}.

\subsection{Proof of Proposition~\ref{Proposition: BE bound for T2}}
We start by noting that $f_{\mathrm{cov}}(x)$ can be written as
\begin{align*}
	f_{\mathrm{cov}}(x) = x^\top (S - I)x - \mE_{P_{0}}[X^\top (S - I) X | \mathcal{X}_2],
\end{align*}
where $S = m_2^{-1} \sum_{j=m_1+1}^n X_j X_j^\top$. Based on this reformulation and using the moment expression for a Gaussian quadratic form \citep[e.g. Lemma 3 of][]{cai2013optimal}, we have
\begin{align*}
	& \mE_{P_{0}}[|f_{\mathrm{cov}}(X)|^2 | \mathcal{X}_{2}] =  \mV_{P_{0}}[f_{\mathrm{cov}}(X)|  \mathcal{X}_{2}] = 2\tr \{ (S - I)^2 \}, \\[.5em]
	& \mE_{P_{0}}[|f_{\mathrm{cov}}(X)|^4|\mathcal{X}_{2}] = 48 \tr\{ (S - I)^4 \}  + 12 \tr^2 \{ (S - I)^2 \}.
\end{align*} 
Therefore the right-hand side of (\ref{Eq: BE bound for covariance 2}) becomes
\begin{align*}
	\frac{\mE_{P_{0}}[|f_{\mathrm{cov}}(X)|^4|\mathcal{X}_{2}]}{\{ \mE_{P_{0}}[|f_{\mathrm{cov}}(X)|^2 | \mathcal{X}_{2}] \}^2 m_1} ~  = ~   \frac{12 \tr\{ (S - I)^4 \}}{m_1 \tr^2 \{ (S - I)^2 \}} + \frac{3}{m_1} \quad \text{a.s.}
\end{align*}
Furthermore, by Cauchy-Schwarz inequality, we have $\tr\{ (S - I)^4\} \leq \tr^2 \{ (S - I)^2 \}$, which leads to
\begin{align}
	\frac{\mE_{P_{0}}[|f_{\mathrm{cov}}(X)|^4|\mathcal{X}_{2}]}{\{ \mE_{P_{0}}[|f_{\mathrm{cov}}(X)|^2 | \mathcal{X}_{2}] \}^2 m_1} \leq \frac{15}{m_1} \quad \text{a.s.}
\end{align}
For a general distribution $P$ in $\mathcal{P}_{0,d,\mathrm{cov}}$, Proposition A.1 of \cite{chen2010tests} yields that 
\begin{align*}
	& \mE_{P}[|f_{\mathrm{cov}}(X)|^2 | \mathcal{X}_{2}] \geq  c \tr \{ (S - I)^2 \}, \\[.5em]
	& \mE_{P}[|f_{\mathrm{cov}}(X)|^4|\mathcal{X}_{2}] \leq  C \tr^2 \{ (S - I)^2 \},
\end{align*}
where constant $c$ is given in Definition~\ref{Definition: covariance model} and $C$ is another constant given in Proposition A.1 of \cite{chen2010tests}. Therefore, for some constant $C'>0$, it holds
\begin{align}
	\sup_{P \in \mathcal{P}_{0,d,\mathrm{cov}}} \frac{\mE_{P}[|f_{\mathrm{cov}}(X)|^4|\mathcal{X}_{2}]}{\{ \mE_{P}[|f_{\mathrm{cov}}(X)|^2 | \mathcal{X}_{2}] \}^2 m_1} \leq \frac{C'}{m_1} \quad \text{a.s.}
\end{align}
This result together with the Berry--Esseen bound~(\ref{Eq: BE bound for covariance}) completes the proof of Proposition~\ref{Proposition: BE bound for T2}.

\subsection{Proof of Theorem~\ref{Theorem: asymptotic power of covariance test}}
To study the power of $\phi_{\mathrm{cov}}$, we need to understand the limiting behavior of $T_{\mathrm{cov}}$ under the alternative. To this end, let us modify $f_{\mathrm{cov}}(x)$ in (\ref{Eq: g function}) and define
\begin{align*}
	f_{\mathrm{cov}}^\ast(x) := f_{\mathrm{cov}}(x) - \tr\{ (\Sigma - I) (S - I) \},
\end{align*} 
where $S = m_2^{-1} \sum_{j=m_1+1}^n X_j X_j^\top$. By putting $\overline{f}_{\mathrm{cov}}^\ast = m_1^{-1}\sum_{i=1}^{m_1} f_{\mathrm{cov}}^\ast(X_i)$, consider a studentized statistic
\begin{align*}
	T_{\mathrm{cov}}^{(a)} := \frac{\sqrt{m_1} \overline{f}_{\mathrm{cov}}^\ast}{\sqrt{m_1^{-1} \sum_{i=1}^{m_1} \{ f_{\mathrm{cov}}^\ast(X_i) - \overline{f}_{\mathrm{cov}}^\ast\}^2 }} = T_{\mathrm{cov}} - \frac{\sqrt{m_1} \tr\{ (\Sigma - I) (S - I) \}}{\sqrt{m_1^{-1} \sum_{i=1}^{m_1} \{ f_{\mathrm{cov}}^\ast(X_i) - \overline{f}_{\mathrm{cov}}^\ast\}^2 }}.
\end{align*}
In the rest of this proof, we prove the following claims:
\begin{enumerate}\setlength{\itemindent}{0.073in}
	\item[\textbf{(D1).}] $\displaystyle T_{\mathrm{cov}}^{(a)} \convD N(0,1).$
	\item[\textbf{(D2).}] $\displaystyle \frac{m_1^{-1} \sum_{i=1}^{m_1} \{ f_{\mathrm{cov}}^\ast(X_i) - \overline{f}_{\mathrm{cov}}^\ast \}^2}{2\tr \{ (\Sigma(S - I))^2\}}= 1 + o_P(1)$.
	\item[\textbf{(D3).}] $\displaystyle \frac{\sqrt{m_1} \tr\{ (\Sigma - I) (S - I) \}}{\sqrt{2\tr \{ (\Sigma(S - I))^2\}}} = \frac{n\sqrt{\kappa(1-\kappa)}\|\Sigma - I\|_F^2}{\sqrt{2}{\tr}(\Sigma^2)} + o_P(1) = O_P(1)$.
\end{enumerate}
Once these claims are proved, the power function of $\phi_{\mathrm{cov}}$ approximates
\begin{align*}
	\mE_P[\phi_{\mathrm{cov}}] ~=~ & P\left( T_{\mathrm{cov}}^{(a)} > z_{1-\alpha} - \frac{\sqrt{m_1} \tr\{ (\Sigma - I) (S - I) \}}{\sqrt{m_1^{-1} \sum_{i=1}^{m_1} \{ f_{\mathrm{cov}}(X_i) - \overline{f}_{\mathrm{cov}}\}^2 }} \right) \\[.5em]
	=~ & P\left( T_{\mathrm{cov}}^{(a)} > z_{1-\alpha} - \frac{\sqrt{m_1} \tr\{ (\Sigma - I) (S - I) \}}{\sqrt{2\tr \{ (\Sigma(S - I))^2\} }}\{1+o_P(1)\} \right) \\[.5em]
	=~ & \Phi \left( z_\alpha + \frac{n\sqrt{\kappa(1-\kappa)}\|\Sigma - I\|_F^2}{\sqrt{2}{\tr}(\Sigma^2)} \right) + o(1),
\end{align*}
which completes the proof. The rest of the proof is dedicated to verifying these claims. 

\vskip 1em

\noindent \textbf{$\bullet$ Verification of (D1).} Since $\mE_P[f_{\mathrm{cov}}^\ast(X)|\mathcal{X}_2] = 0$, the conditional Berry--Esseen bound for a studentized statistic implies
\begin{align*}
	\sup_{-\infty < z < \infty} \big| P(T_{\mathrm{cov}}^{(a)}  \leq z | \mathcal{X}_{2}) - \Phi(z) \big| \leq  C_1 \frac{\mE_P[|f_{\mathrm{cov}}^\ast(X)|^3|\mathcal{X}_{2}]}{\{\mE_P[|f_{\mathrm{cov}}^\ast(X)|^2 | \mathcal{X}_{2}] \}^{3/2} \sqrt{m_1}},
\end{align*}
for some $C_1 >0$. Note that $f_{\mathrm{cov}}^\ast(X)$ can be succinctly written as
\begin{align*}
	f_{\mathrm{cov}}^\ast(x) = x^\top (S - I) x - \mE_{P}[X^\top (S - I)X|\mathcal{X}_2]
\end{align*} 
and its second and fourth moments can be computed as
\begin{equation}
	\begin{aligned}   \label{Eq: fourth moment of g1}
		& \mE_P[|f_{\mathrm{cov}}^\ast(X)|^2 | \mathcal{X}_{2}] =  \mV_P[f_{\mathrm{cov}}^\ast(X)|  \mathcal{X}_{2}] = 2\tr \{ (\Sigma(S - I))^2\}, \\[.5em]
		& \mE_P[|f_{\mathrm{cov}}^\ast(X)|^4|\mathcal{X}_{2}] = 48 \tr\{ (\Sigma (S - I))^4 \}  + 12 \tr^2 \{ (\Sigma(S - I))^2 \}.
	\end{aligned} 
\end{equation}
Therefore, similarly to Proposition~\ref{Proposition: BE bound for T2}, 
\begin{align*}
	\sup_{-\infty < z < \infty} \big| P(T_{\mathrm{cov}}^{(a)}  \leq z) - \Phi(z) \big| \leq  \frac{C_2}{\sqrt{m_1}},
\end{align*}
which shows that $T_{\mathrm{cov}}^{(a)}$ converges to the standard normal distribution. 

\vskip 1em

\noindent \textbf{$\bullet$ Verification of (D2).} Our goal is to prove that for any $\epsilon > 0$
\begin{align} \label{Eq: goal of D2}
	P \left( \bigg| \frac{m_1^{-1} \sum_{i=1}^{m_1} \{ f_{\mathrm{cov}}^\ast(X_i) - \overline{f}_{\mathrm{cov}}^\ast \}^2}{2\tr \{ (\Sigma(S - I))^2\}} - 1 \bigg| \geq \epsilon  \right) \rightarrow 0 \quad \text{as $m_1 \rightarrow \infty$.}
\end{align}
To verify this claim, we follow the same steps used in the proof of Theorem~\ref{Theorem: asymptotic power}. By the law of total expectation and conditional Chebyshev's inequality,
\begin{align*}
	& P \bigg( \bigg|  \frac{m_1^{-1} \sum_{i=1}^{m_1} \{ f_{\mathrm{cov}}^\ast(X_i) - \overline{f}_{\mathrm{cov}}^\ast \}^2}{2\tr \{ (\Sigma(S - I))^2\}} - 1\bigg|  > \epsilon \bigg) \\[.5em]
	\leq ~ & \frac{1}{\epsilon^2} \mE_P \bigg[ \mE_P \bigg\{ \bigg( \frac{m_1^{-1} \sum_{i=1}^{m_1} \{ f_{\mathrm{cov}}^\ast(X_i) - \overline{f}_{\mathrm{cov}}^\ast \}^2 -2\tr \{ (\Sigma(S - I))^2\}}{2\tr \{ (\Sigma(S - I))^2\}} \bigg)^2 \bigg| \mathcal{X}_2 \bigg\} \bigg]
\end{align*} 
and the bias-variance decomposition along with the variance formula of the sample variance~\citep[e.g.][]{lee1990u} gives
\begin{align*}
	& \mE_P \bigg\{ \bigg( \frac{m_1^{-1} \sum_{i=1}^{m_1} \{ f_{\mathrm{cov}}^\ast(X_i) - \overline{f}_{\mathrm{cov}}^\ast \}^2 -2\tr \{ (\Sigma(S - I))^2\}}{2\tr \{ (\Sigma(S - I))^2\}} \bigg)^2 \bigg| \mathcal{X}_2 \bigg\} \\[.5em]
	=~ & \frac{1}{m_1^2} + \mV_P \bigg\{  \frac{m_1^{-1} \sum_{i=1}^{m_1} \{ f_{\mathrm{cov}}^\ast(X_i) - \overline{f}_{\mathrm{cov}}^\ast \}^2}{2\tr \{ (\Sigma(S - I))^2\}}\bigg| \mathcal{X}_2 \bigg\} \\[.5em]
	=~ & \frac{1}{m_1^2} + \frac{(m_1-1)^2}{m_1^2} \Bigg\{ \frac{\mE_P[|f_{\mathrm{cov}}^\ast(X_1)|^4 | \mathcal{X}_2]}{4 m_1\tr^2 \{ (\Sigma(S - I))^2\}}  - \frac{(m_1-3)}{m_1(m_1-1)} \Bigg\}.
\end{align*}
Now leveraging the expression~(\ref{Eq: fourth moment of g1}) and Cauchy-Schwarz inequality yield 
\begin{align*}
	\frac{\mE_P[|f_{\mathrm{cov}}^\ast(X_1)|^4 | \mathcal{X}_2]}{4 m_1\tr^2 \{ (\Sigma(S - I))^2\}}  ~\leq~ \frac{15}{m_1} \quad \text{a.s.}
\end{align*}
Putting things together, we see that the convergence result~(\ref{Eq: goal of D2}) holds and thus (D2) follows.

\vskip 1em

\noindent \textbf{$\bullet$ Verification of (D3).} We begin by analyzing the numerator of 
\begin{align} \label{Eq: Term in (D3)}
	\frac{\sqrt{m_1} \tr\{ (\Sigma - I) (S - I) \}}{\sqrt{2\tr \{ (\Sigma(S - I))^2\}}}. 
\end{align}
The first two moments of the numerator are computed as
\begin{align*}
	\mE_P[ \tr\{ (\Sigma - I) (S - I) \}]  ~=~ & \|\Sigma - I \|_F^2, \\[.5em]
	\mV_P[ \tr\{ (\Sigma - I) (S - I) \}] ~=~ & \frac{2}{m_2} \tr \{ (\Sigma - I) \Sigma (\Sigma - I) \Sigma\} \\[.5em]
	~= ~& \frac{2}{m_2} \sum_{i=1}^d \lambda_i^2(\Sigma) \{ \lambda_i(\Sigma) - 1 \}^2 \\[.5em]
	~\leq~ & \frac{2 \lambda_{d}^2(\Sigma)}{m_2} \| \Sigma - I \|_F^2, 
\end{align*}
where we recall that $\lambda_1(\Sigma) \leq \ldots \leq \lambda_d(\Sigma)$ are eigenvalues of $\Sigma$. Therefore
\begin{align*}
	\tr\{ (\Sigma - I) (S - I) \} & ~=~  \|\Sigma - I \|_F^2 + O_P\left( \! \sqrt{\frac{2 \lambda_d^2(\Sigma)}{m_2}} \| \Sigma - I \|_F \right) \\[.5em] 
	& ~=~  \|\Sigma - I \|_F^2 + O_P\big(m_2^{-1/2}\big),
\end{align*}
where the second equality follows since the eigenvalues are uniformly bounded and $\| \Sigma - I \|_F = O(1)$. Next we turn to the denominator of the quantity~(\ref{Eq: Term in (D3)}). First of all, we compute the expectation
\begin{align*}
	\mE_P[ 2\tr \{ (\Sigma(S - I))^2\} ] ~=~ & 2 \mE_P[ \tr\{ (\Sigma(S - \Sigma))^2 \} ] + 4 \mE [\tr\{ \Sigma(\Sigma - I) \Sigma (S-\Sigma) \}] \\[.5em]
	+ ~ & 2\sum_{i=1}^d \lambda_i^2(\Sigma) \{ \lambda_i(\Sigma) - 1 \}^2 \\[.5em] \label{Eq: expectation of the denominator}
	= ~ & \frac{2\tr(\Sigma^4)}{m_2}  + \frac{2\tr^2(\Sigma^2)}{m_2} + 2\sum_{i=1}^d \lambda_i^2(\Sigma) \{\lambda_i(\Sigma) - 1\}^2,
\end{align*}
where the second term asymptotically dominates the other terms under the given conditions. It can be seen that the variance is bounded by
\begin{align*}
	\mV_P[ 2\tr \{ (\Sigma(S - I))^2\} ] ~ \lesssim ~ & \mV_P[ \tr\{ (\Sigma(S - \Sigma))^2 \} ] + \mV_P[ \tr\{ \Sigma(\Sigma - I) \Sigma (S-\Sigma) \} ] \\[.5em]
	\lesssim ~ &  \frac{1}{m_2^3} \Big\{ \tr(\Sigma^8) + \tr^2(\Sigma^4) + \tr^2(\Sigma^2)\tr(\Sigma^4)  \Big\} \\[.5em]
	+ ~& \frac{1}{m_2^2} \Big\{ \tr^2(\Sigma^8) + \tr(\Sigma^{16})   \Big\} \\[.5em]
	+ ~& \frac{1}{m_2} \tr(\Sigma^8) + \frac{1}{m_2} \tr \{(\Sigma^3(\Sigma - I) )^2\}.
\end{align*}
Using the bounded eigenvalue condition along with the conditions on $\kappa$ and $\tau$, the variance is further upper bounded by
\begin{align*}
	\mV_P[ 2\tr \{ (\Sigma(S - I))^2\} ] ~ \lesssim ~  \frac{d^3}{m_2^3} + \frac{d^2}{m_2^2} + \frac{d}{m_2} = O(1).
\end{align*}
To summarize the results, we have obtained that 
\begin{align*}
	\tr\{ (\Sigma - I) (S - I) \}  ~=~   \|\Sigma - I \|_F^2 + O_P\big(m_2^{-1/2}\big)
\end{align*}
and 
\begin{align*}
	2m_1^{-1}\tr \{ (\Sigma(S - I))^2\} ~=~ &  \frac{2\tr(\Sigma^4)}{m_1m_2}  + \frac{2\tr^2(\Sigma^2)}{m_1m_2} + \frac{2}{m_1}\sum_{i=1}^d \lambda_i^2(\Sigma) \{\lambda_i(\Sigma) - 1\}^2 + O_P \big( m_1^{-1} \big) \\[.5em]
	=~ &  \frac{2 \tr^2(\Sigma^2)}{n^2\kappa(1-\kappa)} + O_P \big(n^{-1} \big).
\end{align*}
Therefore the quantity~(\ref{Eq: Term in (D3)}) approximates
\begin{align*}
	\frac{\sqrt{m_1} \tr\{ (\Sigma - I) (S - I) \}}{\sqrt{2\tr \{ (\Sigma(S - I))^2\}}} = \frac{n \sqrt{\kappa(1-\kappa)}\|\Sigma - I \|_F^2}{\sqrt{2} \tr(\Sigma^2)} + o_P(1).
\end{align*}
Moreover the first term of the right-hand side is bounded, which verifies (D3). This completes the proof of Theorem~\ref{Theorem: asymptotic power of covariance test}.

\subsection{Proof of Proposition~\ref{Proposition: Uniform type I and II error control of the covariance test}}

The proof follows similar lines of argument to the proof of Theorem~\ref{Theorem: Uniform type I and II error control}. First of all, type I error control is a direct consequence of Proposition~\ref{Proposition: BE bound for T2} and so we focus on type II error control. Due to the continuity, the type II error rate can be written as
\begin{align*}
	\mE_{P}[1 - \phi_{\mathrm{cov}}] = P \Bigg( \overline{f}_{\mathrm{cov}} \leq z_{1-\alpha} \sqrt{\frac{1}{m_1^2} \sum_{i=1}^{m_1} \big\{ f_{\mathrm{cov}}(X_i) - \overline{f}_{\mathrm{cov}} \big\}^2 } \Bigg),
\end{align*}
from which it suffices to check sufficient condition~(\ref{Eq: Sufficient Condition for Power}) in Lemma~\ref{Lemma: Sufficient conditions for uniform type II error control} to ensure type II error control. First note that the kernel considered in Proposition~\ref{Proposition: Uniform type I and II error control of the covariance test} is $h(x,y) =\tr\{(xx^\top - I)(yy^\top - I)\}$. Given this kernel, let us compute the conditional expectation and the conditional variance in sufficient condition~(\ref{Eq: Sufficient Condition for Power}). In fact, the conditional expectation can be bounded as
\begin{align*}
	\mE_P[\mV_P[h(X_1,X_2)|X_2]]  ~ \overset{\mathrm{(i)}}{\leq} ~ & \mV_P[h(X_1,X_2)] \\[.5em] \overset{\mathrm{(ii)}}{=} ~& 2 \tr^2(\Sigma^2) + 2 \tr(\Sigma^4) + 4 \tr(\Sigma^2(\Sigma - I)).
\end{align*}
where step~(i) uses the law of total variance and step~(ii) follows by Lemma 4 of \cite{cai2013optimal}. Since $\tr(\Sigma^4) \leq \tr^2(\Sigma^2)$ and $\tr^2(\Sigma^2) \leq 8 \| \Sigma - I\|_F^4 + 8d^2$, which can be verified by applying the basic inequality $(a+b)^2 \leq 2 a^2 + 2 b^2$ for all $a,b \in \mathbb{R}$ twice, we can further bound 
\begin{align*}
	2 \tr^2(\Sigma^2) + 2 \tr(\Sigma^4) + 4 \tr(\Sigma^2(\Sigma - I)) ~\leq~ & 4 \tr^2(\Sigma^2)  + 4 \sum_{i=1}^d \lambda_i^2(\lambda_i - 1)^2 \\[.5em]
	\leq ~ & 32 \|\Sigma - I\|_F^4 + 32 d^2 +  4 \sum_{i=1}^d \lambda_i^2(\lambda_i - 1)^2. 
\end{align*}
On the other hand, by letting $X_1,X_2,X_3 \overset{i.i.d.}{\sim} N(0,\Sigma)$, the conditional variance can be computed as
\begin{align*}
	\mV_P[\mE_P[h(X_1,X_2)|X_2]] ~=~ & \mE_P[\tr^2\{(\Sigma - I) (X_2X_2^\top - \Sigma)\}] \\[.5em]
	= ~ & \mathrm{cov}_P[h(X_2,X_1), h(X_2,X_3)] ~= 2 \tr\{\Sigma^2(\Sigma - I)^2\},
\end{align*}
where the last equality holds due to Lemma~4 of \cite{cai2013optimal}. Moreover, since the maximum eigenvalue of $\Sigma$, namely $\lambda_{\mathrm{max}}$, satisfies $\lambda_{\mathrm{max}} \leq \| \Sigma - I \|_F + 1$ and 
\begin{align*}
	\tr\{\Sigma^2(\Sigma - I)^2\} ~=~ \sum_{i=1}^d \lambda_i^2(\lambda_i - 1)^2 ~\leq~ & \lambda_{\mathrm{max}}^2 \sum_{i=1}^d (\lambda_i -1 )^2  ~=~ \lambda_{\mathrm{max}}^2 \|\Sigma - I\|_F^2 \\[.5em]
	\leq ~ &  2 \| \Sigma - I\|_F^4 + 2 \| \Sigma - I\|_F^2. 
\end{align*}
Hence, we see by using these bounds that sufficient condition~(\ref{Eq: Sufficient Condition for Power}) is fulfilled if
\begin{align*}
	& \mE_P[h(X_1,X_2)] = \| \Sigma - I \|_F^2 \\[.5em]
	\geq ~ & \psi_n \max\{z_{1-\alpha},1\} \max \Bigg\{  \sqrt{\frac{40 \| \Sigma - I\|_F^4 + 8 \|\Sigma - I\|_F^2 + 32d^2}{m_1m_2}}, \\
	&~~~~~~~~~~~~~~~~~~~~~~~~~~~~~~~~~~~~~~~~~~~~~~ \sqrt{\biggl(\frac{4}{m_1} + \frac{4}{m_2} \biggr)(\|\Sigma - I\|_F^4 + \|\Sigma - I\|_F^2)} \Bigg\}.
\end{align*}
Since $\kappa_\ast \leq m_1/n \leq \kappa^\ast$ and $\alpha$ is fixed, this condition holds if the following three conditions are satisfied all together:
\begin{equation}
\begin{aligned} \label{Eq: three conditions}
	& \| \Sigma - I \|_F^2 ~\geq~ C(\alpha,\kappa_\ast,\kappa^\ast) \psi_n \frac{\| \Sigma - I\|_F^2}{n}, \\[.5em]
	&  \| \Sigma - I \|_F^2 ~\geq~ C(\alpha,\kappa_\ast,\kappa^\ast) \psi_n \frac{\| \Sigma - I\|_F}{n}, \\[.5em]
	&   \| \Sigma - I \|_F^2 ~\geq~ C(\alpha,\kappa_\ast,\kappa^\ast) \psi_n \frac{d}{n},
\end{aligned}
\end{equation}
where $C(\alpha,\kappa_\ast,\kappa^\ast)$ is a positive constant that only depends on $\alpha,\kappa_\ast,\kappa^\ast$. Notice that we have freedom of choosing $\psi_n$ and so we choose $\psi_n > 0$ such that 
\begin{align*}
	& K_n \geq C(\alpha,\kappa_\ast,\kappa^\ast)  \psi_n, \ K_n \geq  \sqrt{C(\alpha,\kappa_\ast,\kappa^\ast) \psi_n} \quad \text{and} \\[.5em]
	& 1 \geq  C(\alpha,\kappa_\ast,\kappa^\ast)  \frac{\psi_n}{n}.
\end{align*}
Then under the condition that $\| \Sigma - F \|_F \geq K_n \sqrt{d/n}$, the above three conditions~(\ref{Eq: three conditions}) are fulfilled. Therefore we conclude that 
\begin{align*}
	\lim_{n \rightarrow \infty} \sup_{P \in \mathcal{Q}_{\mathrm{cov},d}(\epsilon_n)} \mE_{P} [1 -  \phi_{\mathrm{cov}}] =0.
\end{align*}
This completes the proof of Proposition~\ref{Proposition: Uniform type I and II error control of the covariance test}.

\subsection{Proof of Theorem~\ref{Theorem: Gaussian Approximation}}
Our strategy to prove Theorem~\ref{Theorem: Gaussian Approximation} is as follows. We begin by considering the conditionally standardized statistic, that is
\begin{align*}
	T_{\text{stand},h} = \frac{\sqrt{m} \overline{f}_h}{\sqrt{\mE_{P}[f_h^2(X_1)|\mathcal{X}_2]}},
\end{align*}
and show that it converges to $N(0,1)$ in distribution. Next we establish that the difference between $T_{\text{stand},h}$ and $T_{h}$ defined in (\ref{Eq: studentized statistic with h}) is asymptotically negligible and hence $T_{h}$ converges to $N(0,1)$ as well. 

In the proof below, the probability measure $P$, kernel $h$ and other associated quantities such as a sequence of eigenvalues $\{\lambda_i\}_{i=1}^\infty$ are allowed to change with $n$. Hence these are implicitly dependent on $n$. For notational convenience, however, the dependence on $n$ will be suppressed. We also assume that $2m = n$ in this proof, which does not affect the asymptotic result. 

%Thus the result holds for a row-wise i.i.d.~triangular array. 

\vskip 1em

\noindent \textbf{$\bullet$ Verification of Lyapunov condition for $T_{\text{stand},h} $.}  Notice that $\overline{f}_h$ is the sample mean of i.i.d.~random variables conditional on $\mathcal{X}_2$. Therefore Lyapunov condition shows that if 
\begin{align} \label{Eq: convergence in prob}
	\frac{\mE_P[f_h^4(X_1)|\mathcal{X}_2]}{\{ \mE_P[f_h^2(X_1) | \mathcal{X}_{2}] \}^2 m} \convP 0,
\end{align}
the standardized statistic $T_{\text{stand},h}$ converges to $N(0,1)$ conditional on $\mathcal{X}_2$. Furthermore, since the limiting distribution is continuous, Lemma 2.11 of \cite{van2000asymptotic} shows
\begin{align*}
	\sup_{-\infty < z < \infty} \big| P(T_{\text{stand},h}  \leq z |\mathcal{X}_2) - \Phi(z) \big| \convP 0.
\end{align*}
Additionally, since $X_n \convP 0$ implies $\mE_P[X_n] \rightarrow 0$ when $X_n$ is bounded, we also have
\begin{align*}
	\sup_{-\infty < z < \infty} \big| P(T_{\text{stand},h}  \leq z) - \Phi(z) \big| \rightarrow 0.
\end{align*}
Hence it is enough to show (\ref{Eq: convergence in prob}) to establish the unconditional central limit theorem for $T_{\text{stand},h}$. To prove that Lyapunov condition~(\ref{Eq: convergence in prob}) holds under the given conditions, let us consider the following decomposition:
\begin{align} \label{Eq: decomposition I}
	\frac{\mE_P[f_h^4(X_1)|\mathcal{X}_{2}]}{\{ \mE_P[f_h^2(X_1) | \mathcal{X}_{2}] \}^2 m}  = \underbrace{\frac{m \mE_P[f_h^4(X_1)|\mathcal{X}_{2}]}{\{ \mE_P [h^2(X_1,X_2)] \}^2}}_{(I)} \times \underbrace{\frac{\{ \mE_P [h^2(X_1,X_2)] \}^2}{\{ \mE_P[f_h^2(X_1) | \mathcal{X}_{2}] \}^2 m^2}}_{(II)}.
\end{align}
Note that the expectation of the first term $(I)$ is
\begin{align*}
	\mE_P[(I)]  ~=~ &  \frac{m\mE_P[ \mE_P[f_h^4(X_1)|\mathcal{X}_{2}]]}{\{ \mE_P [h^2(X_1,X_2)] \}^2}  \\[.5em]
	= ~ & \frac{m}{\{ \mE_P [h^2(X,Y)] \}^2} \Bigg\{ \frac{1}{m^3} \mE_P[h^4(X_1,X_2)] + \frac{3m(m-1)}{m^4} \mE_P[h^2(X_1,X_3)h^2(X_2,X_3)] \Bigg\}.
\end{align*}
Therefore, by Markov's inequality, the first term converges to zero in probability when
\begin{align} \label{Eq: Asymptotic normality condition I}
	\frac{\mE_P[h^4(X_1,X_2)]n^{-1} + \mE_P[h^2(X_1,X_3)h^2(X_2,X_3)] }{n\{ \mE_P [h^2(X_1,X_2)] \}^2} \rightarrow 0.
\end{align}
In fact this is exactly the condition~(\ref{Eq: normality condition 2}). Hence all we need to show is that the second term $(II)$ is stochastically bounded. In other words, we would like to show that for any $\epsilon > 0$, there exist positive constants $N_\epsilon, M_\epsilon$ such that $\mP\bigl((II) \geq M_\epsilon \bigr) \leq g(\epsilon)$ for all $n \geq N_\epsilon$, where $g: (0, \infty) \mapsto (0, \infty)$ is a surjective function. In what follows, we derive two tail bounds for $(II)$ and combine them to prove the claim.

\vskip 1em

\noindent \textbf{First bound.} Fix $\epsilon > 0$ and recall that $h(x,y)$ can be written as
\begin{align*}
	h(x,y) = \sum_{k=1}^\infty \lambda_k \phi_k(x) \phi_k(y).
\end{align*}
Therefore we may write 
\begin{align*}
	\frac{m \mE_P[f_h^2(X_1) | \mathcal{X}_{2}]}{\mE_P[h^2(X_1,X_2)]} = \sum_{k=1}^\infty \widetilde{\lambda}_k  \Bigg(  \bigg\{ \frac{1}{\sqrt{m}} \sum_{j=m+1}^n \phi_k(X_j) \bigg\} \Bigg)^2,
\end{align*}
where $\widetilde{\lambda}_k = \lambda_k^2 / \sum_{k'=1}^\infty \lambda_{k'}^2$. Then we notice that 
\begin{align*}
	(II)^{1/2} ~ = ~  & \frac{1}{\sum_{k=1}^\infty \widetilde{\lambda}_k \big\{ \frac{1}{\sqrt{m}} \sum_{j=m+1}^n \phi_k(X_j) \big\}^2} \\
	\leq ~ & \frac{1}{\widetilde{\lambda}_1  \big\{ \frac{1}{\sqrt{m}} \sum_{j=m+1}^n \phi_1(X_j)\big\}^2}.
\end{align*}
Using this inequality, it can be seen that for any $t_1 > 0$,
\begin{align*}
	P\bigl((II) \geq t_1 \bigr) ~=~& P\bigl(t_1^{-1/2} \geq (II)^{-1/2} \bigr) \\[.5em]
	\leq ~ & P \biggl(\widetilde{\lambda}_1^{-1} t_1^{-1/2} \geq
	 \bigg\{ \frac{1}{\sqrt{m}} \sum_{j=m+1}^n \phi_1(X_j)\bigg\}^2 \biggr). 
\end{align*}
Since $\sum_{j=m+1}^n \phi_1(X_j)$ is the sum of i.i.d.~random variables where each $\phi_1(X_j)$ has mean zero and variance one, an application of the Berry--Esseen bound with the finite fourth moment yields
\begin{align*}
	 P \biggl(\widetilde{\lambda}_1^{-1} t_1^{-1/2} \geq
	\bigg\{ \frac{1}{\sqrt{m}} \sum_{j=m+1}^n \phi_1(X_j)\bigg\}^2 \biggr) \leq P \bigl( \widetilde{\lambda}_1^{-1} t_1^{-1/2} \geq \xi^2 \bigr) + C \sqrt{\frac{\mE_P[\phi_1^4(X)]}{n}},
\end{align*}
where $\xi \sim N(0,1)$ and $C>0$ is some universal constant. To further upper bound the above terms, note that the eigendecomposition (\ref{Eq: Another representation of h}) together with Fubini--Tonelli's theorem yields
\begin{align*}
	\frac{\mE_P[h^2(X_1,X_3)h^2(X_2,X_3)]}{n\{\mE_P[h^2(X_1,X_2)]\}^2} ~=~ & \frac{\sum_{k=1}^\infty \lambda_k^4 \mE_P[\phi_k^4(X)] + \sum_{k=1,k'=1, k \neq k'}^\infty\lambda_k^2 \lambda_{k'}^2 \mE_P[\phi_k^2(X)\phi_{k'}^2(X)]}{n\{\sum_{k=1}^\infty \lambda_i^2 \}^2} \\[.5em]
	\geq~ &  \widetilde{\lambda}_1^2 \frac{\mE_P[\phi_1^4(X)]}{n}.
\end{align*}
Since the above term converges to zero under the condition~(\ref{Eq: normality condition 2}), we may assume that there exists $N_{1,\epsilon}$ such that 
\begin{align*}
	\widetilde{\lambda}_1^2 \frac{\mE_P[\phi_1^4(X)]}{n} \leq \epsilon ~\Longleftrightarrow~  \frac{\mE_P[\phi_1^4(X)]}{n} \leq \epsilon \widetilde{\lambda}_1^{-2}
\end{align*}
for all $n \geq N_{1,\epsilon}$. Additionally, observe that
\begin{align*}
	P \bigl( \widetilde{\lambda}_1^{-1} t_1^{-1/2} \geq \xi^2 \bigr) = 2 \int_0^{\widetilde{\lambda}_1^{-1/2} t_1^{-1/4}} \frac{1}{\sqrt{2\pi}} e^{-\frac{x^2}{2}} dx \leq \frac{2}{\sqrt{2\pi}}\widetilde{\lambda}_1^{-1/2} t_1^{-1/4}. 
\end{align*}
Therefore for all $n \geq N_{1,\epsilon}$ and $t_1>0$, it holds that 
\begin{align} \label{Eq: Bound 1}
	P\bigl((II) \geq t_1 \bigr) \leq  \frac{2}{\sqrt{2\pi}} \frac{1}{\widetilde{\lambda}_1^{1/2} t_1^{1/4}} + C \frac{\sqrt{\epsilon}}{\widetilde{\lambda}_1}.
\end{align}

\vskip 1em

\noindent \textbf{Second bound.} Next we derive another upper bound suited for the case with small $\widetilde{\lambda}_1$. Again fix $\epsilon > 0$ and write
\begin{align*}
	g(x,y) = \mE_P[h(x,X)h(y,X)].
\end{align*}
Then after some algebra, we obtain the identity
\begin{align*}
	\frac{m\mE_P[f_h^2(X_1) | \mathcal{X}_{2}]}{\mE_P[h^2(X_1,X_2)]} ~ = ~ & \frac{1}{m} \sum_{j=m+1}^n \frac{g(X_j,X_j)}{\mE_P[h^2(X_1,X_2)]} + \frac{1}{m} \sum_{m+1 \leq j_1 \neq j_2 \leq n} \frac{g(X_{j_1}, X_{j_2})}{\mE_P[h^2(X_1,X_2)]}.
\end{align*}
The mean and the variance of the above quantity satisfy
\begin{align} \nonumber
	& \mE_P \bigg\{  \frac{m\mE_P[f_h^2(X_1) | \mathcal{X}_{2}]}{\mE_P[h^2(X_1,X_2)]} \bigg\} ~=~ 1 \quad \text{and} \\[.5em] \nonumber
	& \mV_P \bigg\{ \frac{m\mE_P[f_h^2(X_1) | \mathcal{X}_{2}]}{\mE_P[h^2(X_1,X_2)]} \bigg\} ~ \lesssim ~   \frac{\mE_P[h^2(X_1,X_3)h^2(X_2,X_3)]}{m\{\mE_P[h^2(X_1,X_2)] \}^2} + \frac{\mE_P[ g^2(X_1,X_2)]}{\{\mE_P[h^2(X_1,X_2)] \}^2}. 
\end{align}
Under the condition~(\ref{Eq: normality condition 2}), the first term in the upper bound converges to zero. Thus we may assume that there exists $N_{2,\epsilon}$ such that the first term is smaller than $\epsilon$ for all $n \geq N_{2,\epsilon}$. Also note the second term satisfies 
\begin{align*}
	\frac{\mE_P[ g^2(X_1,X_2)]}{\{\mE_P[h^2(X_1,X_2)] \}^2} = \frac{\sum_{k=1}^\infty \lambda_k^4 }{\{ \sum_{k=1}^\infty \lambda_k^2 \}^2} \leq \frac{\lambda_1^2}{\sum_{k=1}^\infty \lambda_k^2} = \widetilde{\lambda}_1.
\end{align*}
Having these observations in mind, Chebyshev's inequality with $t_2 \in (2,\infty)$ yields
\begin{align} \nonumber
	P\bigl((II) \geq t_2 \bigr) ~=~ & P\bigl(t_2^{-1/2} \geq (II)^{-1/2} \bigr) = P\bigl(t_2^{-1/2} - 1 \geq (II)^{-1/2} - 1 \bigr) \\[.5em] \nonumber
	 \leq ~ & \frac{1}{(1 - 1/\sqrt{t_2})^2}\mV_P \bigg\{ \frac{m\mE_P[f_h^2(X_1) | \mathcal{X}_{2}]}{\mE_P[h^2(X_1,X_2)]} \bigg\} \\[.5em]
	 \leq ~ & C' \{\epsilon + \widetilde{\lambda}_1\}, \label{Eq: Bound 2}
\end{align}
which holds for any $n \geq N_{2,\epsilon}$ and some positive constant $C'$.

\vskip 1em

\noindent \textbf{Combining two bounds. } Combining bounds~\eqref{Eq: Bound 1} and ~\eqref{Eq: Bound 2} yields that for any $n \geq \max \{N_{1,\epsilon}, N_{2,\epsilon}\}$ and $t > \max\{2,1/\epsilon^2\}$, 
\begin{align*}
	P\bigl((II) \geq t \bigr)  ~\leq~ & \min \bigg\{\frac{2}{\sqrt{2\pi}} \frac{1}{\widetilde{\lambda}_1^{1/2} t^{1/4}} + C \frac{\sqrt{\epsilon}}{\widetilde{\lambda}_1}, \, C'(\epsilon + \widetilde{\lambda}_1) \bigg\} \\[.5em]
	\leq ~ & \min \bigg\{ C'' \sqrt{\epsilon} \biggl( \frac{1}{\widetilde{\lambda}_1} + \frac{1}{\widetilde{\lambda}_1^{1/2}} \biggr), \, C'(\epsilon + \widetilde{\lambda}_1) \bigg\} 
\end{align*}
where $C'' >0$ is some constant. Suppose that $\widetilde{\lambda}_1 \leq \epsilon^{1/4}$. Then $P\bigl((II) \geq t \bigr) \leq C'(\epsilon + \epsilon^{1/4})$. On the other hand, if $\widetilde{\lambda}_1 > \epsilon^{1/4}$, then $P\bigl((II) \geq t \bigr) \leq C''(\epsilon^{1/4} + \epsilon^{3/8})$. Therefore, by taking $M_\epsilon > \max\{2, 1/\epsilon^2\}$ and $N_\epsilon \geq \max \{N_{1,\epsilon}, N_{2,\epsilon}\}$, it holds that for any $n \geq N_\epsilon$
\begin{align*}
	P\bigl((II) \geq M_\epsilon \bigr) \leq g(\epsilon) := 2(C' + C'') \max \{\epsilon, \epsilon^{1/4}\}.
\end{align*}
Since such $g$ is surjective, this proves that $(II)$ is stochastically bounded, which in turn verifies the claim~(\ref{Eq: convergence in prob}).

\vskip 1em

\noindent \textbf{$\bullet$ Central Limit Theorem for $T_{h}$.} So far we have established that $T_{\text{stand},h}$ converges to $N(0,1)$. Therefore, once we prove 
\begin{align} \label{Eq: convergence in prob to one}
	\frac{m^{-1} \sum_{i=1}^m \{ f_h(X_i) - \overline{f}_h\}^2}{\mE_P[f_h^2(X_1) | \mathcal{X}_2]} \convP 1,
\end{align}
Slutsky's theorem together with the continuous mapping theorem implies that $T_{h}$ converges to $N(0,1)$ as well. To show this, let us consider the decomposition
\begin{align*}
	&\frac{(m-1)^{-1} \sum_{i=1}^m \{ f_h(X_i) - \overline{f}_h\}^2 - \mE_P[f_h^2(X_1) | \mathcal{X}_2]}{\mE_P[f_h^2(X_1) | \mathcal{X}_2]}  \\[.5em]
	= ~ &  \underbrace{\frac{\sum_{i=1}^m \{ f_h(X_i) - \overline{f}_h\}^2  - (m-1) \mE_P[f_h^2(X_1) | \mathcal{X}_2]}{\mE_P[h^2(X_1,X_2)]}}_{(I)'} \times \underbrace{\frac{\mE_P[h^2(X_1,X_2)]}{(m-1) \mE_P[f_h^2(X_1) | \mathcal{X}_2]}}_{(II)'}.
\end{align*}
From the previous calculations for verifying~(\ref{Eq: convergence in prob}), we already know that $(II)' = O_P(1)$. Therefore all we need to show is $(I)' \convP 0$. Note that the expectation of $(I)'$ is equal to zero and the variance of $(I)'$ can be bounded by
\begin{align*}
	\mV_P[(I)'] ~=~ &  \frac{(m-1)^2}{\{\mE_P[h^2(X_1,X_2)] \}^2} \mE \Bigg\{ \mV_P \bigg[ \frac{1}{m-1} \sum_{i=1}^m \{ f_h(X_i) - \overline{f}_h\}^2 \Big| \mathcal{X}_2 \bigg] \Bigg\} \\[.5em]
	\leq ~ & \frac{m\mE_P[ f_h^4(X_1)]}{\{\mE_P[h^2(X_1,X_2)] \}^2},
\end{align*}
where the inequality uses the variance formula for the sample variance \citep[e.g.][]{lee1990u}. The upper bound is equivalent to the expectation of $(I)$ in (\ref{Eq: decomposition I}), which converges to zero as we proved before. As a result, we see that $(I)' \times (II)' \convP 0$ that further implies our claim~(\ref{Eq: convergence in prob to one}) since $(m-1)/m \rightarrow 1$. 

\vskip 1em

\noindent \textbf{$\bullet$ Concluding the proof.}
Following our previous argument, which holds for a row-wise i.i.d.~triangular array, we have 
\begin{align} \label{Eq: pointwise KS distance}
	\sup_{-\infty < z < \infty} \big| P_n(T_{h}  \leq z) - \Phi(z) \big| \rightarrow 0,
\end{align} 
for any sequence $P_n$ in $\mathcal{P}_{n,h}$. To finish the proof, for each $n$, choose $P_n' \in \mathcal{P}_{n,h}$ such that 
\begin{align*}
	\sup_{P \in \mathcal{P}_{n,h}} \sup_{-\infty < z < \infty} \big| P(T_{h}  \leq  z) - \Phi(z) \big| \leq \sup_{-\infty < z < \infty} \big| P_{n'}(T_{h}  \leq z) - \Phi(z) \big| + n^{-1}.
\end{align*}
Now it follows that the right-hand side converges to zero based on the result~(\ref{Eq: pointwise KS distance}) and thus we arrive at the conclusion of Theorem~\ref{Theorem: Gaussian Approximation}.

\subsection{Proof of Proposition~\ref{Proposition: MMD test}}

We again assume that $2m=n$ for simplicity, which does not affect our asymptotic results. We start with type I error control and then move to type II error control. For this proof, we often need to take the expectation involving random variables from different distributions and so we suppress dependence on $P$ in the expectation to simplify the notation. We also let $h_{\text{Gau},\nu_n} := h$ for simplicity. 

\vskip 1em

\noindent \textbf{$\bullet$ Type I error control.} We use Theorem~\ref{Theorem: Gaussian Approximation} to prove type I error control. More specifically we verify that the kernel~(\ref{Eq: Gaussian kernel}) satisfies $\mE[h(X,y)]=0$, the condition~(\ref{Eq: normality condition 2}) and $\lambda_1^2 / \sum_{k=1}^\infty \lambda_k^2 \rightarrow 0$ under the given assumptions. The first condition is trivially true by the definition of the kernel. The latter two conditions were in fact verified in \cite{li2019optimality} (see the proof of their theorem 1) by noting that the condition~(\ref{Eq: normality condition 1}) is equivalent to $\lambda_1^2 / \sum_{k=1}^\infty \lambda_k^2 \rightarrow 0$ since 
\begin{align*}
	\Bigg[ \frac{\lambda_1^2}{\sum_{k=1}^\infty \lambda_k^2} \Bigg]^2 \leq \frac{\sum_{k=1}^\infty \lambda_k^4}{\big\{ \sum_{k=1}^\infty \lambda_k^2 \big\}^2} \leq \frac{\lambda_1^2}{\sum_{k=1}^\infty \lambda_k^2}.
\end{align*}
Hence the test statistic converges to the standard normal distribution and the result follows.

\vskip 1em

\noindent \textbf{$\bullet$ Type II error control.} To prove type II error control, we recall that the type II error is
\begin{align*}
	P( T_{\text{MMD}} \leq z_{1-\alpha})  = P \bigg\{ \overline{f}_h \leq z_{1-\alpha} \bigg[m^{-2} \sum_{i=1}^m \{ f_h(X_i) - \overline{f}_h \}^2 \bigg]^{1/2} \bigg\}, 
\end{align*}
where $f_h$ is defined in (\ref{Eq: general linear function}) with $h(x,y) =h_{\text{Gau},\nu_n}(x,y)$ and $\overline{f}_h = m^{-1} \sum_{i=1}^m f_h(X_i)$. Throughout this proof, we use the notation $a_n \lesssim b_n$ for some positive sequences $a_n,b_n$ if there exists some constant $C>0$ independent of $n$ such that $a_n \leq C b_n$. We proceed the proof by checking the sufficient condition in Lemma~\ref{Lemma: Sufficient conditions for uniform type II error control}. In particular, by looking at sufficient condition~(\ref{Eq: Sufficient Condition for Power}), we only need to deal with 
\begin{align*}
	\mE[h(X_1,X_2)], \ \ \mE[\mV[h(X_1,X_2)|X_1]] \  \text{and} \  \mV[\mE[h(X_1,X_2)|X_1]].
\end{align*}
To start with the mean, note that the expectation of $\overline{f}_h$ is the same as the expectation of the U-statistic in \cite{li2019optimality}. Therefore we can use the result of \cite{li2019optimality} and obtain
\begin{align} \label{Eq: lower bound for the mean}
	\mE[\overline{f}_h] ~\gtrsim~ \nu_n^{-d/2} \| p - q \|_{L_2}^2 \quad \text{for sufficiently large $n$.}
\end{align}
Next, the conditional expectation is bounded as
\begin{align} \nonumber
	\mE[\mV[h(X_1,X_2)|X_1]] ~=~ & \mE[\{ \mE[G_{\nu_n}(X_1,Y)|X_1] - \mE[G_{\nu_n}(X,Y)] \}^2] \\[.5em] \nonumber
	\overset{\mathrm{(i)}}{\leq} ~ & 2 \mE[G_{\nu_n}^2(X,Y)] \\[.5em] \nonumber
	\overset{\mathrm{(ii)}}{\leq} ~ &  \mE[G_{\nu_n}^2(X,X')]  +  \mE[G_{\nu_n}^2(Y,Y')] \\[.5em]  \label{Eq: bound on condi. exp.}
	\overset{\mathrm{(iii)}}{\lesssim}  ~ & \frac{M^2}{\nu^{d/2}}, 
\end{align}
where step~(i) uses Jensen's inequality, step~(ii) follows by using the inequality that $p \times q \leq p^2/2 + q^2/2$, and the last step uses the result of \cite{li2019optimality}. 

On the other hand, the conditional variance is bounded as
\begin{align} \nonumber
	& \mV[\mE[h(X_1,X_2)|X_1]] \\[.5em] \nonumber
	=~ & \mV[\{ \mE[G_{\nu_n}(X_1,X_2)|X_1] - \mE[G_{\nu_n}(X_1,Y)|X_1] + \mE[G_{\nu_n}(X,Y)] - \mE[G_{\nu_n}(X,X')] \}] \\[.5em] \nonumber
	\leq ~ & \mE[\{ \mE[G_{\nu_n}(X_1,X_2)|X_1] - \mE[G_{\nu_n}(X_1,Y)|X_1] \}^2 ] \\[.5em] \label{Eq: bound on condi. var.}
	\lesssim ~ & \frac{M}{\nu^{3d/4}} \| p - q \|_{L_2}^2,
\end{align}
where the last inequality follows by using the result in \cite{li2019optimality}. Since we assume $\alpha$ is fixed and $m_1=m_2$, the bounds in (\ref{Eq: lower bound for the mean}), (\ref{Eq: bound on condi. exp.}) and (\ref{Eq: bound on condi. var.}) guarantee that if 
\begin{align*}
	 \nu_n^{-d/2} \| p - q \|_{L_2}^2 \gtrsim  \psi_n \max\biggl\{ \frac{M}{n \nu^{d/4}}, \ \frac{M^{1/2}}{n^{1/2}\nu^{3d/8}} \| p - q \|_{L_2} \bigg\} \quad \text{for sufficiently large $n$, }
\end{align*}
then, sufficient condition~(\ref{Eq: Sufficient Condition for Power}) is satisfied. Now recall that we choose $\nu_n \asymp n^{4/(d+4s)}$. Therefore, the above condition is fulfilled when $\|p-q\|_{L_2} \gtrsim \max\{\psi_n, \sqrt{\psi}_n \} \nu_n^{d/8} / \sqrt{n}$. Finally, we choose $\psi_n$ such that $\|p-q\|_{L_2} \geq K_n n^{-\frac{2s}{4s+d}} \gtrsim \max\{\psi_n, \sqrt{\psi}_n \} \nu_n^{d/8} / \sqrt{n}$, which proves the desired claim. 

\subsection{Proof of Proposition~\ref{Proposition: limiting distribution of cross-splitting}}
By the multivariate central limit theorem, we have 
\begin{align*}
	\sqrt{m} \begin{pmatrix}
		\widehat{\mu}_1 \\
		\widehat{\mu}_2
	\end{pmatrix} \convD N \left( \begin{bmatrix}
		0 \\
		0
	\end{bmatrix}, \begin{bmatrix}
		\Sigma & 0 \\
		0 & \Sigma
	\end{bmatrix}
	\right).
\end{align*}
For nonzero $x, y \in \mathbb{R}^d$, let us define $g(x,y) := x^\top y / \sqrt{x^\top \Sigma x} + x^\top y / \sqrt{y^\top \Sigma y}$. Then the continuous mapping theorem yields
\begin{align} \label{Eq: continuous mapping theorem}
	\frac{m \widehat{\mu}_1^\top \widehat{\mu}_2}{\sqrt{ m \widehat{\mu}_1^\top \Sigma \widehat{\mu}_1 }} + \frac{m \widehat{\mu}_1^\top \widehat{\mu}_2}{\sqrt{ m \widehat{\mu}_2^\top \Sigma \widehat{\mu}_2}} \convD \frac{\boldsymbol{\xi}_1^\top \Sigma \boldsymbol{\xi}_2}{\sqrt{\boldsymbol{\xi}_1^\top \Sigma^2 \boldsymbol{\xi}_1}} + \frac{\boldsymbol{\xi}_1^\top \Sigma \boldsymbol{\xi}_2}{\sqrt{\boldsymbol{\xi}_2^\top \Sigma^2 \boldsymbol{\xi}_2}}.
\end{align}
Next let $S_1 := m^{-1} \sum_{i=1}^m (X_i - \widehat{\mu}_1)(X_i - \widehat{\mu}_1)^\top$ and $S_2 := m^{-1} \sum_{i=m+1}^n (X_i - \widehat{\mu}_2)(X_i - \widehat{\mu}_2)^\top$. Since the dimension is fixed, the law of large numbers shows that $\|S_i - \Sigma \|_{\text{op}} = o_P(1)$ for $i=1,2$ where $\| A \|_{\text{op}}$ denotes the operator norm of a matrix $A$. Since we can write 
\begin{align*}
	T_{\mathrm{mean}}^{\text{cross}} = \frac{m \widehat{\mu}_1^\top \widehat{\mu}_2}{\sqrt{ m \widehat{\mu}_1^\top (S_2 - \Sigma) \widehat{\mu}_1 + m \widehat{\mu}_1^\top \Sigma \widehat{\mu}_1}} +\frac{m \widehat{\mu}_1^\top \widehat{\mu}_2}{\sqrt{ m \widehat{\mu}_2^\top (S_1 - \Sigma) \widehat{\mu}_2+ m \widehat{\mu}_2^\top \Sigma \widehat{\mu}_2}}
\end{align*}
and $|m \widehat{\mu}_i^\top (S_j - \Sigma) \widehat{\mu}_i| \leq m \widehat{\mu}_i^\top \widehat{\mu}_i \times \|S_j - \Sigma \|_{\text{op}} =  o_P(1)$ for $i \neq j \in \{1,2\}$, 
\begin{align*}
	T_{\mathrm{mean}}^{\text{cross}} ~=~ &  \frac{m \widehat{\mu}_1^\top \widehat{\mu}_2}{\sqrt{m \widehat{\mu}_1^\top \Sigma \widehat{\mu}_1 + o_P(1)}} + \frac{m \widehat{\mu}_1^\top \widehat{\mu}_2}{\sqrt{m \widehat{\mu}_2^\top \Sigma \widehat{\mu}_2 + o_P(1)}} \\[.5em]
	~=~ & \frac{m \widehat{\mu}_1^\top \widehat{\mu}_2}{\sqrt{ m \widehat{\mu}_1^\top \Sigma \widehat{\mu}_1 }} + \frac{m \widehat{\mu}_1^\top \widehat{\mu}_2}{\sqrt{ m \widehat{\mu}_2^\top \Sigma \widehat{\mu}_2}} + o_P(1),
\end{align*}
by the Taylor expansion. Finally Slutsky's theorem along with the convergence~(\ref{Eq: continuous mapping theorem}) proves the claim. 

\vskip 2em

\subsection{Proof of Lemma~\ref{Lemma: Sufficient conditions for uniform type II error control}} \label{Section: Proof of Lemma: Sufficient conditions for uniform type II error control}
We start by collecting preliminary results and then combine these thorough concentration bounds to prove the claim. 

\vskip1em
\noindent \textbf{Preliminaries.} Let us first consider sample mean of $f_h(X_1),\ldots,f_h(X_{m_1})$, i.e.~$\overline{f}_h$, and compute its mean and variance. Due to the independence between $\mathcal{X}_1=\{X_i\}_{i=1}^{m_1}$ and $\mathcal{X}_2=\{X_i\}_{i=m_1+1}^{n}$, the mean of $\overline{f}_h$ is trivially $\mE_P[h(X_1,X_2)]$. On the other hand, by the law of total variance, the variance of $\overline{f}_h$ is 
\begin{align*}
	\mV_P[\overline{f}_h] ~=~ \underbrace{\mE_P[\mV_P[\overline{f}_h|\mathcal{X}_2]]}_{:=(I)} + \underbrace{\mV_P[\mE_P[\overline{f}_h|\mathcal{X}_2]]}_{:=(II)}.
\end{align*}
Note that due to the conditional independence of $f_h(X_1),\ldots,f_h(X_{m_1})$, the first term becomes
\begin{align*}
	(I) ~=~ &  \mE_P\biggl[ \mV_P \biggl[ \frac{1}{m_1} \sum_{i=1}^{m_1} f_h(X_i) | \mathcal{X}_2 \biggr] \biggr] = \frac{1}{m_1} \mE_P \bigl[ \mV_P\bigl[ f_h(X_1) | \mathcal{X}_2 \bigr] \bigr] \\[.5em]
	= ~ & \frac{1}{m_1} \mE_P \biggl[ \mV_P\biggl[ \frac{1}{m_2} \sum_{i=m_1+1}^n h(X_i,X_1) | \mathcal{X}_2 \biggr] \biggr] \\[.5em]
	= ~ & \frac{1}{m_1m_2} \mE_P[\mV_P[h(X_1,X_2)|X_1]] \\[.5em]
	& + \frac{1}{m_1} \mE_P\biggl[ \frac{m_2(m_2-1)}{m_2^2} \mathrm{cov}_P\bigl\{ h(X_1,X_2), h(X_1,X_3) \big| X_2,X_3 \bigr\} \biggr] \\[.5em]
	\leq ~ & \frac{1}{m_1m_2} \mE_P[\mV_P[h(X_1,X_2)|X_1]] + \frac{1}{m_1} \mE_P\bigl[  \mathrm{cov}_P\bigl\{ h(X_1,X_2), h(X_1,X_3) \big| X_2,X_3 \bigr\} \bigr].
\end{align*}
Also note that 
\begin{align*}
	& \mE_P\bigl[  \mathrm{cov}_P \bigl\{ h(X_1,X_2), h(X_1,X_3) \big| X_2,X_3 \bigr\} \bigr] \\[.5em]
	= ~ & \mE_P\bigl[ \{h(X_1,X_2) - \mE_P[h(X_1,X_2)|X_2]\} \{h(X_1,X_3) - \mE_P[h(X_1,X_3)|X_3]\} \bigr] \\[.5em]
	= ~ & \mE_P[\{ \mE_P[h(X_1,X_2)|X_1]\}^2] -  \{\mE_P[\mE_P[h(X_1,X_2)|X_1]]\}^2 \\[.5em]
	=  ~ &  \mV_P[\mE_P[h(X_1,X_2)|X_1]]. 
\end{align*}
Therefore the first term is bounded by
\begin{align*}
	(I) ~ \leq ~  \frac{1}{m_1m_2} \mE_P[\mV_P[h(X_1,X_2)|X_1]] + \frac{1}{m_1}\mV_P[\mE_P[h(X_1,X_2)|X_1]].
\end{align*}
Similarly, the second term is 
\begin{align*}
	(II) ~=~ & \mV_P[\mE_P[f_h(X_1)|\mathcal{X}_2]] \\[.5em]
	= ~ & \mV_P\biggl[ \frac{1}{m_2} \sum_{i=m_1+1}^{n} \mE_P[f(X_i,X_1)|X_i] \biggr] \\[.5em]
	= ~ & \frac{1}{m_2} \mV_P[\mE_P[h(X_1,X_2)|X_1]]. 
\end{align*}
Combining the above bounds yields
\begin{align} \label{Eq: Variance Upper Bound}
	\mV_P[\overline{f}_h] ~ \leq ~ \frac{1}{m_1m_2} \mE_P[\mV_P[h(X_1,X_2)|X_1]] + \biggl( \frac{1}{m_1} + \frac{1}{m_2} \biggr)  \mV_P[\mE_P[h(X_1,X_2)|X_1]].
\end{align}
Next we compute the expected value of the sample standard deviation:
\begin{align*}
	& \mE_P\biggl[  \frac{1}{m_1^2} \sum_{i=1}^{m_1} \{ f_h(X_i) - \overline{f}_h \}^2 \biggr] ~=~  \frac{1}{m_1} \mE_P\bigl[ \{ f_h(X_1) - \overline{f}_h \}^2  \bigr] \\[.5em]
	= ~ & \frac{1}{m_1} \mE_P\biggl[ \bigg\{ \biggl( 1 - \frac{1}{m_1} \biggr) f_h(X_1) - \frac{1}{m_1} \sum_{i=2}^{m_1} f_h(X_i) \bigg\}^2\biggr] \\[.5em]
	= ~ & \frac{(m_1-1)^2}{m_1^3} \mE_P [\mV_P[f_h(X_1|\mathcal{X}_2)]] + \frac{m_1-1}{m_1^3}  \mE_P [\mV_P[f_h(X_1|\mathcal{X}_2)]] \\[.5em]
	= ~ & \frac{m_1-1}{m_1^2} \mE_P [\mV_P[f_h(X_1|\mathcal{X}_2)]].
\end{align*}
Using the previous calculation, we then obtain
\begin{align*}
	& \mE_P\biggl[ \frac{1}{m_1^2} \sum_{i=1}^{m_1} \{ f_h(X_i) - \overline{f}_h \}^2 \biggr] \\[.5em]
	~ \leq ~ & \frac{1}{m_1m_2} \mE_P[\mV_P[h(X_1,X_2)|X_1]] + \frac{1}{m_1}\mV_P[\mE_P[h(X_1,X_2)|X_1]].
\end{align*}

\vskip 1em

\noindent \textbf{Combining pieces.} For simplicity, let us write
\begin{align*}
	& \mu_P:= \mE_P[h(X_1,X_2)], \ \tau_P:=\mE_P[\mV_P[h(X_1,X_2)|X_1]] \quad \text{and} \\[.5em] & \sigma_P:=\mV_P[\mE_P[h(X_1,X_2)|X_1]].
\end{align*}
Note that Markov's inequality ensures that the following event $\mathcal{A}$:
\begin{align}
	\mathcal{A} = \Bigg\{  \frac{1}{m_1^2} \sum_{i=1}^{m_1} \{ f_h(X_i) - \overline{f}_h \}^2 < \frac{2}{\beta} \bigg[ \frac{\tau_P}{m_1m_2} + \sigma_P \biggl( \frac{1}{m_1} +\frac{1}{m_2} \biggr) \bigg]  \Bigg\}
\end{align}
holds with probability at least $\beta/2$ for some $\beta > 0$ specified later. Then, by the union bound argument, the type II error can be bounded as
\begin{align*}
	& P \Biggl( \overline{f}_h \leq  z_{1-\alpha} \sqrt{ \frac{1}{m_1^2} \sum_{i=1}^{m_1} \{ f_h(X_i) - \overline{f}_h \}^2} \Biggr) \\[.5em]
	\leq ~ & P \Biggl( \overline{f}_h \leq  \max\{z_{1-\alpha},1\} \sqrt{ \frac{1}{m_1^2} \sum_{i=1}^{m_1} \{ f_h(X_i) - \overline{f}_h \}^2} \Biggr) \\[.5em]
	\leq ~ & P \Biggl( \overline{f}_h \leq  \max\{z_{1-\alpha},1\} \sqrt{ \frac{1}{m_1^2} \sum_{i=1}^{m_1} \{ f_h(X_i) - \overline{f}_h \}^2}, \ \mathcal{A} \Biggr) + P(\mathcal{A}^c) \\[.5em]
	~ \leq  ~ & P \biggl( \overline{f}_h \leq  \max\{z_{1-\alpha},1\} \sqrt{\frac{2}{\beta} \bigg[ \frac{\tau_P}{m_1m_2} + \sigma_P \biggl( \frac{1}{m_1} +\frac{1}{m_2} \biggr) \bigg]} \biggr) + \frac{\beta}{2}.
\end{align*}
Recall the condition of the lemma: 
\begin{align} \label{Eq: Lower Bound on Mean}
	\mu_P \geq \psi_n \max\{z_{1-\alpha},1\}\max \bigg\{ \sqrt{\frac{\tau_P}{m_1m_2}}, \ \sqrt{ \sigma_P \biggl( \frac{1}{m_1} +\frac{1}{m_2} \biggr)}  \bigg\},
\end{align}
where $\psi_n$ is a positive sequence that goes to infinity as $\min\{m_1,m_2\} \rightarrow \infty$. Therefore, by letting $\beta = 1/\psi_n$, we can take $N \in \mathbb{N}$, independent of $P \in \mathcal{Q}_1$, such that for all $\min\{m_1,m_2\} \geq N$, it holds that
\begin{align*}
	& \psi_n \max\{z_{1-\alpha},1\} \max \bigg\{ \sqrt{\frac{\tau_P}{m_1m_2}}, \ \sqrt{ \sigma_P \biggl( \frac{1}{m_1} +\frac{1}{m_2} \biggr)}  \bigg\} \\[.5em] 
	\geq ~ & 2  \max\{z_{1-\alpha},1\} \sqrt{\frac{2}{\beta} \bigg[ \frac{\tau_P}{m_1m_2} +  \sigma_P \biggl( \frac{1}{m_1} +\frac{1}{m_2} \biggr) \bigg]}.
\end{align*}  
Under this event, the first term in the upper bound of the type II error satisfies
\begin{align*}
	& P \biggl( \overline{f}_h \leq  \max\{z_{1-\alpha},1\} \sqrt{\frac{2}{\beta} \bigg[ \frac{\tau_P}{m_1m_2} +  \sigma_P \biggl( \frac{1}{m_1} +\frac{1}{m_2} \biggr) \bigg]} \biggr) \\[.5em]
	~ = ~ & P  \biggl( \overline{f}_h - \mu_P \leq  \max\{z_{1-\alpha},1\}\sqrt{\frac{2}{\beta} \bigg[ \frac{\tau_P}{m_1m_2} +  \sigma_P \biggl( \frac{1}{m_1} +\frac{1}{m_2} \biggr) \bigg]} - \mu_P \biggr) \\[.5em]
	\leq ~ & P \biggl(  \overline{f}_h - \mu_P \leq  -\frac{\mu_P}{2} \biggr) ~\leq ~ \frac{4\mV_P[\overline{f}_h]}{\mu_P^2},
\end{align*}
where the last step uses Chebyshev's inequality. Now using the upper bound on the variance of $\overline{f}_h$ in (\ref{Eq: Variance Upper Bound}) and also the lower bound on $\mu_P$ in (\ref{Eq: Lower Bound on Mean}), the type II error is less than or equal to 
\begin{align*}
	P \Biggl( \overline{f}_h \leq  z_{1-\alpha} \sqrt{ \frac{1}{m_1^2} \sum_{i=1}^{m_1} \{ f_h(X_i) - \overline{f}_h \}^2} \Biggr)  \leq \frac{8}{\psi_n^2 \max\{z_{1-\alpha}^2,1\}} + \frac{1}{2 \psi_n}.
\end{align*}
Notice that the right-hand side does not depend on $P \in \mathcal{Q}_1$ and $\psi_n \rightarrow \infty$. Therefore, it holds 
\begin{align*}
	\lim_{\min\{m_1,m_2\} \rightarrow \infty }\sup_{P \in \mathcal{Q}_1} P \Biggl( \overline{f}_h \leq  z_{1-\alpha} \sqrt{ \frac{1}{m_1^2} \sum_{i=1}^{m_1} \{ f_h(X_i) - \overline{f}_h \}^2} \Biggr) = 0,
\end{align*}
as desired. This completes the proof of Lemma~\ref{Lemma: Sufficient conditions for uniform type II error control}.

\end{document}